\documentclass[draft,preprint,12pt,numbers,sort&compress]{elsarticle}
\usepackage{}
\usepackage{mathrsfs}
\usepackage{amssymb}
\usepackage{amsthm}
\usepackage{mathrsfs}
\usepackage[centertags]{amsmath}
\usepackage{amsfonts}

\usepackage{color}

\input amssym.def
\input amssym.tex

\date{}

\newtheorem{Theorem}{Theorem}[section]

\newtheorem{Lemma}{Lemma}[section]

\newcommand\R{\mbox{\bf R}}

\newcommand\SR{\mbox{\scriptsize\bf R}}

\newcommand{\definition}{{\lower .5ex
  \hbox{$\>\>\stackrel{\triangle}{=}\>\>$} }}

\begin{document}

\baselineskip=22pt
\thispagestyle{empty}

\mbox{}
\bigskip

\begin{center}
{\Large \bf Global well-posedness of the Cauchy problem}\\[0.5ex]
{\Large \bf for a fifth-order KP-I equation in }\\[0.5ex]
{\Large \bf anisotropic Sobolev spaces}

\bigskip
\bigskip

{Yongsheng Li\footnote{Emails: yshli@scut.edu.cn (YS Li);
yanwei19821115@sina.cn (W Yan);
zhangyimin@whut.edu.cn (YM Zhang)}$^a$,\quad
Wei Yan$^{b}$,\quad
Yimin Zhang$^c$}\\[1ex]

{$^a$School of Mathematics, South  China  University of  Technology,}\\
{Guangzhou, Guangdong 510640,  P. R. China}\\[1ex]

{$^b$School of Mathematics and Information Science,}\\
{Henan Engineering Laboratory for
 Big Data Statistical Analysis and Optimal Control,}\\
{Henan  Normal University, Xinxiang, Henan 453007, P. R. China}\\[2ex]

{$^c$Department of Mathematics, Wuhan University of  Technology, Wuhan430070, China}\\
{Wuhan, Hubei 430070,  P. R. China}\\[1ex]

\end{center}

\bigskip
\bigskip

{\bf Abstract}
In this paper, we consider the Cauchy problem for
the fifth-order KP-I equation
\begin{align*}
   u_{t}+\partial_{x}^{5}u+\partial_{x}^{-1}\partial_{y}^{2}u
    +\frac{1}{2}\partial_{x}(u^{2})=0.
\end{align*}
Firstly, we establish the  local well-posedness of the problem   in the anisotropic
Sobolev spaces $H^{s_{1},\>s_{2}}(\R^{2})$
  with $s_{1}>-\frac{9}{8}$ and $s_{2}\geq 0$. Secondly,  we establish
the  global well-posedness of  the problem in $H^{s_{1},\>0}(\R^{2})$
with $s_{1}>-\frac{4}{7}$. Our result improves considerably the results of
 Saut and  Tzvetkov (J. Math.\ Pures Appl.\ 79(2000), 307--338.) and Li and Xiao
 (J. Math.\ Pures Appl.\ 90(2008), 338--352.) and  Guo, Huo and Fang
 (J. Diff.\ Eqns.\  263 (2017), 5696--5726).

\pagebreak

\bigskip


\bigskip
\noindent {\bf Short Title:} Cauchy problem for fifth-order  KP-I equation

\bigskip
\noindent {\bf Corresponding Author:} Wei Yan

\bigskip
\noindent {\bf Email Address:} yanwei19821115@sina.cn

\bigskip

\leftskip 0 true cm \rightskip 0 true cm

\newpage{}

\begin{center}

{\Large \bf Global well-posedness of the Cauchy problem}\\[0.5ex]
{\Large \bf for a fifth-order KP-I equation in }\\[0.5ex]
{\Large \bf anisotropic Sobolev spaces}

\bigskip
\bigskip

{Yongsheng Li\footnote{Emails: yshli@scut.edu.cn (YS Li);
yanwei19821115@sina.cn (W Yan);
zhangyimin@whut.edu.cn (YM Zhang)}$^a$,\quad
Wei Yan$^{b}$,\quad
Yimin Zhang$^c$}\\[1ex]

{$^a$School of Mathematics, South  China  University of  Technology,}\\
{Guangzhou, Guangdong 510640,  P. R. China}\\[1ex]

{$^b$School of Mathematics and Information Science,}\\
{Henan Engineering Laboratory for
 Big Data Statistical Analysis and Optimal Control,}\\
{Henan  Normal University, Xinxiang, Henan 453007, P. R. China}\\[2ex]
{$^c$Department of Mathematics, Wuhan University of  Technology, Wuhan430070, China}\\
{Wuhan, Hubei 430070,  P. R. China}
\end{center}

{\bf Abstract}

In this paper, we consider the Cauchy problem  for the
fifth-order KP-I  equation
\begin{align*}
 u_{t}+\partial_{x}^{5}u+\partial_{x}^{-1}\partial_{y}^{2}u+\frac{1}{2}\partial_{x}(u^{2})=0.
 \end{align*}
Firstly, we establish the  local well-posedness of the problem
in the anisotropic Sobolev spaces $H^{s_{1},\>s_{2}}(\R^{2})$
  with $s_{1}>-\frac{9}{8}$ and $s_{2}\geq 0$. Secondly,  we establish
the  global well-posedness of  the problem in $H^{s_{1},\>0}(\R^{2})$ with $s_{1}>-\frac{4}{7}$.
Our result improves considerably the results of
 Saut and  Tzvetkov (J. Math.\ Pures Appl.\ 79(2000), 307--338.) and
 Li and Xiao (J. Math.\ Pures Appl.\ 90(2008), 338--352.) and  Guo, Huo and Fang
 (J. Diff.\\ Eqns.\  263 (2017), 5696--5726).

\bigskip

\bigskip

{\large\bf 1. Introduction}
\bigskip

\setcounter{Theorem}{0} \setcounter{Lemma}{0}

\setcounter{section}{1}

This paper is devoted to studying the Cauchy problem for the fifth-order KP-I equation
\begin{eqnarray}
&& u_{t}+\partial_{x}^{5}u+\partial_{x}^{-1}\partial_{y}^{2}u
+\frac{1}{2}\partial_{x}(u^{2})=0,\label{1.01}\\
&&u(x,y,0)=u_{0}(x,y)\label{1.02}
\end{eqnarray}
in anisotropic Sobolev space $H^{s_{1},s_{2}}(\R^{2})$.

(\ref{1.01}) appears as a model describing certain long dispersive
 waves (see \cite{AS,KB1991,KB-1991}).
It is considered as the higher-order version of the following KP equation
\begin{eqnarray}
 u_{t}+\alpha \partial_{x}^{3}u+\partial_{x}^{-1}\partial_{y}^{2}u
 +\frac{1}{2}\partial_{x}(u^{2})=0,\label{1.03}
\end{eqnarray}
where the coefficient $\alpha$ may be either positive or negative. The KP
equation (\ref{1.03}) occurs in physical contexts as models for the propagation
of dispersive long waves with weak
transverse effects and is regarded as the
two-dimensional extensions of the Korteweg-de-Vries
equation (see \cite{KP}).
When $\alpha<0$, (\ref{1.03}) is known as the KP-I equation.
When $\alpha>0$, (\ref{1.03}) is known as the KP-II equation.

Several people have studied its Cauchy problem for (\ref{1.03}), see \cite{Bourgain-GAFA-KP,Hadac2008,Hadac2009,IMCPDE,IMEJDE,IMT,ILM-CPDE,IM2011,KL,TDCDS,TADE,TT,TCPDE,TDIE,TIMRN}
for the KP-II equation  (\ref{1.03}) with $\alpha>0$, and see \cite{BS,CKS-GAFA,CIKS,GPW,
 HN,IMS-NA,IMS1992,IMS-JMAA,IMS,IKT,ILP,Kenig2004,LJDE,MST-Duke,MST2002,MST2004,MST2011,Z}
 for the  KP-I equation (\ref{1.03}) with $\alpha<0$.

 For the KP-II equation,
  by using  the  Fourier  restriction  norm method, Bourgain \cite{Bourgain-GAFA-KP} established
   the global well-posedness of  its Cauchy problem
 in $L^{2}(\R^{2})$ and $L^{2}(\mathbf{T}^{2}).$ Takaokao and  Tzvetkov \cite{TT} and
 Isaza and  Mej\'{\i}a \cite{IMCPDE} established the local well-posedness
  of KP-II equation in $H^{s_{1},s_{2}}(\R^{2})$
 with $s_{1}>-\frac13$ and $s_{2}\geq0.$ Takaoka \cite{TDCDS} established the local well-posedness
 of KP-II equation in $H^{s_{1},0}(\R^{2})$
 with $s_{1}>-\frac12$  with the assumption that
 \begin{eqnarray*}
 \left\||\xi|^{-\frac{1}{2}+\epsilon}\mathscr{F}_{x}u_{0}\right\|_{L^{2}}<\infty
 \end{eqnarray*}
 for  the suitable chosen $\epsilon$. By introducing some resolution spaces,
 Hadac  et al. \cite{Hadac2009}  established the small data global well-posedness
 and scattering result of  KP-II equation in the homogeneous anisotropic Sobolev space
$\dot{H}^{-\frac{1}{2},\>0}(\R^{2})$ defined in \cite{Hadac2009} and arbitrary
 large initial data local well-posedness
 in both homogeneous Sobolev space $\dot{H}^{-\frac{1}{2},\>0}(\R^{2})$ and
 inhomogeneous anisotropic  Sobolev space
$H^{-\frac{1}{2},\>0}(\R^{2})$. Recently, by using new bilinear estimates,
   Koch and Li \cite{KL} established the global well-posedness and scattering
    for the KP-II equation in three space dimensions with small initial data.

For the KP-I equation, Kenig, Molinet, Saut and Tzvetkov
studied its Cauchy problem and periodic boundary value problem and showed that
the problems are globally well-posed
in the  second energy spaces on both $\R^{2}$ and $\mathbf{T^{2}}$
(see \cite{Kenig2004,MST2002,MST2004}).
 Molinet et al. \cite{MST-Duke} proved that the Picard
  iterative method does not work for the KP-I equation in standard Sobolev space and in anisotropic
Sobolev space,  since the flow map fails to be real-analytic at the origin in these spaces.
Ionescu et al. \cite{IKT} established the  global  well-posedness of  KP-I in the natural energy space
$E^{1}$  with the aid of  some  resolution  spaces and  bootstrap  inequality and  the energy estimates.
Molinet et al. \cite{MST2007} established the local well-posedness of  the Cauchy problem for
the KP-I equation in $H^{s,\>0}(\R^{2})$ with $s>\frac{3}{2}.$
Guo et al. \cite{GPW} established the  local  well-posedness of the Cauchy problem for the KP-I equation
in $H^{1,\>0}(\R^{2}).$ Zhang \cite{Z} established the local well-posedness of the
  periodic  KP-I initial value problem
 in the Besov type space
$
B_{2,1}^{\frac{1}{2}}(\mathbf{T}^{2}).
$

 Saut and Tzvetkov \cite{ST1999} established the global well-posedness of
 the Cauchy problem for  the fifth order KP-II equation
 \begin{eqnarray}
&& u_{t}-\partial_{x}^{5}u+\alpha \partial_{x}^{3}u+\partial_{x}^{-1}\partial_{y}^{2}u
+\frac{1}{2}\partial_{x}(u^{2})=0,\alpha\in \R,\label{1.04}
\end{eqnarray}
 in $H^{s_{1},s_{2}}(\R^{2})$ with $s_{1}>-\frac14,s_{2}\geq0$. Isaza et al. \cite{ILM-CPAA}
 established  the local well-posedness of
the fifth-order KP-II equation
in $H^{s_{1},s_{2}}(\R^{2})$ with $s_{1}>-\frac54,s_{2}\geq0$ and globally well-posed
in $H^{s_{1},0}(\R^{2})$ with $s_{1}>-\frac47$ with the aid of $I$-method.

By using  the  Fourier  restriction  norm method  and the Cauchy-Schwartz inequalities
as well as some calculus
inequalities,  Saut and Tzvetkov \cite{ST2000} established the global well-posedness of
Cauchy problem for the fifth order KP-I equation (\ref{1.01})
with initial data $u_{0} \in L^{2}(\R^{2})$ and finite energy.
By using the Fourier restriction  norm method and the dyadic decomposed Strichartz estimates,
 Chen et al. \cite{CLM} established the local well-posedness of  the problem (\ref{1.01})(\ref{1.02})
 in the interpolated  energy  space $E^{s}$ with $0<s\leq 1,$ where
\begin{align*}
E^{s}=\left\{u_{0}\in E^{s}:\|u_{0}\|_{E^{s}}
=\left\|\left(1+|\xi|^{2}+\left| {\mu}/{\xi}\right|\right)^{s}
\mathscr{F}_{xy}u_{0}(\xi,\mu)\right\|_{L^{2}}<\infty\right\}.
\end{align*}
In particular, Chen et al. established the global well-posedness of the problem  (\ref{1.01})(\ref{1.02})
 in the  energy  space $E^{1}$.
By using the Fourier restriction norm method and sufficiently exploiting
the geometric structure of the resonant set of (\ref{1.01}) to deal with the high-high frequency
interaction,  Li and Xiao \cite{LX} established the global well-posedness of
 the Cauchy problem  (\ref{1.01})(\ref{1.02})
in $L^{2}(\R^{2}).$ Guo et al. \cite{GHF} established the local well-posedness of the Cauchy problem for
(\ref{1.01}) in $H^{s,0}(\R^{2})$ with $s\geq -\frac{3}{4}.$
Yan et al. \cite{YLHD} proved that the Cauchy problem for (\ref{1.01}) is locally well-posed in $H^{s,0}(\R^{2})$ with $s>-\frac{3}{4}$
and the Cauchy problem for (\ref{1.01}) is globally well-posed in $H^{s,0}(\R^{2})$ with
$s>-\frac{6}{23}$ with the aid of $I$-method introduced in \cite{CKSTT2001,CKSTT2003}.
The method of \cite{YLHD} establishing local well-posedness is different from the method of \cite{GHF}.
Saut and Tzvetkov \cite{ST-CMP} have proved that the Cauchy problem for (\ref{1.01}) posed on $\mathbf{T}\times \R$ is globally
well-posed in the energy space. Compared to the fifth order KP-II equation, the structure of the fifth order KP-I equation is complicated.
The  reason is that
the resonant function of the fifth-order KP-I equation does not possess
 the same good property as its of fifth-order KP-II equation. More precisely,
  the resonant function of the fifth order KP-I equation
is
\begin{align}
R_{{\rm I}}(\xi_{1},\xi_{2},\mu_{1},\mu_{2})
:=&\phi(\xi,\mu)-\phi(\xi_{1},\mu_{1})-\phi(\xi_{2},\mu_{2})
       \nonumber\\
= & \frac{\xi_{1}\xi_{2}}{\xi}\left[5\xi^{2}(\xi_{1}^{2}+\xi_{1}\xi_{2}+\xi_{2}^{2})
      -\left(\frac{\mu_{1}}{\xi_{1}}-\frac{\mu_{2}}{\xi_{2}}\right)^{2}\right]
      \label{1.05}
\end{align}
and
the resonant function of the fifth order KP-II equation
is
\begin{align}
R_{{\rm I\!I}}(\xi_{1},\xi_{2},\mu_{1},\mu_{2})
:=& \phi(\xi,\mu)-\phi(\xi_{1},\mu_{1})-\phi(\xi_{2},\mu_{2})
      \nonumber\\
=&  \frac{\xi_{1}\xi_{2}}{\xi}\left[5\xi^{2}(\xi_{1}^{2}
       +\xi_{1}\xi_{2}
         +\xi_{2}^{2})
          +\left(\frac{\mu_{1}}{\xi_{1}}
           -\frac{\mu_{2}}{\xi_{2}}\right)^{2}\right].
            \label{1.06}
\end{align}
We remark that $R_{{\rm I}}(\xi_{1},\xi_{2},\mu_{1},\mu_{2})=0$ gives a surface, while  $R_{{\rm I\!I}}(\xi_{1},\xi_{2},\mu_{1},\mu_{2})$ will never be zero away from the origin.

In this paper, motivated by \cite{CKS-GAFA,ST2000,LX,ILM-CPAA},
by using the Fourier restriction norm method introduced in \cite{Beals,Bourgain93,KM,RR}
and developed in \cite{KPV1993,KPV1996}, the Cauchy-Schwartz inequality and Strichartz estimates
 as well as  suitable splitting of domains,
  we establish the local well-posedness of  the Cauchy problem for the fifth-order KP-I equation
 in the anisotropic Sobolev spaces $H^{s_{1},\>s_{2}}(\R^{2})$
  with $s_{1}>-\frac{9}{8}$ and $s_{2}\geq 0$; combining the local well-posness result
  of this paper with the I-method introduced in \cite{CKSTT2001,CKSTT2003},
   we established the  global well-posedness of the problem  in $H^{s_{1},\>0}(\R^{2})$ with $s_{1}>-\frac{4}{7}$.
Thus, our result considerably improves the result of \cite{GHF,LX,ST2000}.

We introduce some notations before presenting the main results.
 Throughout this paper, we assume that
$C$ is a positive constant which may vary from line to line. $a\sim b$ means that there exist constants $C_{j}>0(j=1,2)$ such that $C_{1}|b|\leq |a|\leq C_{2}|b|$.
 $a\gg b$ means that there exist a positive constant $C^{\prime}$ such that  $|a|> C^{\prime}|b|.$ $0<\epsilon\ll1$ means that $0<\epsilon<{10^{-4}}$.
 We define
  \begin{align*}
  &\langle\cdot\rangle:=1+|\cdot|,\\
  &\phi(\xi,\mu):=\xi^{5}+\frac{\mu^{2}}{\xi},\\
  &\sigma:=\tau+\phi(\xi,\mu),\sigma_{j}=\tau_{j}+\phi(\xi_{j},\mu_{j})(j=1,2),\\
  &\mathscr{F}u(\xi,\mu,\tau):=\frac{1}{(2\pi)^{\frac{3}{2}}}\int_{\SR^{3}}e^{-ix\xi-iy\mu-it\tau}u(x,y,t)dxdydt,\\
  &\mathscr{F}_{xy}f(\xi,\mu):=\frac{1}{2\pi}\int_{\SR^{2}}e^{-ix\xi-iy\mu}f(x,y)dxdy,\\
  &\mathscr{F}^{-1}u(\xi,\mu,\tau):=\frac{1}{(2\pi)^{\frac{3}{2}}}\int_{\SR^{3}}e^{ix\xi+iy\mu+it\tau}u(x,y,t)dxdydt,\\
  &D_{x}^{a}u(x,y,t):=\frac{1}{(2\pi)^{\frac{3}{2}}}\int_{\SR^{3}}|\xi|^{a}\mathscr{F}u(\xi,\mu,\tau)e^{ix\xi+iy\mu+it\tau}d\xi d\mu d\tau,\\
  &P^{2}u(x,y,t):=\frac{1}{(2\pi)^{\frac{3}{2}}}\int_{|\xi|\geq2}\int_{\SR^{2}}\mathscr{F}u(\xi,\mu,\tau)e^{ix\xi+iy\mu+it\tau}d\xi d\mu d\tau,\\
 & W(t)f:=\frac{1}{2\pi}\int_{\SR^{2}}e^{ix\xi+iy\mu+it\phi(\xi,\mu)}\mathscr{F}_{xy}f(\xi,\mu)d\xi d\mu.
  \end{align*}
   Let $\eta$ be a bump function with compact support in $[-2,2]\subset \R$
  and $\eta=1$ on $(-1,1)\subset \R$.
 For each integer $j\geq1$, we define $\eta_{j}(\xi)=\eta(2^{-j}\xi)-\eta(2^{1-j}\xi),$
  $\eta_{0}(\xi)=\eta(\xi),$
 $\eta_{j}(\xi,\mu,\tau)=\eta_{j}(\sigma),$ thus, $\sum\limits_{j\geq0}\eta_{j}(\sigma)=1.$
 $\psi(t)$ is a smooth function
 supported in $[0,2]$ and equals
 $1$ in $[0,1]$.
 Let $I\subset \R^{d}$, $\chi_{I}(x)=1$ if $x\in I$; $\chi_{I}(x)=0$
  if $x$ does not belong to $I$.

\noindent We define
 \begin{align*}
 \|f\|_{L_{t}^{r}L_{xy}^{p}}:=\left(\int_{\SR}\left(\int_{\SR^{2}}|f|^{p}dxdy\right)^{\frac{r}{p}}dt\right)^{\frac{1}{r}}.
 \end{align*}
 We denote by $H^{s_{1},s_{2}}(\R^{2})$ the anisotropic Sobolev space as follows:
 \begin{align*}
 H^{s_{1},s_{2}}(\R^{2})
 :=\left\{u_{0}\in \mathscr{S}^{'}(\R^{2}):
 \> \|u_{0}\|_{H^{s_{1},s_{2}}(\SR^{2})}
 =\left\|\langle\xi\rangle^{s_{1}}
 \langle\mu\rangle^{s_{2}}\mathscr{F}_{xy}u_{0}(\xi,\mu)\right\|_{L_{\xi\mu}^{2}}\right\}.
 \end{align*}
The Bourgain space
$
  X_{b}^{s_{1},s_{2}}
$ is defined by
$$
  X_{b}^{s_{1},s_{2}}
   := \left\{u\in  \mathscr{S}^{'}(\R^{3})
    :\, \|u\|_{X_{b}^{s_{1},s_{2}}}
     =  \Big\|\langle\xi\rangle^{s_{1}} \langle\mu\rangle^{s_{2}}
      \left\langle\sigma\right\rangle^{b}\mathscr{F}u(\xi,\mu,\tau)
       \Big\|_{L_{\tau\xi\mu}^{2} }<\infty\right\}.
$$
The space $ X_{b}^{s_{1},s_{2}}([0,T])$ denotes the restriction
 of $X_{b}^{s_{1},s_{2}}$ onto the finite time interval $[0,T]$ and
is equipped with the norm
 \begin{equation*}
    \|u\|_{X_{b}^{s_{1},s_{2}}([0,T])} =\inf \left\{\|g\|_{X_{b}^{s_{1},s_{2}}}
    :g\in X_{b}^{s_{1},s_{2}}, u(t)=g(t)
 \>\> {\rm for} \>  t\in [0,T]\right\}.
 \end{equation*}
For $s<0$ and $N\in\mathbf{ N}^{+},N\geq20$, we define an operator $I_{N}$  by
 $\mathscr{F}I_{N}u(\xi,\mu,\tau)=M(\xi)\mathscr{F}u(\xi,\mu,\tau)$,
where $M(\xi)=1$ if $|\xi|<N$; $M(\xi)=\left({|\xi|}/{N}\right)^{s}$ if $|\xi|\geq N.$

The main results of this paper are as follows.

\begin{Theorem}\label{Thm1}{\rm (Local well-posedness)} \
Let $|\xi|^{-1}\mathscr{F}_{xy}u_{0}(\xi,\mu)\in \mathscr{S}^{'}(\R^{2})$. Then,   the Cauchy problem  for (\ref{1.01}) is  locally well-posed in
 $H^{s_{1},\>s_{2}}(\R^{2})$ with $s_{1}>-{9}/{8},\>s_{2}\geq0.$
\end{Theorem}
\noindent{\bf Remark 1.}
We only consider the case of $-{9}/{8}<s_{1}<0,\> s_{2}\geq0$.
For $s_{1}\geq0,\>s_{2}\geq0$ the local well-posedness is proved by Li and Xiao \cite{LX}.
Lemmas 3.1 and 3.2 are the key ingredients in establishing the bilinear estimates in Lemmas 4.1 and 4.2.
Once Lemma 4.1 is proven to be valid, then we can combine it and Lemma 2.6 with the fixed point argument to obtain
the local wellposedness. Since the phase function $\phi(\xi,\mu)$
is singular at $\xi=0$, to  define  the  derivative  of  $W(t)u_{0}$,
the requirement $|\xi|^{-1}\mathscr{F}_{xy}u_{0}(\xi,\mu)\in \mathscr{S}^{'}(\R^{2})$  is necessary.

\begin{Theorem}\label{Thm2}{\rm(Global well-posedness)} \
 Let $|\xi|^{-1}\mathscr{F}_{xy}u_{0}(\xi,\mu)\in \mathscr{S}^{'}(\R^{2})$.
 Then  the Cauchy problem  for (\ref{1.01}) is  globally well-posed in $H^{s_{1},\>0}(\R^{2})$
with $s_{1}>-{4}/{7}.$
\end{Theorem}

\noindent{\bf Remark 2.}
We only consider the case of $-{4}/{7}<s_{1}<0,\>s_{2}\geq0$.
The case of $s_{1}\geq0,\>s_{2}\geq0$ is proved by Li and Xiao \cite{LX}.
For the fifth order KP-II equation, Isaza, L\'opez and  Mej\'{\i}a \cite{IM2011} obtained
the same result about the global well-posedness, that is, the Cauchy problem for the fifth order KP-II equation
is also globally well-posed in $H^{s_{1},\>s_{2}}(\R^{2})$ with $s_{1}>-\frac{4}{7},\>s_{2}\geq0.$

The rest of the paper is arranged as follows. In Section 2,  we give some
preliminaries. In Section 3, we establish two  $L^{2}$ bilinear estimates.
In Section 4, we establish three bilinear estimates.
 In Section 5, we prove
the local well-posedness. In Section 6, we firstly prove Lemma 6.1 which is a variation of Theorem 1.1,
then, we apply Lemmas 6.1, 4.2 and 2.6 to prove Theorem 1.2.

\bigskip
\bigskip

 \noindent{\large\bf 2. Preliminaries }

\setcounter{equation}{0}

\setcounter{Theorem}{0}

\setcounter{Lemma}{0}

\setcounter{section}{2}

This section is devoted to present Lemmas 2.1--2.6.

\begin{Lemma}\label{Lemma2.1}
 Let $b>|a|\geq 0$. Then, we have
 \begin{align}
 &\int_{-b}^{b}\frac{dx}{\langle x+a\rangle^{\frac{1}{2}}}\leq
  Cb^{\frac{1}{2}},\label{2.01}\\
 &\int_{\SR}\frac{dt}{\langle t\rangle^{\gamma}\langle t-a\rangle^{\gamma}}
 \leq C\langle a\rangle^{-\gamma},\gamma>1,\label{2.02}\\
 &\int_{\SR}\frac{dt}{\langle t\rangle^{\gamma}|t-a|^{\frac{1}{2}}}
 \leq C\langle a\rangle^{-\frac{1}{2}},\gamma\geq1,\label{2.03}\\
 &\int_{-K}^{K}\frac{dx}{|x|^{\frac{1}{2}}|a-x|^{\frac{1}{2}}}
 \leq C\frac{K^{\frac{1}{2}}}{|a|^{\frac{1}{2}}}.\label{2.04}
 \end{align}
 \end{Lemma}
 {\bf Proof.} The  conclusion of (\ref{2.01}) is given in (2.4) of Lemma 2.1 in \cite{ILM-CPAA}.
(\ref{2.02})-(\ref{2.03})  can be seen in Proposition 2.2 of \cite{ST2000}. (\ref{2.04})
can be seen in  \cite[Page 6562]{Hadac2008} .

This completes the proof of Lemma 2.1.
\begin{Lemma}\label{Lemma2.2}
Let $T\in (0,1)$ and $s_{1},s_{2} \in \R$ and $-\frac{1}{2}<b^{\prime}\leq0\leq b\leq b^{\prime}+1$.
Then, for $h\in X_{b^{\prime}}^{s_{1},s_{2}},$  we have
\begin{align}
&\left\|\psi(t)S(t)\phi\right\|_{X_{b}^{s_{1},s_{2}}}\leq C\|\phi\|_{H^{s_{1},\>s_{2}}},\label{2.05}\\
&\left\|\psi\left(\frac{t}{T}\right)\int_{0}^{t}S(t-\tau)h(\tau)d\tau\right\|_{X_{b}^{s_{1},\>s_{2}}}\leq C
T^{1+b^{\prime}-b}\|h\|_{X_{b^{\prime}}^{s_{1},\>s_{2}}}.\label{2.06}
\end{align}
\end{Lemma}

For the proof of Lemma 2.2, we refer  readers to \cite{G,Bourgain93,KPV1993}
 and     \cite[Lemma  1.7 and Lemma 1.9]{Grunrock}.

\begin{Lemma}\label{Lemma2.3}
Let $b>\frac{1}{2}$. Then,
\begin{align}
\left\|D_{x}^{\frac{1}{4}}u\right\|_{L_{t}^{4}L_{xy}^{4}(\SR^{3})}
\leq C\|u\|_{X_{b}^{0,0}}.\label{2.07}
\end{align}
\end{Lemma}

For the proof of Lemma 2.3, we refer  readers to \cite[Theorem 3.1]{Hadac2008}.

\begin{Lemma}\label{Lemma2.4}
Let
\begin{align*}
|\sigma-\sigma_{1}-\sigma_{2}|
=&\left|\xi\xi_{1}\xi_{2}(5\xi^{2}-5\xi\xi_{1}+5\xi_{1}^{2})-\frac{\xi_{1}\xi_{2}}{\xi}
\left|\frac{\mu_{1}}{\xi_{1}}-\frac{\mu_{2}}{\xi_{2}}\right|^{2}\right|\nonumber\\
\geq   &
\frac{\left|\xi\xi_{1}\xi_{2}(5\xi^{2}-5\xi\xi_{1}+5\xi_{1}^{2})\right|}{4}
\end{align*}
and
\begin{align*}
&\quad\mathscr{F}P_{\frac{1}{4}}(u_{1},u_{2})(\xi,\mu,\tau)
 \nonumber\\
& =\int_{\SR^{3}}
   \chi_{|\xi_{1}|\leq \frac{|\xi_{2}|}{4}}(\xi_{1},\mu_{1},\tau_{1},\xi,\mu,\tau)
   \prod\limits_{j=1}^{2}
   \mathscr{F}u_{j}(\xi_{j},\mu_{j},\tau_{j})d\xi_{1}d\mu_{1}d\tau_{1}.
\end{align*}
For $b>\frac{1}{2}$, we have
\begin{align}
\left\|P_{\frac{1}{4}}(u_{1},u_{2})\right\|_{L_{txy}^{2}}
\leq C\left\||D_{x}|^{\frac{1}{2}}u_{1}\right\|_{X_{b}^{0,0}}\left\||D_{x}|^{-1}u_{2}\right\|_{X_{b}^{0,0}}.\label{2.08}
\end{align}
\end{Lemma}

\noindent {\bf Proof.} Let
\begin{align*}
f_{1}(\xi_{1},\mu_{1},\tau_{1})=|\xi_{1}|^{\frac{1}{2}}\langle \sigma_{1}\rangle^{b}\mathscr{F}u_{1}(\xi_{1},\mu_{1},\tau_{1}),
f_{2}(\xi_{2},\mu_{2},\tau_{2})=|\xi_{2}|^{-1}\langle \sigma_{2}\rangle^{b}\mathscr{F}u_{2}(\xi_{2},\mu_{2},\tau_{2}).
\end{align*}
 To obtain
(\ref{2.08}), it suffices to prove that
\begin{align}
&\bigg\|\int_{\SR^{3}}
\frac{|\xi_{1}|^{-\frac{1}{2}}|\xi_{2}|f_{1}(\xi_{1},\mu_{1},\tau_{1})
f_{2}(\xi_{2},\mu_{2},\tau_{2})}
{\prod\limits_{j=1}^{2}\langle\sigma_{j}\rangle^{b}}
d\xi_{1}d\mu_{1}d\tau_{1}\bigg\|_{L_{\tau\xi\mu}^{2}}\nonumber\\
\leq & C\prod_{j=1}^{2}\|f_{j}\|_{L_{\tau\xi\mu}^{2}}.\label{2.09}
\end{align}
To obtain (\ref{2.09}), by duality, it suffices to prove that
\begin{align}
&  \biggl\|\int_{\SR^{3}}
\frac{|\xi_{1}|^{-\frac{1}{2}}|\xi_{2}|f_{1}(\xi_{1},\mu_{1},\tau_{1})
f_{2}(\xi_{2},\mu_{2},\tau_{2})f(\xi,\mu,\tau)}
{\prod\limits_{j=1}^{2}\langle\sigma_{j}\rangle^{b}}
d\xi_{1}d\mu_{1}d\tau_{1}\biggr\|_{L_{\xi\mu\tau}^{2}}\nonumber\\
\leq &  C\|f\|_{L_{\tau\xi\mu}^{2}}
\prod_{j=1}^{2}\|f_{j}\|_{L_{\tau\xi\mu}^{2}}.\label{2.010}
\end{align}
We define
\begin{align}
I(\xi,\mu,\tau):=\int_{\SR^{3}}
\frac{|\xi_{1}|^{-1}|\xi_{2}|^{2}}{\prod\limits_{j=1}^{2}\langle\sigma_{j}\rangle^{b}}
d\xi_{1}d\mu_{1}d\tau_{1}.\label{2.011}
\end{align}
For fixed $(\xi,\mu,\tau),$ we make the change of variables ${\rm L}:(\xi_{1},\mu_{1},\tau_{1})
\longrightarrow (\Delta,\sigma_{1},\sigma_{2})$,
where
\begin{align*}
&\Delta:=\xi\xi_{1}\xi_{2}
(5\xi^{2}-5\xi\xi_{1}+5\xi_{1}^{2}),\\
&\sigma_{1}:=\tau_{1}+\phi(\xi_{1},\mu_{1}),\sigma_{2}:=\tau_{2}+\phi(\xi_{2},\mu_{2}).
\end{align*}
By using a direct computation, since $\sigma=\tau+\phi(\xi,\mu),$ we have
\begin{align}
\sigma_{1}+\sigma_{2}-\sigma=-\Delta+\frac{(\xi_{1}\mu_{2}-\mu_{1}\xi_{2})^{2}}{\xi\xi_{1}\xi_{2}}.\label{2.012}
\end{align}
Thus, we have  the Jacobian determinant
equals
\begin{align}
   \frac{\partial(\Delta,\sigma_{1},\sigma_{2})}
        {\partial(\xi_{1},\mu_{1},\tau_{1})}
=&   -10\left(\xi_{1}^{2}-\xi_{2}^{2}\right)
    \left(\xi_{1}^{2}+\xi_{2}^{2}\right)
    \left(\frac{\mu_{1}}{\xi_{1}}-\frac{\mu_{2}}{\xi_{2}}\right)\nonumber\\
=&  -10\left(\xi_{1}^{2}-\xi_{2}^{2}\right)
    \left(\xi_{1}^{2}+\xi_{2}^{2}\right)
    (\sigma_{1}+\sigma_{2}-\sigma+\Delta)^{\frac{1}{2}}
    \left(\frac{\xi}{\xi_{1}\xi_{2}}\right)^{\frac{1}{2}}.
    \label{2.013}
\end{align}
Notice that it is possible to divide the integration into a finite number of open subsets $W_{i}$ such that
${\rm L}$ is an injective $C^{1}$-function in $W_{i}$ with non-zero Jacobian determinant.
From (\ref{2.013}), since $|\xi_{1}|\leq \frac{|\xi_{2}|}{4}$ and $|\Delta|\sim |\xi_{1}||\xi_{2}|^{4},$ we have
\begin{align}
   \left|\frac{\partial(\Delta,\sigma_{1},\sigma_{2})}{\partial(\xi_{1},\mu_{1},\tau_{1})}\right|
=&  10\left|\left(\xi_{1}^{2}-\xi_{2}^{2}\right)
   \left(\xi_{1}^{2}+\xi_{2}^{2}\right)\left(\frac{\mu_{1}}{\xi_{1}}-\frac{\mu_{2}}{\xi_{2}}\right)\right|\nonumber\\
=&  10\left|\left(\xi_{1}^{2}-\xi_{2}^{2}\right)
   \left(\xi_{1}^{2}+\xi_{2}^{2}\right)
   (\sigma_{1}+\sigma_{2}-\sigma+\Delta)^{\frac{1}{2}}
   \left(\frac{\xi}{\xi_{1}\xi_{2}}\right)^{\frac{1}{2}}\right|\nonumber\\
\sim & |\xi_{1}|^{-1}|\xi_{2}|^{2}|\Delta|^{\frac{1}{2}}
    \left|\sigma_{1}+\sigma_{2}-\sigma+\Delta\right|^{\frac{1}{2}}.
    \label{2.014}
\end{align}
Since $\left|\sigma_{1}+\sigma_{2}-\sigma\right|\geq \frac{|\Delta|}{4}$,
 by using the change of variables
 $(\xi_{1},\mu_{1},\tau_{1})\longrightarrow
 (\Delta,\sigma_{1},\sigma_{2})$ and (\ref{2.04}), we have
\begin{align}
I(\xi,\mu,\tau):=&\int_{\SR^{3}}\chi_{|\xi_{1}|\leq \frac{|\xi_{2}|}{4}}
\frac{|\xi_{2}|^{2}|\xi_{1}|^{-1}}{\prod\limits_{j=1}^{2}\langle\sigma_{j}\rangle^{2b}}
d\xi_{1}d\mu_{1}d\tau_{1}\nonumber\\
\leq &C\int_{\SR^{3}}\frac{\chi_{|\Delta|\leq 4|\sigma_{1}+\sigma_{2}-\sigma|
  d\Delta d\sigma_{1}d\sigma_{2}}}{|\Delta|^{\frac{1}{2}}
\left|\sigma_{1}+\sigma_{2}-\sigma+\Delta\right|^{\frac{1}{2}}\prod\limits_{j=1}^{2}\langle\sigma_{j}\rangle^{2b}}\nonumber\\
=&\int_{\SR^{2}}\left(\int_{\SR}\frac{\chi_{|\Delta|\leq 4|\sigma_{1}+\sigma_{2}-\sigma|d\Delta }}{|\Delta|^{\frac{1}{2}}
\left|\sigma_{1}+\sigma_{2}-\sigma+\Delta\right|^{\frac{1}{2}}}\right)\frac{d\sigma_{1}d\sigma_{2}}
{\prod\limits_{j=1}^{2}\langle\sigma_{j}\rangle^{2b}}\nonumber\\
\leq & C\int_{\SR^{2}}\frac{d\sigma_{1}d\sigma_{2}}
{\prod\limits_{j=1}^{2}\langle\sigma_{j}\rangle^{2b}}\leq C.\label{2.015}
\end{align}
Combining (\ref{2.010}) with (\ref{2.015}), by using the Cauchy-Schwartz inequality twice, we have
\begin{align}
& \int_{\SR^{3}}\chi_{|\xi_{1}|\leq \frac{|\xi_{2}|}{4}}
\frac{|\xi_{2}||\xi_{1}|^{-\frac{1}{2}}
f_{1}(\xi_{1},\mu_{1},\tau_{1})
f_{2}(\xi_{2},\mu_{2},\tau_{2})f(\xi,\mu,\tau)}
{\prod\limits_{j=1}^{2}\langle\sigma_{j}\rangle^{b}}
d\xi_{1}d\mu_{1}d\tau_{1}\nonumber\\
\leq& C\left[\sup \limits_{\xi,\mu,\tau}I(\xi,\mu,\tau)\right]^{\frac{1}{2}}
\|f\|_{L_{\tau\xi\mu}^{2}}
\prod_{j=1}^{2}\|f_{j}\|_{L_{\tau\xi\mu}^{2}}\nonumber\\
\leq& C\|f\|_{L_{\tau\xi\mu}^{2}}
\prod_{j=1}^{2}\|f_{j}\|_{L_{\tau\xi\mu}^{2}} .\label{2.016}
\end{align}

This completes the proof of Lemma 2.4.

\begin{Lemma}\label{Lemma2.5}
Let
\begin{align*}
|\sigma-\sigma_{1}-\sigma_{2}|
=&\left|\xi\xi_{1}\xi_{2}(5\xi^{2}-5\xi\xi_{1}+5\xi_{1}^{2})-\frac{\xi_{1}\xi_{2}}{\xi}
\left|\frac{\mu_{1}}{\xi_{1}}-\frac{\mu_{2}}{\xi_{2}}\right|^{2}\right|\nonumber\\   \geq & \frac{\left|\xi\xi_{1}\xi_{2}(5\xi^{2}
      -5\xi\xi_{1}+5\xi_{1}^{2})\right|}{4},
\end{align*}
 $b>\frac{1}{2}$ and
$
G(\xi_{1},\mu_{1},\tau_{1},\xi,\mu,\tau)=f_{1}(\xi_{1},\mu_{1},\tau_{1})
f_{2}(\xi-\xi_{1},\mu-\mu_{1},\tau-\tau_{1})f(\xi,\mu,\tau),
 $
 we have
\begin{align}
&\biggl|\int_{\SR^{6}}\frac{|\xi_{1}|^{-\frac{1}{2}}
|\xi_{2}|G(\xi_{1},\mu_{1},\tau_{1},\xi,\mu,\tau)}
{\prod\limits_{j=1}^{2}
\langle\sigma_{j}\rangle^{b}}d\xi_{1}d\mu_{1}d\tau_{1}d\xi d\mu d\tau\biggr|
\leq   C\|f\|_{L_{\tau\xi\mu}^{2}}
  \prod\limits_{j=1}^{2}\|f_{j}\|_{L_{\tau\xi\mu}^{2}} , \label{2.017}\\
&\left|\int_{\SR^{6}}\frac{|\xi_{1}|^{-\frac{1}{2}}
    |\xi|G(\xi_{1},\mu_{1},\tau_{1},\xi,\mu,\tau)}
     {\langle\sigma_{1}\rangle^{b}
      \langle\sigma\rangle^{b}}d\xi_{1}d\mu_{1}d\tau_{1}d\xi d\mu d\tau\right|
\leq   C\|f\|_{L_{\tau\xi\mu}^{2}}
      \prod_{j=1}^{2}\|f_{j}\|_{L_{\tau\xi\mu}^{2}} , \label{2.018}\\
&\left|\int_{\SR^{6}}\frac{|\xi|^{-\frac{1}{2}}
  |\xi_{2}|G(\xi_{1},\mu_{1},\tau_{1},\xi,\mu,\tau)}
   {\langle\sigma\rangle^{b}\langle\sigma_{2}\rangle^{b}}d\xi_{1}d\mu_{1}d\tau_{1}d\xi d\mu d\tau\right|
\leq   C\|f\|_{L_{\tau\xi\mu}^{2}}
     \prod\limits_{j=1}^{2}\|f_{j}\|_{L_{\tau\xi\mu}^{2}} , \label{2.019}\\
&\left|\int_{\SR^{6}}\frac{|\xi|^{-\frac{1}{2}}
  |\xi_{1}|G(\xi_{1},\mu_{1},\tau_{1},\xi,\mu,\tau)}
   {\langle\sigma\rangle^{b}\langle\sigma_{1}\rangle^{b}}d\xi_{1}d\mu_{1}d\tau_{1}d\xi d\mu d\tau\right|
\leq   C\|f\|_{L_{\tau\xi\mu}^{2}}
       \prod\limits_{j=1}^{2}\|f_{j}\|_{L_{\tau\xi\mu}^{2}} .\label{2.020}
\end{align}
\end{Lemma}
\noindent{\bf Proof.}  We firstly prove (\ref{2.017}). When $\frac{|\xi_{2}|}{4}\geq |\xi_{1}|$, from Lemma 2.4,
we have  (\ref{2.017}) is valid. When $\frac{|\xi_{2}|}{4}< |\xi_{1}|$, since $|\xi_{1}|^{-\frac{1}{2}}|\xi_{2}|\leq C
|\xi_{1}|^{\frac{1}{4}}|\xi_{2}|^{\frac{1}{4}}$, from Lemma 2.3, we know that (\ref{2.017}) is valid.
Let $\xi_{1}=\xi_{1}^{\prime},\mu_{1}=\mu_{1}^{\prime},\tau_{1}=\tau_{1}^{\prime}$
and $-\xi_{2}=\xi^{\prime},-\tau_{2}=\tau^{\prime},-(\mu-\mu_{1})=\mu^{\prime}$
and
 $-\xi=\xi^{\prime}-\xi_{1}^{\prime},-\mu=\mu^{\prime}-\mu_{1}^{\prime},-\tau=\tau^{\prime}-\tau_{1}^{\prime}$
and $\sigma_{2}^{\prime}=\tau_{2}^{\prime}-\phi^{\prime}(\xi_{2}^{\prime},\mu_{2}^{\prime}),\sigma_{1}=\sigma_{1}^{\prime}
=\tau_{1}^{\prime}-\phi(\xi_{1}^{\prime},\mu_{1}^{\prime})$. Thus, $-\sigma=\sigma_{2}^{\prime},\sigma_{1}=\sigma_{1}^{\prime}.$
Let
\begin{align*}
H(\xi_{1}^{\prime},\mu_{1}^{\prime},\tau_{1}^{\prime},\xi^{\prime},\mu^{\prime},
\tau^{\prime})=f_{1}(\xi_{1}^{\prime},\mu_{1}^{\prime},\tau_{1}^{\prime})
f_{2}(-\xi^{\prime},-\mu^{\prime},-\tau^{\prime})f(-\xi_{2}^{\prime},
-\mu_{2}^{\prime},-\tau_{2}^{\prime}).
\end{align*}
To obtain (\ref{2.018}), it suffices to prove that
\begin{align}
&\left|\int_{\SR^{6}}\frac{|\xi_{1}^{\prime}|^{-\frac{1}{2}}
  |\xi_{2}^{\prime}||H(\xi_{1}^{\prime},\mu_{1}^{\prime},
  \tau_{1}^{\prime},\xi^{\prime},\mu^{\prime},
   \tau^{\prime})}
   {\langle\sigma_{1}^{\prime}\rangle^{b}\langle\sigma_{2}^{\prime}\rangle^{b}}
d\xi_{1}^{\prime}d\mu_{1}^{\prime}d\tau_{1}^{\prime}d\xi^{\prime}
d\mu ^{\prime} d\tau^{\prime}\right|\nonumber\\
\leq & C\|f\|_{L_{\tau\xi\mu}^{2}}
   \prod\limits_{j=1}^{2}
   \|f_{j}\|_{L_{\tau\xi\mu}^{2}}\label{2.021}.
\end{align}
Obviously, (\ref{2.021}) follows from (\ref{2.017}). By using a
  proof similar to (\ref{2.018}), we obtain that (\ref{2.019})-(\ref{2.020}) are valid.

This ends the proof of Lemma 2.5.

\begin{Lemma}\label{Lemma2.6}
Let $0<b_{1}<b_{2}<\frac12$. Then, we have
\begin{align}
&\left\|\chi_{I}(\cdot)u\right\|_{X_{b_{1}}^{0,0}}\leq C\left\|u\right\|_{X_{b_{2}}^{0,0}},\label{2.022}\\
&\left\|\chi_{I}(\cdot)u\right\|_{X_{-b_{2}}^{0,0}}\leq C\left\|u\right\|_{X_{-b_{1}}^{0,0}}.\label{2.023}
\end{align}
\end{Lemma}

For the proof of Lemma 2.6, we refer the readers to  Lemma 3.1. of \cite{IMEJDE}.

\bigskip
\bigskip

\noindent{\large\bf 3. $L^{2}$-bilinear estimates}

\setcounter{equation}{0}

 \setcounter{Theorem}{0}

\setcounter{Lemma}{0}

 \setcounter{section}{3}
 Inspired by the idea of  Lemma 5.1 of \cite{IKT}, we give the proof of Lemma 3.1.
 For $k\in \mathbf{Z}$ and  $l,j\in \R$, we define
 \begin{align*}
&\hspace{-1cm} D_{k,l,j}:=\left\{(\xi,\mu,\tau):|\xi| \in [2^{k-1},2^{k+1}],|\mu|\leq 2^{l},|\tau+\phi(\xi,\mu)|<2^{j}\right\},\nonumber\\
&\hspace{-1cm}
D_{k,\infty,j}:=\bigcup\limits_{l\in\mathbf{Z}}D_{k,l,j}.
\end{align*}
\begin{Lemma}\label{Lemma3.1}
Assume $\alpha\in \R$ and $k_{1},k_{2},k_{3}\in\mathbf{Z}$, $k_{\rm max}={\rm max}\left\{k_{1},k_{2},k_{3}\right\}$
 and $k_{\rm min}={\rm min}\left\{k_{1},k_{2},k_{3}\right\}$
 and $j_{1},j_{2},j_{3}\in \mathbf{Z}_{+},j_{\rm max}={\rm max}\left\{j_{1},j_{2},j_{3}\right\}$
and $f_{i}:\R^{3}\rightarrow \R$  are  $L^{2}$  functions supported in $D_{k_{i},\infty,j_{i}}(i=1,2,3)$. We assume that
\begin{align}
&R_{{\rm I}}(\xi_{1},\xi_{2},\mu_{1},\mu_{2})\leq 2^{k_{1}+k_{2}+k_{3}+2k_{\rm max}-60},\label{3.01}\\
&j_{\rm max}\leq k_{1}+k_{2}+k_{3}+2k_{\rm max}-60.\label{3.02}
\end{align}
\noindent {\rm (1)} If
 $|k_{1}-k_{2}|\leq 5$, $k_{1}\geq20,$
then, we have
\begin{align}
&\int_{\SR^{3}}(f_{1}*f_{2})f_{3}d\xi_{1}d\mu_{1}d\tau_{1}
\leq C2^{\frac{j_{1}+j_{2}+j_{3}}{2}}2^{-\frac{7}{4}(k_{1}+k_{2})+\frac{k_{3}}{2}}
   \prod\limits_{j=1}^{3}\|f_{j}\|_{L^{2}} .\label{3.03}
\end{align}
{\rm (2)} If
 $k_{2}-10\geq k_{1}$ and  $|k_{2}-k_{3}|\leq 5$, $k_{2}\geq20,$
then, we have
\begin{align}
&\int_{\SR^{3}}(f_{1}*f_{2})f_{3}d\xi_{1}d\mu_{1}d\tau_{1}
\leq C2^{\frac{j_{1}+j_{2}+j_{3}}{2}}2^{-\frac{5}{2}(k_{2}+k_{3})-\frac{k_{1}}{2}}
   \prod\limits_{j=1}^{3}\|f_{j}\|_{L^{2}} .\label{3.04}
\end{align}
\end{Lemma}
\noindent{\bf Proof.} \
First we prove \eqref{3.03}. From  (5.4) of \cite{IKT}, we have
\begin{align}
\int_{\SR^{3}}(f_{1}*f_{2})f_{3}=\int_{\SR^{3}}(\tilde{f}_{1}*\tilde{f}_{3})f_{2}=\int_{\SR^{3}}(\tilde{f}_{2}*\tilde{f}_{3})f_{1}.\label{3.05}
\end{align}
where
$\tilde{f}_{i}=f_{i}(-\xi,-\mu,-\tau)(i=1,2)$.
 Due to the symmetry in  (\ref{3.05}),  without loss of  generality,
we may assume $j_{3}={\rm max}\left\{j_{1},j_{2},j_{3}\right\}.$
We define $$f_{i}^{\#}(\xi_{i},\mu_{i},\theta_{i}):=f_{i}(\xi_{i},\mu_{i},\theta_{i}-\phi(\xi_{i},\mu_{i}))(i=1,2,3).$$
Obviously, $\|f_{i}^{\#}\|_{L^{2}}=\|f_{i}\|_{L^{2}}.$
The left-hand side of (\ref{3.03}) can be rewritten as follows:
\begin{align}
\int_{\SR^{6}}\biggl(\prod\limits_{i=1}^{2}f_{i}^{\#}(\xi_{i},\mu_{i},\theta_{i})\biggr)
f_{3}^{\#}(\xi_{1}+\xi_{2},\mu_{1}+\mu_{2},R_{{\rm I}}(\xi_{1},\mu_{1},\xi_{2},\mu_{2}))
d\xi_{1}d\xi_{2}d\mu_{1}d\mu_{2}d\theta_{1}d\theta_{2},\label{3.06}
\end{align}
where $R_{{\rm I}}(\xi_{1},\xi_{2},\mu_{1},\mu_{2})$ is the resonant function
defined as in (\ref{1.05}).
The functions $f_{i}^{\#}$ $(i=1,2)$ are supported in the sets
$$
   \left\{(\xi_{i},\mu_{i},\theta_{i}):
       \xi _{i}\in \tilde{I}_{k_{i}},\mu \in \R,
         |\theta_{i}|\leq 2^{j_{i}}\right\}
$$
 and $f_{3}^{\#}$ is supported in the set
$$
  \left\{(\xi_{3},\mu_{3},\theta_{3}):
   \xi_{3} \in \tilde{I}_{k_{3}},\mu_{3} \in \R,
    |\theta_{3}|\leq 2^{j_{3}}\right\}.
$$

We will prove that if $g_{i}:\R^{2}\longrightarrow\R$ are $L^{2}$ functions supported in $\tilde{I}_{k_{i}}\times \R(i=1,2)$
and $g:\R^{2}\longrightarrow\R$ are $L^{2}$ functions supported in $\tilde{I}_{k}\times \R\times[-2^{j},2^{j}]$,
then, we have
\begin{align}
&\int_{\SR^{4}}\biggl(
    \prod\limits_{j=1}^{2}g_{j}(\xi_{j},\mu_{j})\biggr)
g(\xi_{1}+\xi_{2},\mu_{1}+\mu_{2},R_{{\rm I}}(\xi_{1},\mu_{1},\xi_{2},\mu_{2}))
   d\xi_{1}d\xi_{2}d\mu_{1}d\mu_{2}\nonumber\\
\leq & C2^{\frac{j}{2}}2^{-\frac{7}{4}(k_{1}+k_{2})+\frac{k_{3}}{2}}\|g\|_{L^{2}}
\prod\limits_{j=1}^{2}\|g_{j}\|_{L^{2}} .\label{3.07}
\end{align}

When $j\leq k_{1}+k_{2}+k_{3}+2k_{\rm max}-60$ is valid, we may assume
 that the integral in the left-hand side of (\ref{3.07})
is taken over the set
\begin{align}
\mathcal {R}_{++}=\left\{(\xi_{1},\mu_{1},\xi_{2},\mu_{2}):\xi_{1}+\xi_{2}\geq0,
 \frac{\mu_{1}}{\xi_{1}}-\frac{\mu_{2}}{\xi_{2}}\geq0\right\}\label{3.08}
\end{align}
since other case can be proved similarly to case (\ref{3.08}).
\noindent We make the changes of variables
\begin{align}
\mu_{1}=\sqrt{5}\xi_{1}^{3}+\beta_{1}\xi_{1},\mu_{2}=-\sqrt{5}\xi_{2}^{3}+\beta_{2}\xi_{2}.\label{3.09}
\end{align}
From (\ref{3.09}), we have
\begin{align}
\frac{\mu_{1}}{\xi_{1}}-\frac{\mu_{2}}{\xi_{2}}=\sqrt{5}\left[\xi_{1}^{2}+\xi_{2}^{2}\right]+\beta_{1}-\beta_{2}
=\sqrt{5}(\xi_{1}^{2}+\xi_{2}^{2})+\beta_{1}-\beta_{2}\geq0.\label{3.010}
\end{align}
From (\ref{3.010}), we know that
\begin{align}
\beta_{1}-\beta_{2}\geq-\sqrt{5}(\xi_{1}^{2}+\xi_{2}^{2}).\label{3.011}
\end{align}
By using  the assumption upon $g$ and the definition of $R_{{\rm I}}(\xi_{1},\mu_{1},\xi_{2},\mu_{2})$, we infer that
\begin{align}
 \left|(\beta_{1}-\beta_{2})^{2}+2\sqrt{5}(\beta_{1}-\beta_{2})(\xi_{1}^{2}+\xi_{2}^{2})
   +5\xi_{1}\xi_{2}(4\xi_{1}\xi_{2}-3\xi^{2})\right|
\leq 2^{j+k-k_{1}-k_{2}+3}.\label{3.012}
\end{align}
The left hand side of  (\ref{3.07})  can be bounded by
\begin{align}
& C2^{k_{1}+k_{2}}
  \int_{S}h_{1}(\xi_{1},\sqrt{5}\xi_{1}^{3}+\beta_{1}\xi_{1})
h_{2}(\xi_{2},-\sqrt{5}\xi_{2}^{3}+\beta_{2}\xi_{2})\nonumber\\   &\quad\times
h(\xi_{1}+\xi_{2},\sqrt{5}(\xi_{1}^{3}-\xi_{2}^{3})+\beta_{1}\xi_{1}+\beta_{2}\xi_{2},\tilde{R}_{1}(\xi_{1},\beta_{1},\xi_{2},\beta_{2})
d\xi_{1}d\xi_{2}d\beta_{1}d\beta_{2},\label{3.013}
\end{align}
where
\begin{align}
S=\left\{(\xi_{1},\beta_{1},\xi_{2},\beta_{2}):\xi_{1}+\xi_{2}\geq0,\beta_{1}-\beta_{2} \ \ {\rm satisfies}
 \ \ (3.11)-(3.12)\right\}.\label{3.014}
\end{align}
and
\begin{align}
 &  \tilde{R}_{1}(\xi_{1},\beta_{1},\xi_{2},\beta_{2})
     \nonumber\\
=& (\beta_{1}-\beta_{2})^{2}
    +2\sqrt{5}(\beta_{1}-\beta_{2})(\xi_{1}^{2}
     +\xi_{2}^{2})+5\xi_{1}\xi_{2}(4\xi_{1}\xi_{2}-3\xi^{2}).\label{3.015}
\end{align}
We define the functions $h_{i}:\R^{2}\longrightarrow \R$ supported in $\tilde{I}_{k_{i}}\times \R(i=1,2)$
\begin{align}
&h_{1}(\xi_{1},\beta_{1})=2^{\frac{k_{1}}{2}}g_{1}(\xi_{1},\sqrt{5}\xi_{1}^{3}+\beta_{1}\xi_{1}),
 \label{3.016}\\ &h_{2}(\xi_{2},\beta_{2})=2^{\frac{k_{2}}{2}}g_{2}(\xi_{2},-\sqrt{5}\xi_{2}^{3}+\beta_{2}\xi_{2}).\label{3.017}
\end{align}
with $\|h_{i}\|_{L^{2}}\approx\|g_{i}\|_{L^{2}}(i=1,2)$.

To prove (\ref{3.07}),  it suffices to prove that
\begin{align}
& 2^{\frac{k_{1}+k_{2}}{2}}
  \int_{\tilde{S}}h_{1}(\xi_{1},\beta_{1})h_{2}(\xi_{2},\beta_{2})
    \nonumber\\
&\qquad \times
  h(\xi_{1}+\xi_{2},\sqrt{5}(\xi_{1}^{3}-\xi_{2}^{3})+\beta_{1}\xi_{1}+\beta_{2}\xi_{2},
   \tilde{R}_{1}(\xi_{1},\beta_{1},\xi_{2},\beta_{2})d\xi_{1}d\xi_{2}d\beta_{1}d\beta_{2}\nonumber\\
\leq &  C2^{\frac{j}{2}}
   2^{-\frac{7}{4}(k_{1}+k_{2})+\frac{k_{3}}{2}}
  \|h\|_{L^{2}}
   \prod\limits_{j=1}^{2}\|h_{j}\|_{L^{2}} .\label{3.018}
\end{align}
Combining (\ref{3.011}) with (\ref{3.012}),  we have
\begin{align}\nonumber
  \frac{\sqrt{5}(B_{1}-B_{2})}{\xi_{1}^{2}+\xi_{2}^{2}
   +\sqrt{\xi^{2}(\xi^{2}-\xi_{1}\xi_{2})-B_{2}}}
&  \leq \beta_{1}-\beta_{2} \\
&   \leq \frac{\sqrt{5}(B_{1}+B_{2})}
    {\xi_{1}^{2}+\xi_{2}^{2}
     +\sqrt{\xi^{2}(\xi^{2}-\xi_{1}\xi_{2})+B_{2}}}.\label{3.019}
\end{align}
where
\begin{align}
 B_{1}=\xi_{1}\xi_{2}
  \left[3\xi_{1}^{2}+2\xi_{1}\xi_{2}+3\xi_{2}^{2}\right],
\quad B_{2}=\frac{2^{j+k-k_{1}-k_{2}+3}}{5}.\label{3.020}
\end{align}
Now we claim that the following inequality is valid
\begin{align}
\left|\beta-\frac{\sqrt{5}B_{1}}{\xi_{1}^{2}
  +\xi_{2}^{2}+\sqrt{\xi^{2}(\xi^{2}-\xi_{1}\xi_{2})}}\right|\leq B_{3},\label{3.021}
\end{align}
where
\begin{align}
B_{3}=2^{j-\frac{3(k_{1}+k_{2})}{2}+10}.\label{3.022}
\end{align}
When $\xi_{1}\xi_{2}\geq0$, we have
\begin{align}
\sqrt{\xi^{2}(\xi^{2}-\xi_{1}\xi_{2}-\frac{3\alpha}{5})\pm B_{2}}+\sqrt{\xi^{2}(\xi^{2}-\xi_{1}\xi_{2})}\geq \xi^{2}.\label{3.023}
\end{align}
 By using a direct computation,  since
\begin{align}
 B_{1}\leq 3\xi^{4},\quad
   \sqrt{\xi^{2}(\xi^{2}-\xi_{1}\xi_{2})\pm
B_{2}}\sqrt{\xi^{2}(\xi^{2}-\xi_{1}\xi_{2})}
  \geq \frac{B_{1}}{3},\label{3.024}
\end{align}
 we have
\begin{align}
& \left|\frac{B_{1}}
{\xi_{1}^{2}+\xi_{2}^{2}+\sqrt{\xi^{2}(\xi^{2}-\xi_{1}\xi_{2})}}
-\frac{B_{1}}{\xi_{1}^{2}+\xi_{2}^{2}+\sqrt{\xi^{2}(\xi^{2}-\xi_{1}\xi_{2})\pm B_{2}}}\right|\nonumber\\
=& \frac{B_{1}B_{2}}{\left[\xi_{1}^{2}+\xi_{2}^{2}+\sqrt{\xi^{2}(\xi^{2}-\xi_{1}\xi_{2})\pm B_{2}}\right]
\left[\xi_{1}^{2}+\xi_{2}^{2}+\sqrt{\xi^{2}(\xi^{2}-\xi_{1}\xi_{2})}\right]}\nonumber\\   & \qquad \times
\frac{1}{\left[\sqrt{\xi^{2}(\xi^{2}-\xi_{1}\xi_{2})\pm B_{2}}
+\sqrt{\xi^{2}(\xi^{2}-\xi_{1}\xi_{2})}\right]}\nonumber\\
\leq & \frac{B_{3}}{10}.\label{3.025}
\end{align}
By using a direct computation,  we have
\begin{align}
\frac{B_{2}}{\xi_{1}^{2}+\xi_{2}^{2}+\sqrt{\xi^{2}(\xi^{2}-\xi_{1}\xi_{2})\pm B_{2}}}
\leq \frac{B_{3}}{10}.\label{3.026}
\end{align}
Combining (\ref{3.025}) with  (\ref{3.026}), we have
(\ref{3.021}) is valid.

 When $\xi_{1}\xi_{2}\leq0$,
we have
\begin{align}
\sqrt{\xi^{2}(\xi^{2}-\xi_{1}\xi_{2})\pm B_{2}}+\sqrt{\xi^{2}(\xi^{2}-\xi_{1}\xi_{2})}\geq |\xi||\xi_{1}\xi_{2}|^{\frac12}.\label{3.027}
\end{align}
 Since
\begin{align}
&\hspace{-2cm}\left[\xi_{1}^{2}+\xi_{2}^{2}+\sqrt{\xi^{2}(\xi^{2}-\xi_{1}\xi_{2})\pm B_{2}}\right]
\left[\xi_{1}^{2}+\xi_{2}^{2}+\sqrt{\xi^{2}(\xi^{2}-\xi_{1}\xi_{2})}\right]\geq \frac{B_{1}}{4},\label{3.028}
\end{align}
 by  a direct computation we have
\begin{align}
& \left|\frac{B_{1}}
{\xi_{1}^{2}+\xi_{2}^{2}+\sqrt{\xi^{2}(\xi^{2}-\xi_{1}\xi_{2})}}-\frac{B_{1}}{\xi_{1}^{2}
+\xi_{2}^{2}+\sqrt{\xi^{2}(\xi^{2}-\xi_{1}\xi_{2})\pm B_{2}}}\right|\nonumber\\
=& \frac{B_{1}B_{2}}{\left[\xi_{1}^{2}+\xi_{2}^{2}+\sqrt{\xi^{2}(\xi^{2}-\xi_{1}\xi_{2})\pm B_{2}}\right]
\left[\xi_{1}^{2}+\xi_{2}^{2}+\sqrt{\xi^{2}(\xi^{2}-\xi_{1}\xi_{2})}\right]}\nonumber\\   & \qquad \times
\frac{1}{\left[\sqrt{\xi^{2}(\xi^{2}-\xi_{1}\xi_{2})\pm B_{2}}+\sqrt{\xi^{2}
(\xi^{2}-\xi_{1}\xi_{2})}\right]}\nonumber\\
\leq& \frac{B_{3}}{10}.\label{3.029}
\end{align}
By  a direct computation,  we have
\begin{align}
\frac{B_{2}}{\xi_{1}^{2}+\xi_{2}^{2}+\sqrt{\xi^{2}(\xi^{2}-\xi_{1}\xi_{2})\pm B_{2}}}\leq \frac{B_{3}}{10}.\label{3.030}
\end{align}
By  (\ref{3.023})-(\ref{3.030}), we see that (\ref{3.021})  is valid.

To obtain (\ref{3.018}), we make the change of variable $\beta_{1}=\beta_{2}+\beta$.
Thus, (\ref{3.011})-(\ref{3.012}) can be rewritten as follows:
\begin{align}
&\beta\geq-\sqrt{5}(\xi_{1}^{2}+\xi_{2}^{2}),\nonumber\\   &\left|\beta^{2}+2\sqrt{5}\beta(\xi_{1}^{2}+\xi_{2}^{2})+
5(4\xi_{1}^{2}\xi_{2}^{2}-3\xi^{2}\xi_{1}\xi_{2})\right|\leq 2^{j+k-k_{1}-k_{2}+3}.\label{3.031}
\end{align}
Since $|\xi_{j}|\in [2^{k_{j}-1},2^{k_{j}+1}](j=1,2)$, we can assume that  $\xi_{j}=a_{j}2^{k_{j}}(j=1,2),$ where
$\frac{1}{2}\leq |a_{j}|\leq 2.$ Consequently, (\ref{3.021}) can be rewritten as follows:
\begin{align}
\left|\beta-\sqrt{5}f_{1}(a_{1},a_{2},k_{1},k_{2})\right|\leq B_{3},\label{3.032}
\end{align}
where
\begin{align}
 &  f_{1}(a_{1},a_{2},k_{1},k_{2})\nonumber\\
=&  \frac{a_{1}a_{2}2^{k_{1}+k_{2}}
(3a_{1}^{2}4^{k_{1}}+a_{1}a_{2}2^{k_{1}+k_{2}+1}+3a_{2}^{2}4^{k_{2}})}
{a_{1}^{2}4^{k_{1}}+a_{2}^{2}4^{k_{2}}+\sqrt{(a_{1}2^{k_{1}}+a_{2}2^{k_{2}})^{2}
(a_{1}^{2}4^{k_{1}}+a_{1}a_{2}2^{k_{1}+k_{2}}+a_{2}^{2}4^{k_{2}})}}\nonumber\\
=&\frac{B_{1}}{\xi_{1}^{2}+\xi_{2}^{2}+\sqrt{\xi^{2}(\xi^{2}-\xi_{1}\xi_{2})}}.\label{3.033}
\end{align}
Thus, the left hand side of (\ref{3.018}) can be bounded by
\begin{align}
&2^{\frac{k_{1}+k_{2}}{2}}\int_{\tilde{S}}h_{1}(\xi_{1},\beta+\beta_{2})h_{2}(\xi_{2},\beta_{2})
\chi_{[-1,1)}\biggl(\frac{\beta-\sqrt{5}f_{1}(a_{1},a_{2},k_{1},k_{2})}{B_{3}}-m\biggr)
\nonumber\\
&\qquad \times
h(\xi_{1}+\xi_{2},A(\xi_{1},\xi_{2},\beta)+\beta_{2}(\xi_{1}+\xi_{2}),B(\xi_{1},\xi_{2},\beta))
 d\xi_{1}d\xi_{2}d\beta d\beta_{2},\label{3.034}
\end{align}
where
\begin{align}
&\tilde{S}=\left\{(\xi_{1},\xi_{2},\beta,\beta_{2})\in \R^{4},\xi_{1}+\xi_{2}\geq 0,\quad
\beta \quad {\rm satisfies }\quad (3.21) \right\},\label{3.035}\\
&A(\xi_{1},\xi_{2},\beta)=\sqrt{5}\left[\xi_{1}^{3}-\xi_{2}^{3}\right]+\beta\xi_{1},\label{3.036}
\\   &B(\xi_{1},\xi_{2},\beta)
=\frac{\xi_{1}\xi_{2}}{\xi}\left[\beta^{2}+2\sqrt{5}\beta(\xi_{1}^{2}+\xi_{2}^{2})
+5\xi_{1}\xi_{2}(4\xi_{1}\xi_{2}-3\xi^{2})\right].\label{3.037}
\end{align}
Let $j^{\prime}=j-\frac{3(k_{1}+k_{2})}{2}+10$.
Decompose
\begin{align}
   h_{i}(\xi^{\prime},\beta^{\prime})\nonumber
=& \sum_{m\in \mathbf{Z}}h_{i} (\xi^{\prime},\beta^{\prime})
   \chi_{[0,1)}\biggl(\frac{\beta^{\prime}
   -\sqrt{5}f_{1}(a_{1},a_{2},k_{1},k_{2})}{2^{j^{\prime}}}-m\biggr)\\
=& \sum\limits_{m\in \mathbf{Z}}h_{i}^{m}(\xi^{\prime},\beta^{\prime}),
  \quad i=1,2\label{3.038}
\end{align}
 for all $a_{j}\in R,\frac{1}{2}\leq |a_{j}|\leq2(j=1,2).$ Obviously,
 if $m_{1},m_{2}\in \mathbf{Z}, m_{1}\neq m_{2},$ then
 \begin{align*}
 \prod\limits_{i=1}^{2}\chi_{[0,1)}
 \biggl(\frac{\beta^{\prime}
-\sqrt{5}f_{1}(a_{1},a_{2},k_{1},k_{2})}{2^{j^{\prime}}}-m_{i}\biggr)=0.
 \end{align*}
 Thus, for  $m_{1},m_{2}\in \mathbf{Z}, m_{1}\neq m_{2},$ we have  $\prod\limits_{i=1}^{2}h_{i}^{m_{i}}(\xi^{\prime},\beta^{\prime})=0.$
Consequently, we have
 \begin{align}
\|h_{i}(\xi^{\prime},\beta^{\prime})\|_{L_{\xi^{\prime}\beta^{\prime}}^{2}}= \biggl\|\sum\limits_{m\in \mathbf{Z}}h_{i}^{m}(\xi^{\prime},\beta^{\prime})\biggr\|_{L_{\xi^{\prime}\beta^{\prime}}^{2}}
 =\biggl(\sum\limits_{m\in \mathbf{Z}}\|h_{i}^{m}\|_{L_{\xi^{\prime}\beta^{\prime}}^{2}}^{2}\biggr)^{\frac12}.\label{3.039}
 \end{align}
Thus, (\ref{3.034}) is controlled by
\begin{align}
& 2^{\frac{k_{1}+k_{2}}{2}}\sum\limits_{|m-m^{\prime}|\leq 4}\int_{\tilde{S}}h_{1}^{m}(\xi_{1},
\beta+\beta_{2})h_{2}^{m^{\prime}}(\xi_{2},\beta_{2})
\nonumber\\
&\qquad \qquad\times
h(\xi_{1}+\xi_{2},A(\xi_{1},\xi_{2},\beta)+\beta_{2}(\xi_{1}+\xi_{2}),B(\xi_{1},\xi_{2},\beta))
 d\xi_{1}d\xi_{2}d\beta d\beta_{2}.\label{3.040}
\end{align}
To prove (\ref{3.018}),  it suffices to prove that
\begin{align}
&2^{\frac{k_{1}+k_{2}}{2}}\int_{\tilde{S}}h_{1}^{m}(\xi_{1},\beta+\beta_{2})h_{2}^{m^{\prime}}(\xi_{2},\beta_{2})
\chi_{[m^{\prime},m^{\prime}+1)}\left(\frac{\beta_{2}-\sqrt{5}f_{1}(a_{1},a_{2},k_{1},k_{2})}
{2^{j^{\prime}}}\right)\nonumber\\   &\qquad\qquad\times
h(\xi_{1}+\xi_{2},A(\xi_{1},\xi_{2},\beta)+\beta_{2}(\xi_{1}+\xi_{2}),B(\xi_{1},\xi_{2},\beta))
 d\xi_{1}d\xi_{2}d\beta d\beta_{2}\nonumber\\
 \leq & C2^{\frac{j}{2}}2^{-\frac{7}{4}(k_{1}+k_{2})+\frac{k_{3}}{2}}\|h\|_{L^{2}}\|h_{1}^{m}\|_{L^{2}}\|h_{2}^{m^{\prime}}\|_{L^{2}}.\label{3.041}
\end{align}
If (\ref{3.041}) is valid, by using the Cauchy-Schwartz  inequality,  we have
\begin{align}
& 2^{\frac{k_{1}+k_{2}}{2}}
   \sum\limits_{|m-m^{\prime}|\leq 4}\int_{\tilde{S}}h_{1}^{m}(\xi_{1},
   \beta+\beta_{2})h_{2}^{m^{\prime}}(\xi_{2},\beta_{2})
\nonumber\\
&\qquad\qquad\times
   h(\xi_{1}+\xi_{2},A(\xi_{1},\xi_{2},\beta)+\beta_{2}(\xi_{1}+\xi_{2}),B(\xi_{1},\xi_{2},\beta))
    d\xi_{1}d\xi_{2}d\beta d\beta_{2}\nonumber\\
\leq & C2^{\frac{j}{2}}2^{-\frac{7}{4}(k_{1}+k_{2})+\frac{k_{3}}{2}}
   \|h\|_{L^{2}}
    \sum\limits_{|m-m^{\prime}|\leq  4}\|h_{1}^{m}\|_{L^{2}}
    \|h_{2}^{m^{\prime}}\|_{L^{2}}
    \nonumber\\
 =&C2^{\frac{j}{2}}2^{-\frac{7}{4}(k_{1}+k_{2})+\frac{k_{3}}{2}}\|h\|_{L^{2}}\left[\sum\limits_{m\in \mathbf{Z}}\|h_{1}^{m}\|_{L^{2}}\left(\sum\limits_{m-4\leq m^{\prime}\leq m+4}\|h_{2}^{m^{\prime}}\|_{L^{2}}\right)\right]\nonumber\\
 \leq & C2^{\frac{j}{2}}2^{-\frac{7}{4}(k_{1}+k_{2})+\frac{k_{3}}{2}}\|h\|_{L^{2}}\left[\sum\limits_{m\in \mathbf{Z}}\|h_{1}^{m}\|_{L^{2}}^{2}\right]^{\frac{1}{2}}\left[\sum\limits_{m\in \mathbf{Z}}\left(\sum\limits_{m-4\leq m^{\prime}\leq m+4}\|h_{2}^{m^{\prime}}\|_{L^{2}}\right)^{2}\right]^{\frac{1}{2}}\nonumber\\
 \leq & C2^{\frac{j}{2}}2^{-\frac{7}{4}(k_{1}+k_{2})+\frac{k_{3}}{2}}\|h\|_{L^{2}}\|h_{1}\|_{L^{2}}
 \left(\sum_{m\in \mathbf{Z}}\sum_{m-4\leq m^{\prime}\leq m+4}\|h_{2}^{m^{\prime}}\|^2_{L^{2}}\right)^{\frac{1}{2}}\nonumber\\
 \leq & C2^{\frac{j}{2}}2^{-\frac{7}{4}(k_{1}+k_{2})+\frac{k_{3}}{2}}\|h\|_{L^{2}}
 \prod\limits_{i=1}^{2}\|h_{i}\|_{L^{2}}.\label{3.042}
\end{align}
To prove (\ref{3.041}), without loss of generality,  we assume that $|\xi_{1}|\leq |\xi_{2}|$.
To prove (\ref{3.041}), by using the Minkowski inequality with respect to variables $(\xi_{1},\xi_{2},\beta)$ with
\begin{align}
&S^{\prime}=\left\{(\xi_{1},\xi_{2},\beta)\in \R^{3},\xi_{1}+\xi_{2}\geq 0,\quad \beta \quad {\rm satisfies} \quad (3.21) \right\}
,\label{3.043}
\end{align}
 the left-hand side of (\ref{3.040}) is controlled by
\begin{align}
& 2^{\frac{k_{1}+k_{2}}{2}}\int_{\SR}
\chi_{[m,m+1)}\left(\frac{\beta_{2}-\sqrt{5}f_{1}(a_{1},a_{2},k_{1},k_{2})}{2^{j^{\prime}}}\right)
\nonumber\\
&\quad \times\left(\int_{S^{\prime}}|h_{1}^{m}(\xi_{1},\beta+\beta_{2})
h_{2}^{m^{\prime}}(\xi_{2},\beta_{2})|^{2}d\xi_{1}d\xi_{2}d\beta\right)^{\frac12}\nonumber\\   &\quad \times
\left(\int_{S^{\prime}}|h(\xi_{1}+\xi_{2},A(\xi_{1},\xi_{2},\beta)
+\beta_{2}(\xi_{1}+\xi_{2}),B(\xi_{1},\xi_{2},\beta))|^{2}
 d\xi_{1}d\xi_{2}d\beta\right)^{\frac12} d\beta_{2}.\label{3.044}
\end{align}
From (\ref{3.044}), it suffices to prove that
\begin{align}
&\left(\int_{S^{\prime}}|h(\xi_{1}+\xi_{2},A(\xi_{1},\xi_{2},\beta)
+\beta_{2}(\xi_{1}+\xi_{2}),B(\xi_{1},\xi_{2},\beta))|^{2}
 d\xi_{1}d\xi_{2}d\beta\right)^{\frac12}\nonumber\\
 \leq & C2^{-\frac{3(k_{1}+k_{2})}{2}+\frac{k_{3}}{2}}\|h\|_{L^{2}}.\label{3.045}
\end{align}
If (\ref{3.045}) is valid, by using the Cauchy-Schwartz inequality
with  respect to $\beta_{2}$, we have  (\ref{3.044})  is controlled by
\begin{align}
&   C2^{-(k_{1}+k_{2})+\frac{k_{3}}{2}}
   \int_{\SR}\chi_{[m^{\prime},m^{\prime}+1)}
   \bigg(\frac{\beta_{2}-\sqrt{5}f_{1}(a_{1},a_{2},k_{1},k_{2})}{2^{j^{\prime}}}\bigg)
   \nonumber\\
&\qquad\qquad\qquad\qquad\times\|h_{1}^{m}\|_{L^{2}}
   \|h_{2}^{m^{\prime}}(\cdot,\beta_{2})\|_{L_{\xi_{2}}^{2}}
   \|h\|_{L^{2}}d\beta_{2}\nonumber\\
=&  C2^{-(k_{1}+k_{2})+\frac{k_{3}}{2}}
   \int_{m^{\prime}2^{j^{\prime}}
     +\sqrt{5}f_{1}(a_{1},a_{2},k_{1},k_{2})}^{(m^{\prime}+1)2^{j^{\prime}}
      +\sqrt{5}f_{1}(a_{1},a_{2},k_{1},k_{2})}
       \|h_{1}^{m}\|_{L^{2}}\|h_{2}^{m^{\prime}}(\cdot,\beta_{2})\|_{L_{\xi_{2}}^{2}}
        \|h\|_{L^{2}}d\beta_{2}\nonumber\\
\leq&  C2^{-(k_{1}+k_{2})+\frac{k_{3}}{2}}
       2^{\frac{j^{\prime}}{2}}\|h\|_{L^{2}}
         \|h_{1}^{m}\|_{L^{2}}\|h_{2}^{m^{\prime}}\|_{L^{2}}
          \nonumber\\
\leq&  C2^{\frac{j}{2}}2^{-\frac{7}{4}(k_{1}+k_{2})
       +\frac{k_{3}}{2}}\|h\|_{L^{2}}
         \|h_{1}^{m}\|_{L^{2}}
          \|h_{2}^{m^{\prime}}\|_{L^{2}}.\label{3.046}
\end{align}
To prove (\ref{3.045}), we may assume that $\beta_{2}=0$ and make the change of
 variable $\beta=\sqrt{5}\xi_{1}\xi_{2}\nu$.
 From   (\ref{3.021}), we have
\begin{align}
\left|\nu-\frac{3\xi_{1}^{2}+2\xi_{1}\xi_{2}+3\xi_{2}^{2}}{\xi_{1}^{2}+\xi_{2}^{2}
+\sqrt{\xi^{2}(\xi^{2}-\xi_{1}\xi_{2})}}\right|\leq 2^{-20}.\label{3.047}
\end{align}
The left hand side of (\ref{3.044}) is controlled by
\begin{align}
&\hspace{-1.7cm}C2^{\frac{k_{1}+k_{2}}{2}}\left(\int_{S^{\prime\prime}}
|h(\xi_{1}+\xi_{2},H_{1}(\xi_{1},\xi_{2},\nu),
H_{2}(\xi_{1},\xi_{2},\nu))d\xi_{1}d\xi_{2}d\nu\right)^{\frac12},\label{3.048}
\end{align}
where
\begin{align}
& S^{\prime\prime}=\left\{(\xi_{1},\xi_{2},\nu)\in \R^{3}:\xi_{i}\in \tilde{I}_{k_{i}},
\nu\ \ {\rm satisfies} \ \ (3.45)\right\},\label{3.049}\\
&H_{1}(\xi_{1},\xi_{2},\nu)=\sqrt{5}(\xi_{1}^{3}-\xi_{2}^{3}+\nu\xi_{1}^{2}\xi_{2}), \label{3.050}\\   &
H_{2}(\xi_{1},\xi_{2},\nu)=\frac{5(\xi_{1}\xi_{2})^{2}}{\xi}\left(\xi_{1}\xi_{2}\nu^{2}
+2\nu(\xi_{1}^{2}+\xi_{2}^{2})+4\xi_{1}\xi_{2}-3\xi^{2}\right).\label{3.051}
\end{align}
We consider $\xi_{1}\xi_{2}\geq0,\xi_{1}\xi_{2}\leq0,$ respectively.

  Firstly, we consider $\xi_{1}\xi_{2}\geq0.$
  We define the function
\begin{align}
    G(\xi,x,y)
=   2^{-k_{3}+2k_{1}+2k_{2}}\left|h\bigg(\xi,\sqrt{5}\left[\xi_{1}^{3}-\xi_{2}^{3}+x\right],
\frac{5(\xi_{1}\xi_{2})^{2}}{\xi}\left[y+4\xi_{1}\xi_{2}-3\xi^{2}\right]\bigg)\right|^{2},\label{3.052}
\end{align}
where
\begin{align}
x=\xi_{1}^{2}\xi_{2}\nu,\quad
y=\xi_{1}\xi_{2}\nu^{2}+2(\xi_{1}^{2}+\xi_{2}^{2})\nu.\label{3.053}
\end{align}
Obviously, $\|G\|_{L^{1}}=\|h\|_{L^{2}}^{2}.$
From (\ref{3.052}), we have
(\ref{3.048}) can be bounded by
\begin{align}
C2^{-\frac{k_{1}+k_{2}-k_{3}}{2}}\left(\int_{S^{\prime\prime}}|G(\xi_{1}+\xi_{2},\xi_{1}^{2}\xi_{2}\nu,
\xi_{1}\xi_{2}\nu^{2}+2(\xi_{1}^{2}+\xi_{2}^{2})\nu)|
 d\xi_{1}d\xi_{2}d\nu\right)^{\frac12}.\label{3.054}
\end{align}
We make the change of variables $(\xi_{1},\xi_{2},\nu)\longrightarrow
 (\xi_{1}+\xi_{2},\xi_{1}^{2}\xi_{2}\nu,
\xi_{1}\xi_{2}\nu^{2}+2(\xi_{1}^{2}+\xi_{2}^{2})\nu)$, thus the absolute
 value of the Jacobi determinant equals
\begin{align}
|\nu\xi_{1}|\left|\nu\xi_{1}\xi_{2}(\xi_{1}-3\xi_{2})+2\left(\xi_{1}^{3}
-\xi_{1}\xi_{2}^{2}-2\xi_{2}^{3}\right)\right|.\label{3.055}
\end{align}
By using a direct computation, we have  (\ref{3.055}) equals
\begin{align}
|\nu\xi_{1}^{2}\xi_{2}(\xi_{1}-3\xi_{2})|\left|\nu+\frac{2(\xi_{1}^{3}-\xi_{1}\xi_{2}^{2}
-2\xi_{2}^{3})}{\xi_{1}\xi_{2}(\xi_{1}-3\xi_{2})}\right|,\label{3.056}
\end{align}
where  $ \nu $ satisfies (\ref{3.047}).

\noindent By using a direct computation, we have
\begin{align}
\left|\frac{2(\xi_{1}^{3}-\xi_{1}\xi_{2}^{2}-2\xi_{2}^{3})}
{\xi_{1}\xi_{2}(\xi_{1}-3\xi_{2})}\right|\geq 2.\label{3.057}
\end{align}
 From (\ref{3.047}),  we have
\begin{align}
&1-2^{-20}\leq |\nu|\leq \left|\frac{3\xi_{1}^{2}+2\xi_{1}\xi_{2}+3\xi_{2}^{2}}{\xi_{1}^{2}+\xi_{2}^{2}
+\sqrt{\xi^{2}(\xi^{2}-\xi_{1}\xi_{2})}}\right|+2^{-20}\leq \frac{3}{2}+2^{-20}.\label{3.058}
\end{align}
Combining (\ref{3.057})  with (\ref{3.058}), we have
\begin{align}
|\nu|\left|\nu+\frac{2(\xi_{1}^{3}-\xi_{1}\xi_{2}^{2}-2\xi_{2}^{3})}
{\xi_{1}\xi_{2}(\xi_{1}-3\xi_{2})}\right|\geq \frac{1}{4}.\label{3.059}
\end{align}
Combining  (\ref{3.056}) with  (\ref{3.059}), we have
\begin{align}
|\nu\xi_{1}|\left|\nu\xi_{1}\xi_{2}(\xi_{1}-3\xi_{2})+2\left(\xi_{1}^{3}
-\xi_{1}\xi_{2}^{2}-2\xi_{2}^{3}\right)\right|\geq C\xi_{1}^{2}\xi_{2}^{2}.\label{3.060}
\end{align}
Combining (\ref{3.054}) with (\ref{3.060}), we have
(\ref{3.045}) can be bounded by
\begin{align}
C2^{-\frac{3(k_{1}+k_{2})}{2}+\frac{k_{3}}{2}}\|G\|_{L^{1}}\leq C
2^{-\frac{3(k_{1}+k_{2})}{2}+\frac{k_{3}}{2}}\|h\|_{L^{2}}^{2}.\label{3.061}
\end{align}

 Now we consider $\xi_{1}\xi_{2}\leq0.$
We define
\begin{align}
   G(\xi,x,y)
=   2^{-k_{3}+7k_{1}}\left|h\bigg(\xi,\sqrt{5}\left[\xi_{1}^{3}-\xi_{2}^{3}+\xi_{1}^{2}\xi_{2}x\right],
\frac{5(\xi_{1}\xi_{2})^{2}}{\xi}\left[y+4\xi_{1}\xi_{2}-3\xi^{2}\right]\bigg)\right|^{2},\label{3.062}
\end{align}
where
\begin{align}
x=\nu,\quad y=\xi_{1}\xi_{2}\nu^{2}+2(\xi_{1}^{2}+\xi_{2}^{2})\nu.\label{3.063}
\end{align}
Obviously, $\|G\|_{L^{1}}\approx\|h\|_{L^{2}}^{2}.$
From (\ref{3.062}), we have
(\ref{3.048}) can be bounded by
\begin{align}
C2^{-\frac{5k_{1}-k_{3}}{2}}\left(\int_{S^{\prime\prime}}|G(\xi_{1}+\xi_{2},\nu,
\xi_{1}\xi_{2}\nu^{2}+2(\xi_{1}^{2}+\xi_{2}^{2})\nu)|
 d\xi_{1}d\xi_{2}d\nu\right)^{\frac12}.\label{3.064}
\end{align}
 We make the change of variables $(\xi_{1},\xi_{2},\nu)\longrightarrow
 (\xi_{1}+\xi_{2},\nu,
\xi_{1}\xi_{2}\nu^{2}+2(\xi_{1}^{2}+\xi_{2}^{2})\nu)$, thus the absolute
 value of the Jacobi determinant equals
\begin{align}
|\xi_{1}-\xi_{2}||\nu(4-\nu)|.\label{3.065}
\end{align}
From (\ref{3.047}),  we have
\begin{align}
1-2^{-20}
&\leq\frac{3\xi_{1}^{2}+2\xi_{1}\xi_{2}+3\xi_{2}^{2}}{\xi_{1}^{2}+\xi_{2}^{2}
+\sqrt{\xi^{2}(\xi^{2}-\xi_{1}\xi_{2})}}-2^{-20}\nonumber\\   &\leq |\nu|\leq \frac{3\xi_{1}^{2}+2\xi_{1}\xi_{2}+3\xi_{2}^{2}}{\xi_{1}^{2}+\xi_{2}^{2}
+\sqrt{\xi^{2}(\xi^{2}-\xi_{1}\xi_{2})}}+2^{-20}\leq 3+2^{-20}.\label{3.066}
\end{align}
Combining (\ref{3.065})  with (\ref{3.066}), we have
\begin{align}
|\xi_{1}-\xi_{2}||\nu(4-\nu)|\sim 2^{k_{1}}.\label{3.067}
\end{align}
Combining (\ref{3.064}) with (\ref{3.067}), we have
(\ref{3.044}) can be bounded by
\begin{align}
C2^{-\frac{3(k_{1}+k_{2})}{2}+\frac{k_{3}}{2}}\|G\|_{L^{1}}\leq C
2^{-\frac{3(k_{1}+k_{2})}{2}+\frac{k_{3}}{2}}\|h\|_{L^{2}}^{2}.\label{3.068}
\end{align}
Therefore the proof of \eqref{3.03} is completed.

Now we prove \eqref{3.04}. From (\ref{3.05})-(\ref{3.017}), we know that
it suffices to prove
\begin{align}
&2^{\frac{k_{1}+k_{2}}{2}}\int_{\tilde{S}}h_{1}(\xi_{1},\beta_{1})h_{2}(\xi_{2},\beta_{2})
\nonumber\\
&\qquad\qquad\times
h(\xi_{1}+\xi_{2},\sqrt{5}(\xi_{1}^{3}-\xi_{2}^{3})+\beta_{1}\xi_{1}+\beta_{2}\xi_{2},
\tilde{R}_{1}(\xi_{1},\beta_{1},\xi_{2},\beta_{2})d\xi_{1}d\xi_{2}d\beta_{1}d\beta_{2}\nonumber\\
\leq& C2^{\frac{j}{2}}2^{-\frac{k_{1}}{2}-\frac{5k_{2}}{2}}
\|h\|_{L^{2}}
\prod\limits_{j=1}^{2}\|h_{j}\|_{L^{2}} .\label{3.069}
\end{align}
where $h_{i}(i=1,2)$ are defined as in (\ref{3.016})-(\ref{3.017}) and $S$ is defined as in (\ref{3.014})
 and $\tilde{R}_{1}(\xi_{1},\beta_{1},\xi_{2},\beta_{2})$ is defined as in (\ref{3.015})
 and $\beta_{1}-\beta_{2}$ satisfies (\ref{3.019})-(\ref{3.020}).
To obtain (\ref{3.069}),  we make the change of variable $\beta_{1}=\beta_{2}+\beta$.
Now we claim that the following inequality is valid
\begin{align}
\left|\beta-\frac{\sqrt{5}B_{1}}{\xi_{1}^{2}+\xi_{2}^{2}+
\sqrt{\xi^{2}(\xi^{2}-\xi_{1}\xi_{2})}}\right|\leq B_{4},\label{3.070}
\end{align}
where $B_{1}$ is defined as in (\ref{3.021}) and
\begin{align}
B_{4}=2^{j-3k_{2}+10}.\label{3.071}
\end{align}
(\ref{3.070}) can be proved similarly to (\ref{3.021}).
Since $|\xi_{j}|\in [2^{k_{j}-1},2^{k_{j}+1}](j=1,2)$, we can assume that  $\xi_{j}=a_{j}2^{k_{j}}(j=1,2),$ where
$\frac{1}{2}\leq |a_{j}|\leq 2.$ Consequently, (\ref{3.070}) can be rewritten as follows:
\begin{align}
\left|\beta-\sqrt{5}f_{1}(a_{1},a_{2},k_{1},k_{2})\right|\leq B_{4},\label{3.072}
\end{align}
where $f_{1}(a_{1},a_{2},k_{1},k_{2})$ is defined as in (\ref{3.033}).

\noindent Thus, the left hand side of (\ref{3.068}) can be bounded by
\begin{align}
&2^{\frac{k_{1}+k_{2}}{2}}\int_{\tilde{S}}h_{1}(\xi_{1},\beta+\beta_{2})h_{2}(\xi_{2},\beta_{2})
\chi_{[-1,1)}\bigg(\frac{\beta_{2}-\sqrt{5}f_{1}(a_{1},a_{2},k_{1},k_{2})}{B_{4}}-m\bigg)
\nonumber\\
&\qquad\qquad\times
h(\xi_{1}+\xi_{2},A(\xi_{1},\xi_{2},\beta)+\beta_{2}(\xi_{1}+\xi_{2}),B(\xi_{1},\xi_{2},\beta))
 d\xi_{1}d\xi_{2}d\beta d\beta_{2},\label{3.073}
\end{align}
where
$\tilde{S}$ satisfies  (\ref{3.035}) and $A(\xi_{1},\xi_{2},\beta)$ satisfies  (\ref{3.036}) and
  $B(\xi_{1},\xi_{2},\beta)$ satisfies  (\ref{3.037}).
Let $j^{\prime\prime}=j-3k_{2}+10$.
Decompose
\begin{align}
 h_{i}(\xi^{\prime},\beta^{\prime})
=   &\sum\limits_{m\in \mathbf{Z}}h_{i}
(\xi^{\prime},\beta^{\prime})\chi_{[0,1)}\bigg(\frac{\beta^{\prime}
-\sqrt{5}f_{1}(a_{1},a_{2},k_{1},k_{2})}
{2^{j^{\prime\prime}}}-m\bigg)\nonumber\\
=&\sum\limits_{m\in \mathbf{Z}}h_{i}^{m}(\xi^{\prime},\beta^{\prime}),\label{3.074}
\end{align}
 $i=1,2$.
Obviously,
$\|h_{i}\|_{L^{2}}
=\bigg(\sum\limits_{m\in \mathbf{Z}}\|h_{i}^{m}\|_{L^{2}}^{2}\bigg)^{\frac12}.$
Thus, (\ref{3.073}) is controlled by
\begin{align}
& \quad
2^{\frac{k_{1}+k_{2}}{2}}
  \sum\limits_{|m-m^{\prime}|\leq 4}\int_{\tilde{S}}
     h_{1}^{m}(\xi_{1},\beta+\beta_{2})
      h_{2}^{m^{\prime}}(\xi_{2},\beta_{2})
      \nonumber\\
&\qquad\qquad\times
  h(\xi_{1}+\xi_{2},A(\xi_{1},\xi_{2},\beta)+\beta_{2}(\xi_{1}+\xi_{2}),
    B(\xi_{1},\xi_{2},\beta))
     d\xi_{1}d\xi_{2}d\beta d\beta_{2}.\label{3.075}
\end{align}
To prove (\ref{3.069}),  it suffices to prove that
\begin{align}
& 2^{\frac{k_{1}+k_{2}}{2}}\int_{\tilde{S}}h_{1}^{m}(\xi_{1},\beta+\beta_{2})h_{2}^{m^{\prime}}(\xi_{2},\beta_{2})
\chi_{[m^{\prime},m^{\prime}+1)}\left(\frac{\beta_{2}-\sqrt{5}f_{1}(a_{1},a_{2},k_{1},k_{2})}
{2^{j^{\prime\prime}}}\right)\nonumber\\   &\qquad\qquad\times
h(\xi_{1}+\xi_{2},A(\xi_{1},\xi_{2},\beta)+\beta_{2}(\xi_{1}+\xi_{2}),B(\xi_{1},\xi_{2},\beta))
 d\xi_{1}d\xi_{2}d\beta d\beta_{2}\nonumber\\
\leq &
 C2^{\frac{j}{2}}2^{-\frac{k_{1}}{2}-\frac{5k_{2}}{2}}\|h\|_{L^{2}}
\|h_{1}^{m}\|_{L^{2}}\|h_{2}^{m^{\prime}}\|_{L^{2}}.\label{3.076}
\end{align}
From (\ref{3.076}), it suffices to prove that
\begin{align}
&\left(\int_{S^{\prime}}|h(\xi_{1}+\xi_{2},A(\xi_{1},\xi_{2},\beta)
+\beta_{2}(\xi_{1}+\xi_{2}),B(\xi_{1},\xi_{2},\beta))|^{2}
 d\xi_{1}d\xi_{2}d\beta\right)^{\frac12}\nonumber\\
\leq & C
 2^{-\frac{3k_{2}}{2}-k_{1}}\|h\|_{L^{2}}.\label{3.077}
\end{align}
If (\ref{3.077}) is valid, by using the Cauchy-Schwartz inequality  with  respect to $\beta_{2}$,
we have  (\ref{3.075})  is controlled by
\begin{align}
&\quad C2^{-\frac{k_{1}}{2}-k_{2}}\int_{\SR}\chi_{[m^{\prime},m^{\prime}+1)}
\left(\frac{\beta_{2}-\sqrt{5}f_{1}(a_{1},a_{2},k_{1},k_{2})}{2^{j^{\prime\prime}}}\right)
\nonumber\\   &\qquad\qquad\qquad\qquad\times\|h_{1}^{m}\|_{L^{2}}\|h_{2}^{m^{\prime}}(\cdot,\beta_{2})\|_{L_{\xi_{2}}^{2}}\|h\|_{L^{2}}d\beta_{2}\nonumber\\
&=C2^{-\frac{k_{1}}{2}-k_{2}}\int_{m^{\prime}2^{j^{\prime\prime}}
-\sqrt{5}f_{1}(a_{1},a_{2},k_{1},k_{2})}^{(m^{\prime}+1)2^{j^{\prime\prime}}
-\sqrt{5}f_{1}(a_{1},a_{2},k_{1},k_{2})}
\|h_{1}^{m}\|_{L^{2}}\|h_{2}^{m^{\prime}}(\cdot,\beta_{2})\|_{L_{\xi_{2}}^{2}}\|h\|_{L^{2}}d\beta_{2}\nonumber\\
&\leq C2^{-\frac{k_{1}}{2}-k_{2}}2^{\frac{j^{\prime\prime}}{2}}\|h\|_{L^{2}}
\|h_{1}^{m}\|_{L^{2}}\|h_{2}^{m^{\prime}}\|_{L^{2}}\nonumber\\   &\leq C2^{-\frac{5k_{2}}{2}-\frac{k_{1}}{2}}2^{\frac{j}{2}}\|h\|_{L^{2}}
\|h_{1}^{m}\|_{L^{2}}\|h_{2}^{m^{\prime}}\|_{L^{2}}.\label{3.078}
\end{align}
To prove (\ref{3.077}), we may assume that $\beta_{2}=0$ and make the change of
 variable $\beta=\sqrt{5}\xi_{1}\xi_{2}\nu$.
The left hand side of (\ref{3.077}) is controlled by
\begin{align}
&\hspace{-1.7cm}C2^{\frac{k_{1}+k_{2}}{2}}\left(\int_{S^{\prime\prime}}
|h(\xi_{1}+\xi_{2},H_{1}(\xi_{1},\xi_{2},\nu),H_{2}(\xi_{1},\xi_{2},\nu))|^{2}d\xi_{1}d\xi_{2}d\nu\right)^{\frac12},\label{3.079}
\end{align}
where
$S^{\prime\prime}$ is defined as in (\ref{3.049}) and
$H_{1}(\xi_{1},\xi_{2},\nu)$ is defined as in (\ref{3.050}) and
$H_{2}(\xi_{1},\xi_{2},\nu)$ is defined as in (\ref{3.051}).

\noindent We define the function
$G(\xi,x,y)$ as in (\ref{3.052})
and $x,y$ as in (\ref{3.053}).
Obviously, $\|G\|_{L^{1}}\approx\|h\|_{L^{2}}^{2}.$
Thus, we have
(\ref{3.079}) can be bounded by
\begin{align}
C2^{-\frac{k_{1}}{2}}\left(\int_{S^{\prime\prime}}|G(\xi_{1}+\xi_{2},\xi_{1}^{2}\xi_{2}\nu,
\xi_{1}\xi_{2}\nu^{2}+2(\xi_{1}^{2}+\xi_{2}^{2})\nu)|
 d\xi_{1}d\xi_{2}d\nu\right)^{\frac12}.\label{3.080}
\end{align}
We make the change of variables $(\xi_{1},\xi_{2},\nu)\longrightarrow (\xi_{1}+\xi_{2},\xi_{1}^{2}\xi_{2}\nu,
\xi_{1}\xi_{2}\nu^{2}+2(\xi_{1}^{2}+\xi_{2}^{2})\nu)$, thus the absolute value of the Jacobi determinant equals
\begin{align}
2|\nu\xi_{1}\left(\xi_{1}^{3}-\xi_{1}\xi_{2}^{2}-2\xi_{2}^{3}\right)|\left|\frac{\xi_{1}\xi_{2}(\xi_{1}-3\xi_{2})\nu}
{2\left(\xi_{1}^{3}-\xi_{1}\xi_{2}^{2}-2\xi_{2}^{3}\right)}+1\right|
.\label{3.081}
\end{align}
In this case,  by using a direct computation, we have
\begin{align}
\left|\frac{\xi_{1}\xi_{2}(\xi_{1}-3\xi_{2})\nu}
{2\left(\xi_{1}^{3}-\xi_{1}\xi_{2}^{2}-2\xi_{2}^{3}\right)}\right|\leq \frac{1}{32}.\label{3.082}
\end{align}
 From (\ref{3.043}),  we have
\begin{align}
&\hspace{-1.5cm}1-2^{-20}\leq\left|\frac{3\xi_{1}^{2}+2\xi_{1}\xi_{2}+3\xi_{2}^{2}}{\xi_{1}^{2}+\xi_{2}^{2}
+\sqrt{\xi^{2}(\xi^{2}-\xi_{1}\xi_{2})}}\right|-2^{-20}\nonumber\\   &\leq\left|\nu\right|
\leq \left|\frac{3\xi_{1}^{2}+2\xi_{1}\xi_{2}+3\xi_{2}^{2}}{\xi_{1}^{2}+\xi_{2}^{2}
+\sqrt{\xi^{2}(\xi^{2}-\xi_{1}\xi_{2})}}\right|+2^{-20}\leq 3+2^{-20}.\label{3.083}
\end{align}
Combining (\ref{3.082})  with (\ref{3.083}), we have
\begin{align}
&\quad
|\nu\xi_{1}\left(\xi_{1}^{3}-\xi_{1}\xi_{2}^{2}-2\xi_{2}^{3}\right)|\left|\frac{\xi_{1}\xi_{2}(\xi_{1}-3\xi_{2})\nu}
{2\left(\xi_{1}^{3}-\xi_{1}\xi_{2}^{2}-2\xi_{2}^{3}\right)}+1\right|\nonumber\\
&\sim |\xi_{1}\left(\xi_{1}^{3}-\xi_{1}\xi_{2}^{2}-2\xi_{2}^{3}\right)|
\sim 2^{k_{1}+3k_{2}}.\label{3.084}
\end{align}
Combining  (\ref{3.084})  with the fact that $\|G\|_{L^{1}}\approx\|h\|_{L^{2}}^{2}$,
we have
(\ref{3.078}) can be bounded by
\begin{align}
C2^{-\frac{3k_{2}}{2}-k_{1}}\|G\|_{L^{1}}\leq C
2^{-\frac{3k_{2}}{2}-k_{1}}\|h\|_{L^{2}}^{2}.\label{3.085}
\end{align}

This completes the proof of Lemma 3.1.

Inspired by the idea of  Lemma 5.2 of \cite{IKT}, we give the proof of Lemma 3.2.
\begin{Lemma}\label{Lemma3.2}
Assume $\alpha\in \R$ and $k_{1},k_{2},k_{3}\in \mathbf{Z}$, $k_{\rm max}={\rm max}\left\{k_{1},k_{2},k_{3}\right\}\geq20$
 and $k_{\rm min}={\rm min}\left\{k_{1},k_{2},k_{3}\right\}$
 and $j_{1},j_{2},j_{3}\in \mathbf{Z}_{+},j_{\rm max}={\rm max}\left\{j_{1},j_{2},j_{3}\right\}$
and $f_{i}:\R^{3}\rightarrow \R$  are  $L^{2}$  functions supported in $D_{k_{i},\infty,j_{i}}(i=1,2,3)$.
Then, $k_{max}\geq20,$ we have
\begin{align}
&\int_{\SR^{3}}(f_{1}*f_{2})f_{3}d\xi d\mu d\tau
\leq C2^{\frac{j_{1}+j_{2}+j_{3}}{2}}2^{-\frac{j_{\rm max}}{2}
-\frac{k_{\rm max}}{4}-\frac{k_{\rm min}}{4}} \prod\limits_{j=1}^{3}\|f_{j}\|_{L^{2}} .\label{3.086}
\end{align}
\end{Lemma}
\noindent{\bf Proof.}
From (\ref{3.079}), we  assume that $j_{3}={\rm max}\left\{j_{1},j_{2},j_{3}\right\}.$
Then, we have
\begin{align}
\int_{\SR^{3}}(f_{1}*f_{2})f_{3}d\xi d\mu d\tau\leq C\|f_{3}\|_{L^{2}}\|f_{1}*f_{2}\|_{L^{2}}\leq
C\|f_{3}\|_{L^{2}}
\prod_{j=1}^{2}\|\mathscr{F}^{-1}(f_{j})\|_{L^{4}}.\label{3.087}
\end{align}
From Theorem 3.1 of \cite{Hadac2008}, we have
\begin{align}
\left\|\int_{\SR^{2}}|\xi|^{\frac{1}{4}}f_{j}(\xi,\mu)e^{ix\xi+iy\mu}
e^{it\phi(\xi,\mu)}d\xi d\mu\right\|_{L_{xyt}^{4}}\leq C\|f_{j}\|_{L_{\xi\mu}^{2}}(j=1,2).\label{3.088}
\end{align}
From (\ref{3.088}), by using the Cauchy-Schwartz inequality with respect to $\theta$, we have
\begin{align}
&\left\|\int_{\SR^{3}}|\xi_{j}|^{\frac{1}{4}}f_{j}(\xi,\mu,\tau)e^{ix\xi+iy\mu}
e^{it\tau}d\xi d\mu d\tau\right\|_{L_{xyt}^{4}}\nonumber\\
=&\left\|\int_{\SR^{3}}|\xi_{j}|^{\frac{1}{4}}f_{j}^{\#}(\xi,\mu,\theta)e^{ix\xi+iy\mu}
e^{it\theta}e^{it\phi(\xi,\mu)}d\xi d\mu d\theta\right\|_{L_{xyt}^{4}}\nonumber\\
\leq & C\int_{\SR}\|f_{j}^{\#}(\xi,\mu,\theta)\|_{L_{\xi\mu}^{2}}d\theta\leq C
2^{\frac{j_{i}}{2}}\|f_{j}^{\#}(\xi,\mu,\theta)\|_{L_{\xi\mu\theta}^{2}}
\nonumber\\
=&C2^{\frac{j_{i}}{2}}\|f_{j}(\xi,\mu,\tau)\|_{L_{\xi\mu\tau}^{2}}.\label{3.089}
\end{align}
Here $f_{j}^{\#}(j=1,2)$ are defined as in Lemma 3.1.
From (\ref{3.089}),  we have
\begin{align}
\|\mathscr{F}^{-1}(f_{j})\|_{L^{4}}\leq C2^{-\frac{k_{i}}{4}}
2^{\frac{j_{i}}{2}}\|f_{j}(\xi,\mu,\tau)\|_{L_{\xi\mu\tau}^{2}}.\label{3.090}
\end{align}
Inserting (\ref{3.090}) into  (\ref{3.087}) yields
\begin{align}
&\int_{\SR^{3}}(f_{1}*f_{2})f_{3}\leq C\|f_{3}\|_{L^{2}}\|f_{1}*f_{2}\|_{L^{2}}\nonumber\\   &\leq
C2^{\frac{j_{1}+j_{2}}{2}}2^{-\frac{k_{\rm max}+k_{\rm min}}{4}}\|f\|_{L^{2}}
\prod_{j=1}^{2}\|\mathscr{F}^{-1}(f_{j})\|_{L^{4}}.\label{3.091}
\end{align}
Combining the fact with $j_{3}={\rm max}\left\{j_{1},j_{2},j_{3}\right\}$ with (\ref{3.091}),  we have
(\ref{3.086}) is valid.

This completes the proof of Lemma 3.2.

\bigskip
\bigskip

\noindent{\large\bf 4. Bilinear estimates}

\setcounter{equation}{0}

 \setcounter{Theorem}{0}

\setcounter{Lemma}{0}

 \setcounter{section}{4}
 This section is devoted to establishing   Lemmas 4.1-4.3.
 Lemma 4.1 is used to establish  Theorems 1.1. Lemma 4.2 is used to establish  the almost conservation.
Lemma 4.3 is used to establish Lemma 6.1 which is the variant of Theorem 1.1.

 \begin{Lemma}\label{Lemma4.1}
Let $-\frac{9}{8}+16\epsilon\leq s_{1}<0,s_{2}\geq0$ and $u_{j}\in X_{\frac{1}{2}+\epsilon}^{s_{1},s_{2}}(j=1,2)$.
Then, we have
\begin{align}
&\|\partial_{x}(u_{1}u_{2})\|_{X_{-\frac{1}{2}+2\epsilon}^{s_{1},s_{2}}}\leq C
\prod_{j=1}^{2}\|u_{j}\|_{X_{\frac{1}{2}+\epsilon}^{s_{1},s_{2}}}.\label{4.01}
\end{align}
\end{Lemma}
\noindent{\bf Proof.}  To derive (\ref{4.01}),  by duality, it suffices to  show that
\begin{align}
&\left|\int_{\SR^{3}}\bar{u}\partial_{x}(u_{1}u_{2})dxdydt\right|\leq
C\|u\|_{X_{\frac{1}{2}-2\epsilon}^{-s_{1},-s_{2}}}
\prod_{j=1}^{2}
\|u_{j}\|_{X_{\frac{1}{2}+\epsilon}^{s_{1},s_{2}}}
.\label{4.02}
\end{align}
for $u\in X_{\frac{1}{2}-2\epsilon}^{-s_{1},-s_{2}}.$
Let
\begin{align}
&F(\xi,\mu,\tau)=\langle\xi\rangle^{-s_{1}}\langle\mu\rangle^{-s_{2}}
\langle \sigma\rangle^{\frac{1}{2}-2\epsilon}\mathscr{F}u(\xi,\mu,\tau),\nonumber\\   &
F_{j}(\xi_{j},\mu_{j},\tau_{j})=\langle\xi_{j}\rangle^{s_{1}}\langle\mu_{j}\rangle^{s_{2}}
\langle \sigma_{j}\rangle^{\frac{1}{2}+\epsilon}
\mathscr{F}u_{j}(\xi_{j},\mu,\tau_{j})(j=1,2),\label{4.03}
\end{align}
and
\begin{align*}
D:=\left\{(\xi_1,\mu_{1},\tau_1,\xi,\mu,\tau)\in {\rm R^6},
\xi=\sum_{j=1}^{2}\xi_j,\mu=\sum_{j=1}^{2}\mu_{j},\tau=\sum_{j=1}^{2}\tau_j\right\}.
\end{align*}
To derive (\ref{4.02}), from (\ref{4.03}), it suffices to show
\begin{align}
 \int_{D}\frac{|\xi|\langle\xi\rangle^{s_{1}}\langle\mu\rangle^{s_{2}}
F(\xi,\mu,\tau)\prod\limits_{j=1}^{2}F_{j}(\xi_{j},\mu_{j},\tau_{j})}{\langle\sigma_{j}\rangle^{\frac{1}{2}-2\epsilon}
\prod\limits_{j=1}^{2}\langle\xi_{j}\rangle^{s_{1}}\langle\mu_{j}\rangle^{s_{2}}\langle\sigma_{j}\rangle^{\frac{1}{2}+\epsilon}}
dV
\leq C \|F\|_{L_{\tau\xi\mu}^{2}}
\prod_{j=1}^{2}\|F_{j}\|_{L_{\tau\xi\mu}^{2}},\label{4.04}
\end{align}
where $dV:=d\xi_{1}d\mu_{1}d\tau_{1}d\xi d\mu d\tau$.
Without loss of generality, by using the symmetry,  we assume that
$|\xi_{1}|\geq |\xi_{2}|$  and   $F(\xi,\mu,\tau)\geq 0,F_j(\xi_{j},\mu_{j},\tau_{j})\geq 0(j=1,2)$
and
\begin{align*}
D^{*}:=\left\{(\xi_1,\mu_{1},\tau_1,\xi,\mu,\tau)\in D,|\xi_{2}|\geq|\xi_{1}|\right\}.
\end{align*}
Let
\begin{align*}
&\hspace{-0.8cm}\Omega_1=\left\{(\xi_1,\mu_{1},\tau_1,\xi,\mu,\tau)\in D^{*},
 |\xi_2|\leq |\xi_{1}|\leq 80\right\},\\
&\hspace{-0.8cm} \Omega_2=\{ (\xi_1,\mu_{1},\tau_1,\xi,\mu,\tau)\in D^{*},
|\xi_1|\geq 80, |\xi_{1}|\gg|\xi_{2}|,|\xi_{2}|\leq 20\},\\
&\hspace{-0.8cm} \Omega_3=\{ (\xi_1,\mu_{1},\tau_1,\xi,\mu,\tau)\in D^{*},
|\xi_1|\geq 80, |\xi_{1}|\gg|\xi_{2}|,|\xi_{2}|> 20\},\\
&\hspace{-0.8cm}\Omega_4=\{(\xi_1,\mu_{1},\tau_1,\xi,\mu,\tau)\in D^{*},
|\xi_{1}|\geq 80,4|\xi|\leq |\xi_{1}|\sim|\xi_{2}|,|\xi|\leq 20,\xi_{1}\xi_{2}<0\},\\
&\hspace{-0.8cm}\Omega_5=\{(\xi_1,\mu_{1},\tau_1,\xi,\mu,\tau)\in D^{*},
|\xi_{1}|\geq 80,4|\xi|\leq |\xi_{1}|\sim|\xi_{2}|,|\xi|> 20,\xi_{1}\xi_{2}<0\},\\
&\hspace{-0.8cm}\Omega_6=\{(\xi_1,\mu_{1},\tau_1,\xi,\mu,\tau)\in D^{*},
 |\xi_{1}|\geq 80, |\xi_{2}|\sim |\xi_{1}|,\xi_{1}\xi_{2}<0,|\xi|\geq \frac{|\xi_{2}|}{4}\},\\
 &\hspace{-0.8cm}\Omega_7=\left\{(\xi_1,\mu_{1},\tau_1,\xi,\mu,\tau)\in D^{*},
 |\xi_{1}|\geq 80, |\xi_{2}|\sim |\xi_{1}|,\xi_{1}\xi_{2}>0\right\}.
\end{align*}
Obviously, $D^{*}\subset\bigcup\limits_{j=1}^{7}\Omega_{j}.$
We define
\begin{equation}
    K_{1}(\xi_{1},\mu_{1},\tau_{1},\xi,\mu,\tau):=\frac{|\xi|\langle\mu\rangle^{s_{2}}\langle\xi\rangle^{s_{1}}
}{\langle\sigma_{j}\rangle^{\frac{1}{2}-2\epsilon}
\prod\limits_{j=1}^{2}\langle\xi_{j}\rangle^{s_{1}}
\langle\mu_{j}\rangle^{s_{2}}\langle\sigma_{j}\rangle^{\frac{1}{2}+\epsilon}}\label{4.05}
\end{equation}
and
\begin{align*}
{\rm Int}_{j}:=\int_{\Omega_{j}} K_{1}(\xi_{1},\mu_{1},\tau_{1},\xi,\mu,\tau)F(\xi,\mu,\tau)
\prod_{j=1}^{2}F_{j}(\xi_{j},\mu_{j},\tau_{j})
d\xi_{1}d\mu_{1}d\tau_{1}d\xi d\mu d\tau,
\end{align*}
$1\leq j\leq 7, j\in\mathbf{ N}.$
Since $s_{2}\geq0$ and $\mu=\sum\limits_{j=1}^{2}\mu_{j}$, we have
 $\langle \mu\rangle^{s_{2}} \leq \prod\limits_{j=1}^{2}\langle\mu_{j}\rangle^{s_{2}}$,
thus, we have
\begin{align}
K_{1}(\xi_{1},\mu_{1},\tau_{1},\xi,\mu,\tau)\leq\frac{|\xi|\langle\xi\rangle^{s_{1}}}
{\langle\sigma\rangle^{\frac{1}{2}-2\epsilon}
\prod\limits_{j=1}^{2}\langle\xi_{j}\rangle^{s_{1}}\langle\sigma_{j}\rangle^{\frac{1}{2}+\epsilon}}.\label{4.06}
\end{align}

Now we estimate the integrals Int$_{j}$ over the above seven regions one by one.

(I)  Region $\Omega_{1}.$ In this region $|\xi|\leq |\xi_{1}|+|\xi_{2}|\leq 160,$ this case can be proved similarly to case
$low+low\longrightarrow low$ of Pages 344--345 of Theorem 3.1 in \cite{LX}.

(II) Region $\Omega_{2}.$ In this region, we have  $|\xi|\sim |\xi_{1}|$.

  By using the Cauchy-Schwartz inequality with respect to $\xi_{1},\mu_{1},\tau_{1}$,
 from (\ref{4.06}),  we have
\begin{align}
&{\rm Int_{2}}\leq C\int_{\SR^{3}}\frac{|\xi|}{\langle \sigma \rangle^{\frac{1}{2}-2\epsilon}}
\biggl(\int_{\SR^{3}}
   \langle\sigma_{1}\rangle^{-(1+2\epsilon)}
    \langle\sigma_{2}\rangle^{-(1+2\epsilon)}
     d\xi_{1}d\mu_{1}d\tau_{1}
      \biggr)^{\frac{1}{2}}\nonumber\\   &\qquad\qquad\qquad \times\left(\int_{\SR^{3}}
\prod\limits_{j=1}^{2}\left|F_{j}(\xi_{j},\mu_{j},\tau_{j})\right|^{2}
d\xi_{1}d\mu_{1}d\tau_{1}\right)^{\frac{1}{2}}F(\xi,\mu,\tau)d\xi d\mu d\tau.\label{4.07}
\end{align}
By using (\ref{2.02}), we have
\begin{align}\nonumber
& \frac{|\xi|}
       {\langle \sigma \rangle^{\frac{1}{2}-2\epsilon}}
 \biggl(\int_{\SR^{3}}
   \langle\sigma_{1}\rangle^{-(1+2\epsilon)}
    \langle\sigma_{2}\rangle^{-(1+2\epsilon)}
     d\xi_{1}d\mu_{1}d\tau_{1}
      \biggr)^{\frac{1}{2}}\\
\leq& C\frac{|\xi|}{\langle \sigma \rangle^{\frac{1}{2}-2\epsilon}}
\left(\int_{\SR^{2}}\frac{d\xi_{1}d\mu_{1}}{
\langle\tau+\phi(\xi_{1},\mu_{1})
+\phi(\xi_{2},\mu_{2})\rangle^{1+2\epsilon}}
\right)^{\frac{1}{2}}\label{4.08}.
\end{align}
Let $\nu=\tau+\phi(\xi_{1},\mu_{1})+\phi(\xi_{2},\mu_{2})$ and
 $\Delta=\xi\xi_{1}\xi_{2}(5\xi^{2}-5\xi\xi_{1}+5\xi_{1}^{2}),$ since $|\xi_{1}|\gg|\xi_{2}|$,
then we have  the absolute value of Jacobian determinant equals
\begin{align}
   \left|\frac{\partial(\Delta,\nu)}{\partial(\xi_{1},\mu_{1})}\right|
=&  2\left|\frac{\mu_{1}}{\xi_{1}}-\frac{\mu_{2}}{\xi_{2}}\right|
   \left|5(\xi_{1}^{4}-\xi_{2}^{4})-3\alpha(\xi_{1}^{2}-\xi_{2}^{2})\right|\nonumber\\
=&  2\left|\sigma-\nu-\Delta\right|^{\frac{1}{2}}\left|\frac{\xi}{\xi_{1}\xi_{2}}\right|^{\frac{1}{2}}
   \left|5(\xi_{1}^{4}-\xi_{2}^{4})-3\alpha(\xi_{1}^{2}-\xi_{2}^{2})\right|
    \nonumber\\
\sim&
    \left|\sigma-\nu+\delta\right|^{\frac{1}{2}}
    \left|\frac{\xi}{\xi_{1}\xi_{2}}\right|^{\frac{1}{2}}|\xi_{1}|^{4}
.\label{4.09}
\end{align}
Inserting (\ref{4.09}) into (\ref{4.08}), by using (\ref{2.03}), we have
\begin{align}\nonumber
&\frac{|\xi|}{\langle \sigma \rangle^{\frac{1}{2}-2\epsilon}}
\biggl(\int_{\SR^{3}}
   \langle\sigma_{1}\rangle^{-(1+2\epsilon)}
    \langle\sigma_{2}\rangle^{-(1+2\epsilon)}
     d\xi_{1}d\mu_{1}d\tau_{1}
      \biggr)^{\frac{1}{2}}\\
\leq & C\frac{|\xi|}
   {\langle \sigma \rangle^{\frac{1}{2}-2\epsilon}}
    \left(\int_{\SR^{2}}\frac{d\xi_{1}d\mu_{1}}{
\langle\tau+\phi(\xi_{1},\mu_{1})+\phi(\xi_{2},\mu_{2})\rangle^{1+2\epsilon}}\right)^{\frac{1}{2}}\nonumber\\
\leq &\frac{C}{|\xi|\langle \sigma \rangle^{\frac{1}{2}-2\epsilon}}
\left(\int_{\SR^{2}}\frac{d\nu d\Delta}{\left|\sigma-\nu-\Delta\right|^{\frac{1}{2}}
\langle\nu\rangle^{1+2\epsilon}}\right)^{\frac{1}{2}}\nonumber\\
\leq & \frac{C}{|\xi|\langle \sigma \rangle^{\frac{1}{2}-2\epsilon}}
\left(\int_{|\Delta|<20|\xi|^{4}}\frac{d\Delta}{
\langle\Delta-\sigma\rangle^{\frac{1}{2}}}\right)^{\frac{1}{2}}
\label{4.010}.
\end{align}
When $|\sigma|<20|\xi|^{4},$ combining (\ref{4.010}) with (\ref{2.01}), we have
\begin{align}
 \frac{C}{|\xi|\langle \sigma \rangle^{\frac{1}{2}-2\epsilon}}\left(\int_{|\Delta|<20|\xi|^{4}}\frac{d\Delta}{
\langle\Delta-\sigma\rangle^{\frac{1}{2}}}\right)^{\frac{1}{2}}
\leq C\frac{|\xi|}{|\xi|\langle \sigma \rangle^{\frac{1}{2}-2\epsilon}}\leq C\label{4.011}.
\end{align}
When $|\sigma|\geq20|\xi|^{4},$ from   (\ref{4.010}), we have
\begin{align}
 \frac{C}{|\xi|\langle \sigma \rangle^{\frac{1}{2}-2\epsilon}}\left(\int_{|\Delta|<20|\xi|^{4}}\frac{d\Delta}{
\langle\Delta-\sigma\rangle^{\frac{1}{2}}}\right)^{\frac{1}{2}}\leq C\frac{|\xi|^{2}}
{|\xi|\langle \sigma \rangle^{\frac{1}{2}-2\epsilon}}\leq C\frac{|\xi|}
{\langle \sigma \rangle^{\frac{1}{2}-2\epsilon}}\leq C\label{4.012}.
\end{align}
Combining (\ref{4.08}) with (\ref{4.09})-(\ref{4.012}), we have
\begin{align}
&\frac{|\xi|}{\langle \sigma \rangle^{\frac{1}{2}-2\epsilon}}\biggl(\int_{\SR^{3}}
   \langle\sigma_{1}\rangle^{-(1+2\epsilon)}
    \langle\sigma_{2}\rangle^{-(1+2\epsilon)}
     d\xi_{1}d\mu_{1}d\tau_{1}
      \biggr)^{\frac{1}{2}}\leq C\label{4.013}.
\end{align}
Inserting (\ref{4.013}) into (\ref{4.07}), by using the Cauchy-Schwartz
 inequality with respect to $\xi,\mu,\tau$, we have
\begin{align}
{\rm Int_{2}}
&\leq C\int_{\SR^{3}}\frac{|\xi|}{\langle \sigma \rangle^{\frac{1}{2}-2\epsilon}}
\biggl(\int_{\SR^{3}}
   \langle\sigma_{1}\rangle^{-(1+2\epsilon)}
    \langle\sigma_{2}\rangle^{-(1+2\epsilon)}
     d\xi_{1}d\mu_{1}d\tau_{1}
      \biggr)^{\frac{1}{2}}\nonumber\\
      &\qquad\qquad  \times\left(\int_{\SR^{3}}
\prod\limits_{j=1}^{2}\left|F_{j}(\xi_{j},\mu_{j},\tau_{j})\right|^{2}
d\xi_{1}d\mu_{1}d\tau_{1}\right)^{\frac{1}{2}}F(\xi,\mu,\tau)d\xi d\mu d\tau\nonumber\\
&\leq C\int_{\SR^{3}}\left(\int_{\SR^{3}}
\prod\limits_{j=1}^{2}\left|F_{j}(\xi_{j},\mu_{j},\tau_{j})\right|^{2}
d\xi_{1}d\mu_{1}d\tau_{1}\right)^{\frac{1}{2}}F(\xi,\mu,\tau)d\xi d\mu d\tau\nonumber\\   &
\leq C\|F\|_{L_{\tau\xi\mu}^{2}}\
prod_{j=1}^{2}\|F_{j}\|_{L_{\tau\xi\mu}^{2}}.\label{4.014}
\end{align}

 (III)  Region $\Omega_{3}.$ In this region, we have  $|\xi|\sim |\xi_{1}|$.
In this region, we consider
\begin{align}
|\sigma-\sigma_{1}-\sigma_{2}|
&=\left|\xi\xi_{1}\xi_{2}(5\xi^{2}-5\xi\xi_{1}+5\xi_{1}^{2})-\frac{\xi_{1}\xi_{2}}{\xi}
\left|\frac{\mu_{1}}{\xi_{1}}-\frac{\mu_{2}}{\xi_{2}}\right|^{2}\right|\nonumber\\   &\geq
\frac{\left|\xi\xi_{1}\xi_{2}(5\xi^{2}-5\xi\xi_{1}+5\xi_{1}^{2})\right|}{2^{70}}\label{4.015}
\end{align}
and
\begin{align}
|\sigma-\sigma_{1}-\sigma_{2}|
&=\left|\xi\xi_{1}\xi_{2}(5\xi^{2}-5\xi\xi_{1}+5\xi_{1}^{2})-\frac{\xi_{1}\xi_{2}}{\xi}
\left|\frac{\mu_{1}}{\xi_{1}}-\frac{\mu_{2}}{\xi_{2}}\right|^{2}\right|\nonumber\\   &
< \frac{\left|\xi\xi_{1}\xi_{2}(5\xi^{2}-5\xi\xi_{1}+5\xi_{1}^{2})\right|}{2^{70}}\label{4.016},
\end{align}
respectively.

When (\ref{4.015}) is valid, we have  one of the following three cases must occur:
\begin{align}
&|\sigma|:={\rm max}\left\{|\sigma|,|\sigma_{1}|,|\sigma_{2}|\right\}
\geq C\left|\xi\xi_{1}\xi_{2}(5\xi^{2}-5\xi\xi_{1}+5\xi_{1}^{2})\right|,\label{4.017}\\
&|\sigma_{1}|:={\rm max}\left\{|\sigma|,|\sigma_{1}|,|\sigma_{2}|\right\}
\geq C\left|\xi\xi_{1}\xi_{2}(5\xi^{2}-5\xi\xi_{1}+5\xi_{1}^{2})\right|,\label{4.018}\\
&|\sigma_{2}|:={\rm max}\left\{|\sigma|,|\sigma_{1}|,|\sigma_{2}|\right\}
\geq C\left|\xi\xi_{1}\xi_{2}(5\xi^{2}-5\xi\xi_{1}+5\xi_{1}^{2})\right|.\label{4.019}
\end{align}

When (\ref{4.017}) is valid, since  $-\frac{9}{8}+16\epsilon\leq s_{1}<0$, we have
\begin{align}
K_{1}(\xi_{1},\mu_{1},\tau_{1},\xi,\mu,\tau)
&\leq\frac{|\xi|\langle\xi\rangle^{s_{1}}}
{\langle\sigma\rangle^{\frac{1}{2}-2\epsilon}
\prod\limits_{j=1}^{2}\langle\xi_{j}\rangle^{s_{1}}\langle\sigma_{j}\rangle^{\frac{1}{2}+\epsilon}}
\leq C\frac{|\xi|^{-1+8\epsilon}|\xi_{2}|^{-s_{1}-\frac{1}{2}+2\epsilon}}
{\prod\limits_{j=1}^{2}\langle\sigma_{j}\rangle^{\frac{1}{2}+\epsilon}}
\nonumber\\
&\leq C\frac{|\xi|^{-1+8\epsilon}|\xi_{2}|^{\frac{5}{8}-14\epsilon}}
{\prod\limits_{j=1}^{2}\langle\sigma_{j}\rangle^{\frac{1}{2}+\epsilon}}
\leq C\frac{|\xi_{1}||\xi_{2}|^{-\frac{1}{2}}}
{\prod\limits_{j=1}^{2}\langle\sigma_{j}\rangle^{\frac{1}{2}+\epsilon}}
.\label{4.020}
\end{align}
Thus, combining (\ref{2.017}) with (\ref{4.020}),  we have
\begin{align*}
|{\rm Int_{3}}|\leq C\|F\|_{L_{\tau\xi\mu}^{2}}
\prod\limits_{j=1}^{2}\|F_{j}\|_{L_{\tau\xi\mu}^{2}}.
\end{align*}

When (\ref{4.018}) is valid, since  $-\frac{9}{8}+16\epsilon\leq s_{1}<0$
 and
$$
\langle \sigma\rangle^{-\frac{1}{2}+2\epsilon}
 \langle \sigma_{1}\rangle^{-\frac{1}{2}-\epsilon}\leq \langle \sigma\rangle^{-\frac{1}{2}-\epsilon}\langle
  \sigma_{1}\rangle^{-\frac{1}{2}+2\epsilon},
$$
we have
\begin{align}
K_{1}(\xi_{1},\mu_{1},\tau_{1},\xi,\mu,\tau)
&\leq\frac{|\xi|\langle\xi\rangle^{s_{1}}}
{\langle\sigma\rangle^{\frac{1}{2}-2\epsilon}
\prod\limits_{j=1}^{2}\langle\xi_{j}\rangle^{s_{1}}\langle\sigma_{j}\rangle^{\frac{1}{2}+\epsilon}}\nonumber\\
&\leq C\frac{|\xi|^{-1+8\epsilon}|\xi_{2}|^{-s_{1}-\frac{1}{2}+2\epsilon}}
{\langle\sigma\rangle^{\frac{1}{2}+\epsilon}\langle\sigma_{2}\rangle^{\frac{1}{2}+\epsilon}}
\leq C\frac{|\xi|^{-\frac{1}{2}}|\xi_{2}|}
{\langle\sigma\rangle^{\frac{1}{2}+\epsilon}\langle\sigma_{2}\rangle^{\frac{1}{2}+\epsilon}}
.\label{4.021}
\end{align}
Thus, combining (\ref{2.019}) with (\ref{4.021}),  we have
\begin{align*}
|{\rm Int_{3}}|\leq C\|F\|_{L_{\tau\xi\mu}^{2}} \prod\limits_{j=1}^{2}\|F_{j}\|_{L_{\tau\xi\mu}^{2}}.
\end{align*}

When (\ref{4.019}) is valid, since  $-\frac{9}{8}+16\epsilon\leq s_{1}<0$
 and
$$\langle \sigma\rangle^{-\frac{1}{2}+2\epsilon}\langle \sigma_{2}
 \rangle^{-\frac{1}{2}-\epsilon}\leq \langle \sigma\rangle^{-\frac{1}{2}-\epsilon}
 \langle \sigma_{2}\rangle^{-\frac{1}{2}+2\epsilon},
$$
we have
\begin{align}
K_{1}(\xi_{1},\mu_{1},\tau_{1},\xi,\mu,\tau)
&\leq\frac{|\xi|\langle\xi\rangle^{s_{1}}}
{\langle\sigma\rangle^{\frac{1}{2}-2\epsilon}
\prod\limits_{j=1}^{2}\langle\xi_{j}\rangle^{s_{1}}\langle\sigma_{j}\rangle^{\frac{1}{2}+\epsilon}}
\leq C\frac{|\xi|^{-1+8\epsilon}|\xi_{1}|^{-s_{1}-\frac{1}{2}+2\epsilon}}
{\langle\sigma_{1}\rangle^{\frac{1}{2}+\epsilon}\langle\sigma\rangle^{\frac{1}{2}+\epsilon}}\nonumber\\
&\leq C\frac{|\xi_{1}|^{\frac{1}{4}-14\epsilon}|\xi|^{-1+8\epsilon}}
{\langle\sigma_{1}\rangle^{\frac{1}{2}+\epsilon}\langle\sigma\rangle^{\frac{1}{2}+\epsilon}}
\leq C\frac{|\xi_{1}|^{-\frac{1}{2}}|\xi|}
{\langle\sigma_{1}\rangle^{\frac{1}{2}+\epsilon}\langle\sigma\rangle^{\frac{1}{2}+\epsilon}}
.\label{4.022}
\end{align}
Thus, combining (\ref{2.020}) with (\ref{4.022}),  we have
\begin{align*}
|{\rm Int_{3}}|\leq C\|F\|_{L_{\tau\xi\mu}^{2}} \prod\limits_{j=1}^{2}\|F_{j}\|_{L_{\tau\xi\mu}^{2}} .
\end{align*}

When (\ref{4.016}) is valid,
we consider the following two cases
respectively.
\begin{align}
&{\rm max}\left\{|\sigma|,|\sigma_{1}|,|\sigma_{2}|\right\}\geq
\frac{\left|\xi\xi_{1}\xi_{2}(5\xi^{2}-5\xi\xi_{1}+5\xi_{1}^{2})\right|}{2^{80}},\label{4.023}\\
&{\rm max}\left\{|\sigma|,|\sigma_{1}|,|\sigma_{2}|\right\}<
\frac{\left|\xi\xi_{1}\xi_{2}(5\xi^{2}-5\xi\xi_{1}+5\xi_{1}^{2})\right|}{2^{80}}.\label{4.024}
\end{align}

 We dyadically decompose the spectra as
\begin{align*}
\langle\sigma\rangle\sim 2^{j},\quad
\langle\sigma_{1}\rangle\sim 2^{j_{1}},\quad
\langle\sigma_{2}\rangle\sim 2^{j_{2}},\quad
|\xi|\sim 2^{k},\quad
|\xi_{1}|\sim 2^{k_{1}},\quad
|\xi_{2}|\sim 2^{k_{2}}.
\end{align*}
We define
\begin{align*}
&f_{k_{m},j_m}:
=\left|\eta_{k_{m}}(\xi_{m})\eta_{j_{m}}(\sigma_{m})
   F_{l}(\xi_{m},\mu_{m},\tau_{m})\right|,\quad(m=1,2),\\
&f_{k,j}:=\left|\eta_{k}(\xi)\eta_{j}(\sigma)F(\xi,\mu,\tau)\right|.
\end{align*}
When (\ref{4.023}) is valid,  by  using Lemma 3.2, since $-\frac{9}{8}+16\epsilon\leq s_{1}<0,$ we have
\begin{align}
 {\rm Int_{3}}
&\leq C\sum\limits_{k,k_{1},k_{2}>0}
        \sum\limits_{j_{1},j_{2},j\geq0}
        2^{-j(\frac{1}{2}-2\epsilon)-(j_{1}
            +j_{2})(\frac{1}{2}+\epsilon)
            -k_{2}s_{1}+k}
         \int_{\SR^{6}}f_{k,j}
        \prod_{m=1}^{2}f_{k_{m},j_{m}}dV \nonumber\\
&\leq C \sum\limits_{k,k_{1},k_{2}>0}
         \sum\limits_{j_{1},j_{2},j\geq0}
         2^{2j\epsilon-(j_{1}+j_{2})\epsilon
          -\frac{j_{\rm max}}{2}-k_{2}(s_{1}+\frac{1}{4})
          +\frac{3k}{4}}\|f_{k,j}\|_{L^{2}}
  \prod_{m=1}^{2}
    \|f_{k_{m},j_{m}}\|_{L^{2}} ;\label{4.025}
\end{align}
when $j=j_{\rm max},$ from (\ref{4.025}), since $-\frac{9}{8}+16\epsilon\leq s_{1}<0,$ we have
\begin{align}
 {\rm Int_{3}}
 &\leq  C\sum\limits_{k,k_{1},k_{2}>0}
 \sum\limits_{j_{1},j_{2},j\geq0}
  2^{-j(\frac{1}{2}-2\epsilon)-(j_{1}+j_{2})\epsilon
  -k_{2}(s_{1}+\frac{1}{4})+\frac{3k}{4}}
 \int_{\SR^{6}}f_{k,j}\prod_{m=1}^{2}f_{k_{m},j_{m}}dV \nonumber\\
&\leq C
  \sum\limits_{k, k_{1},k_{2}>0}
  2^{-k_{2}(s_{1}+\frac{3}{4}-2\epsilon)
  -k(\frac{5}{4}-8\epsilon)}\|f\|_{L^{2}}
  \prod_{m=1}^{2}\|f_{m}\|_{L^{2}} \nonumber\\
&\leq C\sum\limits_{k, k_{1},k_{2}>0}
  2^{-k_{2}(s_{1}+\frac{9}{8})
  -k(\frac{7}{8}-10\epsilon)}
  \|f\|_{L^{2}} \prod_{m=1}^{2}\|f_{m}\|_{L^{2}} \nonumber\\
& \leq C \|f\|_{L^{2}}
 \prod_{m=1}^{2}\|f_{m}\|_{L^{2}} ;\label{4.026}
\end{align}
when $j_{1}=j_{\rm max},$ from (\ref{4.025}), since $-\frac{9}{8}+16\epsilon\leq s_{1}<0,$ we have
\begin{align}
 {\rm Int_{3}}
 &\leq   C\sum\limits_{k,k_{1},k_{2}>0}
          \sum\limits_{j_{1},j_{2},j\geq0}
           2^{2j\epsilon-(j_{1}+j_{2})\epsilon
           -\frac{j_{\rm max}}{2}-k_{2}(s_{1}+\frac{1}{4})
           +\frac{3k}{4}}\|f_{k,j}\|_{L^{2}}
           \prod_{m=1}^{2}\|f_{k_{m},j_{m}}\|_{L^{2}}
           \nonumber\\
& \leq C \sum\limits_{k,k_{1},k_{2}>0}
  2^{-k_{2}(s_{1}+\frac{3}{4}-\epsilon)
  -k(\frac{5}{4}-\epsilon)}\|f\|_{L^{2}}
  \prod_{m=1}^{2}\|f_{m}\|_{L^{2}} \nonumber\\
&\leq C\sum\limits_{k, k_{1},k_{2}>0}
    2^{-k_{2}(s_{1}+\frac{9}{8})
    -k(\frac{7}{8}-5\epsilon)}
    \|f\|_{L^{2}} \prod_{m=1}^{2}\|f_{m}\|_{L^{2}}\nonumber\\
&\leq C \|f\|_{L^{2}}
 \prod_{m=1}^{2}\|f_{m}\|_{L^{2}} ;\label{4.027}
\end{align}
when $j_{2}=j_{\rm max},$ this case can be proved similarly to case $j_{1}=j_{\rm max}.$

 When (\ref{4.024}) is valid, by using (2) of  Lemma 3.1, since $-\frac{9}{8}+16\epsilon\leq s_{1}<0,$ we have
\begin{align}
 {\rm Int_{3}}
&\leq  C\sum\limits_{k,k_{1},k_{2}>0}
   \sum\limits_{j_{1},j_{2},j\geq0}
   2^{-j(\frac{1}{2}-2\epsilon)
   -(j_{1}+j_{2})(\frac{1}{2}+\epsilon)
   -k_{2}s_{1}+k}
    \int_{\SR^{6}}f_{k,j}\prod_{m=1}^{2}f_{k_{m},j_{m}}dV\nonumber\\
&\leq C \sum\limits_{k,k_{1},k_{2}>0}
        \sum\limits_{j_{1},j_{2},j\geq0}
        2^{2j\epsilon-(j_{1}+j_{2})\epsilon
        -k_{2}(s_{1}+\frac{1}{2})
        -\frac{3}{2}k}\|f_{k,j}\|_{L^{2}}
        \prod_{m=1}^{2}\|f_{k_{m},j_{m}}\|_{L^{2}}\nonumber\\
&\leq C \sum\limits_{k,k_{1},k_{2}>0}
    2^{-k_{2}(s_{1}+\frac{9}{8})
    -\frac{3k}{8}}\|f\|_{L^{2}} \prod_{m=1}^{2}\|f_{m}\|_{L^{2}} \nonumber\\
&\leq C \|f\|_{L^{2}}
      \prod_{m=1}^{2}\|f_{m}\|_{L^{2}} .\label{4.028}
\end{align}

 (IV)  Region $\Omega_{4}.$
In this case, we consider (\ref{4.015}), (\ref{4.016}), respectively.

  When (\ref{4.015}) is valid,  one of (\ref{4.017})-(\ref{4.019}) must occur.

When (\ref{4.017}) is valid, since  $-\frac{9}{8}+16\epsilon\leq s_{1}<0$, we have
\begin{align}
K_{1}(\xi_{1},\mu_{1},\tau_{1},\xi,\mu,\tau)
&\leq\frac{|\xi|\langle\xi\rangle^{s_{1}}}
{\langle\sigma\rangle^{\frac{1}{2}-2\epsilon}
\prod\limits_{j=1}^{2}\langle\xi_{j}\rangle^{s_{1}}\langle\sigma_{j}\rangle^{\frac{1}{2}+\epsilon}}
 \leq C\frac{|\xi|^{\frac{1}{2}+2\epsilon}|\xi_{2}|^{-2s_{1}-2+8\epsilon}}
{\prod\limits_{j=1}^{2}\langle\sigma_{j}\rangle^{\frac{1}{2}+\epsilon}}
   \nonumber\\
&\leq C\frac{|\xi|^{\frac{1}{2}+2\epsilon}
  |\xi_{2}|^{-\frac{1}{2}-24\epsilon}}
{\prod\limits_{j=1}^{2}\langle\sigma_{j}\rangle^{\frac{1}{2}+\epsilon}}
 \leq C\frac{|\xi_{1}|^{-\frac{1}{2}}|\xi_{2}|}
{\prod\limits_{j=1}^{2}\langle\sigma_{j}\rangle^{\frac{1}{2}+\epsilon}}
.\label{4.029}
\end{align}
Thus, combining (\ref{2.017}) with (\ref{4.029}),  we have
\begin{align*}
|{\rm Int_{4}}|
\leq C\|F\|_{L_{\tau\xi\mu}^{2}}
\prod\limits_{j=1}^{2}\|F_{j}\|_{L_{\tau\xi\mu}^{2}} .
\end{align*}

When (\ref{4.018}) is valid, since  $-\frac{9}{8}+16\epsilon\leq s_{1}<0$ and
 $$\langle \sigma\rangle^{-\frac{1}{2}+2\epsilon}\langle \sigma_{1}\rangle^{-\frac{1}{2}-\epsilon}
 \leq \langle \sigma\rangle^{-\frac{1}{2}-\epsilon}\langle \sigma_{1}\rangle^{-\frac{1}{2}+2\epsilon},
 $$
we have
\begin{align}
K_{1}(\xi_{1},\mu_{1},\tau_{1},\xi,\mu,\tau)
&\leq\frac{|\xi|\langle\xi\rangle^{s_{1}}}
{\langle\sigma\rangle^{\frac{1}{2}-2\epsilon}
\prod\limits_{j=1}^{2}\langle\xi_{j}\rangle^{s_{1}}\langle\sigma_{j}\rangle^{\frac{1}{2}+\epsilon}}
\leq C\frac{|\xi|^{\frac{1}{2}+2\epsilon}|\xi_{2}|^{-2s_{1}-2+8\epsilon}}
{\langle\sigma\rangle^{\frac{1}{2}+\epsilon}\langle\sigma_{2}\rangle^{\frac{1}{2}+\epsilon}}\nonumber\\
&\leq C\frac{|\xi|^{\frac{1}{2}+2\epsilon}|\xi_{2}|^{-\frac{1}{2}-24\epsilon}}
{\langle\sigma\rangle^{\frac{1}{2}+\epsilon}\langle\sigma_{2}\rangle^{\frac{1}{2}+\epsilon}}
\leq C\frac{|\xi|^{-\frac{1}{2}}|\xi_{2}|}
{\langle\sigma\rangle^{\frac{1}{2}+\epsilon}\langle\sigma_{2}\rangle^{\frac{1}{2}+\epsilon}}
.\label{4.030}
\end{align}
Thus, combining (\ref{2.019}) with (\ref{4.030}),  we have
\begin{align*}
|{\rm Int_{4}}|\leq C\|F\|_{L_{\tau\xi\mu}^{2}}
\prod\limits_{j=1}^{2}\|F_{j}\|_{L_{\tau\xi\mu}^{2}}.
\end{align*}

When (\ref{4.019}) is valid, this case can be proved similarly to (\ref{4.018}).

When (\ref{4.016}) is valid, we consider (\ref{4.023})  and (\ref{4.024}), respectively.

We dyadically decompose the spectra as
\begin{align*}
\langle\sigma\rangle\sim 2^{j},\quad
\langle\sigma_{1}\rangle\sim 2^{j_{1}},\quad
\langle\sigma_{2}\rangle\sim 2^{j_{2}},\quad
|\xi|\sim 2^{k},\quad
|\xi_{1}|\sim 2^{k_{1}},\quad
|\xi_{2}|\sim 2^{k_{2}}.
\end{align*}
We define
\begin{align*}
&f_{k_{m},j_m}:
=\left|\eta_{k_{m}}(\xi_{m})\eta_{j_{m}}(\sigma_{m})
F_{l}(\xi_{m},\mu_{m},\tau_{m})\right|,\qquad
(m=1,2),\\
&f_{k,j}:=\left|\eta_{k}(\xi)\eta_{j}(\sigma)F(\xi,\mu,\tau)\right|,\quad
dV=d\xi_{1}d\mu_{1}d\tau_{1}d\xi d\mu d\tau.
\end{align*}
When (\ref{4.023}) is valid,  by   Lemma 3.2, since $-\frac{9}{8}+16\epsilon\leq s_{1}<0,$ we have
\begin{align}
 {\rm Int_{4}}
 &\leq  C\sum\limits_{k_{1},k_{2}>0,\>k}
 \sum\limits_{j_{1},j_{2},j\geq0}
2^{-j(\frac{1}{2}-2\epsilon)-(j_{1}+j_{2})
(\frac{1}{2}+\epsilon)-2k_{2}s_{1}+k}
 \int_{\SR^{6}}f_{k,j}\prod_{m=1}^{2}
 f_{k_{m},j_{m}}dV\nonumber\\
  &\leq C
\sum\limits_{k_{1},k_{2}>0,\>k}
\sum\limits_{j_{1},j_{2},j\geq0}
2^{2j\epsilon-(j_{1}+j_{2})\epsilon-\frac{j_{\rm max}}{2}-k_{2}(2s_{1}+\frac{1}{4})
+\frac{3k}{4}}\|f_{k,j}\|_{L^{2}}
\prod_{m=1}^{2}\|f_{k_{m},j_{m}}\|_{L^{2}};\label{4.031}
\end{align}
when $j=j_{\rm max},$ from (\ref{4.031}), since $\frac{9}{8}+16\epsilon\leq s_{1}<0,$ we have
\begin{align}
{\rm Int_{4}}
&\leq  C\sum\limits_{k_{1},k_{2}>0,\>k}
   \sum\limits_{j_{1},j_{2},j\geq0}
2^{-j(\frac{1}{2}-2\epsilon)
-(j_{1}+j_{2})(\frac{1}{2}+\epsilon)
 -2k_{2}s_{1}+k}
 \int_{\SR^{6}}f_{k,j}
 \prod_{m=1}^{2}f_{k_{m},j_{m}}dV \nonumber\\   &\leq C
\sum\limits_{k_{1},k_{2}>0,\>k}
2^{-k_{2}(2s_{1}+\frac{9}{4}-8\epsilon)
+k(\frac{1}{4}+2\epsilon)}\|f\|_{L^{2}}
\prod_{m=1}^{2}\|f_{m}\|_{L^{2}}\nonumber\\
&\leq C
\|f\|_{L^{2}} \prod_{m=1}^{2}\|f_{m}\|_{L^{2}} ;\label{4.032}
\end{align}
when $j_{1}=j_{\rm max},$ from (\ref{4.031}), since $-\frac{9}{8}+16\epsilon\leq s_{1}<0,$ we have
\begin{align}
 {\rm Int_{4}}
 &\leq   C\sum\limits_{k_{1},k_{2}>0,\>k}
 \sum\limits_{j_{1},j_{2},j\geq0}
2^{2j\epsilon-(j_{1}+j_{2})\epsilon
-\frac{j_{\rm max}}{2}
-k_{2}(2s_{1}+\frac{1}{4})
+\frac{3k}{4}}\|f_{k,j}\|_{L^{2}}
\prod_{m=1}^{2}\|f_{k_{m},j_{m}}\|_{L^{2}}\nonumber\\
&\leq C
\sum\limits_{k_{1},k_{2}>0,\>k}
2^{-k_{2}(2s_{1}+\frac{9}{4}-4\epsilon)
+k(\frac{1}{4}+\epsilon)}\|f\|_{L^{2}} \prod_{m=1}^{2}\|f_{m}\|_{L^{2}}\nonumber\\
&\leq C
\|f\|_{L^{2}}\prod_{m=1}^{2}\|f_{m}\|_{L^{2}};\label{4.033}
\end{align}
when $j_{2}=j_{\rm max},$ this case can be proved similarly to case $j_{1}=j_{\rm max}.$

  When (\ref{4.024}) is valid,  by using (1) of  Lemma 3.1, since $-\frac{9}{8}+16\epsilon\leq s_{1}<0,$ we have
\begin{align}
{\rm Int_{4}}
&\leq  C\sum\limits_{k_{1},k_{2}>0,\>k}\sum\limits_{j_{1},j_{2},j\geq0}
2^{-j(\frac{1}{2}-2\epsilon)
-(j_{1}+j_{2})(\frac{1}{2}+\epsilon)-2k_{2}s_{1}+k}
\int_{\SR^{6}}f_{k,j}\prod_{m=1}^{2}f_{k_{m},j_{m}}dV\nonumber\\   &\leq C
\sum\limits_{k_{1},k_{2}>0,\>k}\sum\limits_{j_{1},j_{2},j\geq0}
2^{2j\epsilon-(j_{1}+j_{2})\epsilon
-k_{2}(2s_{1}+\frac{7}{2})+\frac{3k}{2}}\|f_{k,j}\|_{L^{2}}
\prod_{m=1}^{2}\|f_{k_{m},j_{m}}\|_{L^{2}}
\nonumber\\
&\leq C
\sum\limits_{k_{1},k_{2}>0,\>k}
2^{-k_{2}(2s_{1}+\frac{7}{2}-8\epsilon)+k(\frac{3}{2}+\epsilon)}\|f\|_{L^{2}} \prod_{m=1}^{2}\|f_{m}\|_{L^{2}} \nonumber\\   &\leq C
\|f\|_{L^{2}} \prod_{m=1}^{2}\|f_{m}\|_{L^{2}} .\label{4.034}
\end{align}

(V)  Region $\Omega_{5}.$
In this region, we consider (\ref{4.015}), (\ref{4.016}), respectively.

  When (\ref{4.015}) is valid,  one of (\ref{4.017})-(\ref{4.019}) must occur.

 When (\ref{4.017}) is valid, since  $-\frac{9}{8}+16\epsilon\leq s_{1}<0$, we have
\begin{align}
&K_{1}(\xi_{1},\mu_{1},\tau_{1},\xi,\mu,\tau)\leq\frac{|\xi|\langle\xi\rangle^{s_{1}}}
{\langle\sigma\rangle^{\frac{1}{2}-2\epsilon}
\prod\limits_{j=1}^{2}\langle\xi_{j}\rangle^{s_{1}}\langle\sigma_{j}\rangle^{\frac{1}{2}+\epsilon}}
\nonumber\\
&\leq C\frac{|\xi|^{s_{1}+\frac{1}{2}-2\epsilon}|\xi_{2}|^{-2s_{1}-2+8\epsilon}}
{\prod\limits_{j=1}^{2}\langle\sigma_{j}\rangle^{\frac{1}{2}+\epsilon}}
\leq C\frac{|\xi|^{-\frac{5}{8}+14\epsilon}|\xi_{2}|^{\frac{1}{4}-24\epsilon}}
{\prod\limits_{j=1}^{2}\langle\sigma_{j}\rangle^{\frac{1}{2}+\epsilon}}
\leq C\frac{|\xi_{1}|^{-\frac{1}{2}}|\xi_{2}|}
{\prod\limits_{j=1}^{2}\langle\sigma_{j}\rangle^{\frac{1}{2}+\epsilon}}
.\label{4.035}
\end{align}
Thus, combining (\ref{2.017}) with (\ref{4.035}),  we have
\begin{align*}
|{\rm Int_{5}}|\leq C\|F\|_{L_{\tau\xi\mu}^{2}} \prod\limits_{j=1}^{2}\|F_{j}\|_{L_{\tau\xi\mu}^{2}} .
\end{align*}

When (\ref{4.018}) is valid, since  $-\frac{9}{8}+16\epsilon\leq s_{1}<0$ and
 $$\langle \sigma\rangle^{-\frac{1}{2}+2\epsilon}\langle \sigma_{1}\rangle^{-\frac{1}{2}-\epsilon}
 \leq \langle \sigma\rangle^{-\frac{1}{2}-\epsilon}\langle \sigma_{1}\rangle^{-\frac{1}{2}+2\epsilon},$$
we have
\begin{align}
K_{1}(\xi_{1},\mu_{1},\tau_{1},\xi,\mu,\tau)
&\leq\frac{|\xi|\langle\xi\rangle^{s_{1}}}
{\langle\sigma\rangle^{\frac{1}{2}-2\epsilon}
\prod\limits_{j=1}^{2}\langle\xi_{j}\rangle^{s_{1}}\langle\sigma_{j}\rangle^{\frac{1}{2}+\epsilon}}
\leq C\frac{|\xi|^{s_{1}+\frac{1}{2}-2\epsilon}|\xi_{2}|^{-2s_{1}-2+8\epsilon}}
{\langle\sigma_{2}\rangle^{\frac{1}{2}+\epsilon}\langle\sigma\rangle^{\frac{1}{2}+\epsilon}}
\nonumber\\
&\leq C\frac{|\xi|^{-\frac{5}{8}+14\epsilon}|\xi_{2}|^{\frac{1}{4}-24\epsilon}}
{\langle\sigma_{2}\rangle^{\frac{1}{2}+\epsilon}\langle\sigma\rangle^{\frac{1}{2}+\epsilon}}
\leq C\frac{|\xi|^{-\frac{1}{2}}|\xi_{2}|}
{\langle\sigma\rangle^{\frac{1}{2}+\epsilon}\langle\sigma_{2}\rangle^{\frac{1}{2}+\epsilon}}
.\label{4.036}
\end{align}
Thus, combining (\ref{2.019}) with (\ref{4.036}),  we have
\begin{align*}
|{\rm Int_{5}}|\leq C\|F\|_{L_{\tau\xi\mu}^{2}}
\prod\limits_{j=1}^{2}\|F_{j}\|_{L_{\tau\xi\mu}^{2}}.
\end{align*}

When (\ref{4.019}) is valid, this case can be proved similarly to (\ref{4.018}) with the aid of (\ref{2.020}).

 When (\ref{4.016}) is valid, consider (\ref{4.023}), (\ref{4.024}),  respectively.

 We dyadically decompose the spectra  as
\begin{align*}
\langle\sigma\rangle\sim 2^{j},\quad
\langle\sigma_{1}\rangle\sim 2^{j_{1}},\quad
\langle\sigma_{2}\rangle\sim 2^{j_{2}},\quad
|\xi|\sim 2^{k},\quad
|\xi_{1}|\sim 2^{k_{1}},|\xi_{2}|\sim 2^{k_{2}}.
\end{align*}
We define
\begin{align*}
&f_{k_{m},j_m}:=\left|\eta_{k_{m}}(\xi_{m})\eta_{j_{m}}(\sigma_{m})F_{l}(\xi_{m},\mu_{m},\tau_{m})\right|(m=1,2),\\
&f_{k,j}:=\left|\eta_{k}(\xi)\eta_{j}(\sigma)F(\xi,\mu,\tau)\right|,dV=d\xi_{1}d\mu_{1}d\tau_{1}d\xi d\mu d\tau.
\end{align*}

When (\ref{4.023}) is valid, we use Lemma 3.2 to deal with this case.
Thus, by   Lemma 3.2, since $-\frac{9}{8}+16\epsilon\leq s_{1}<0,$ we have
\begin{align}
 {\rm Int_{5}}
 &\leq  C\sum\limits_{k_{1},k_{2}>0,\>k}
 \sum\limits_{j_{1},j_{2},j\geq0}
2^{-j(\frac{1}{2}-2\epsilon)
-(j_{1}+j_{2})(\frac{1}{2}+\epsilon)-2k_{2}s_{1}
+k(1+s_{1})}
 \int_{\SR^{6}}f_{k,j}\prod_{m=1}^{2}f_{k_{m},j_{m}}dV \nonumber\\
 & \leq C
\sum\limits_{k_{1},k_{2}>0,\>k}\sum\limits_{j_{1},j_{2},j\geq0}
2^{2j\epsilon-(j_{1}+j_{2})\epsilon-\frac{j_{\rm max}}{2}-k_{2}(2s_{1}
+\frac{1}{4})+k\left(\frac{3}{4}+s_{1}\right)}
\|f_{k,j}\|_{L^{2}} \prod_{m=1}^{2}\|f_{k_{m},j_{m}}\|_{L^{2}} .
\label{4.037}
\end{align}
When $j=j_{\rm max},$ from (\ref{4.037}), if $\frac{9}{8}+16\epsilon\leq s_{1}<-\frac{3}{4},$ we have
\begin{align}
{\rm Int_{5}}
&\leq  C\sum\limits_{k_{1},k_{2}>0,\>k}
\sum\limits_{j_{1},j_{2},j\geq0}
2^{2j\epsilon-(j_{1}+j_{2})\epsilon-\frac{j_{\rm max}}{2}-k_{2}(2s_{1}+\frac{1}{4})
+k\left(\frac{3}{4}+s_{1}\right)}
\|f_{k,j}\|_{L^{2}}
\prod_{m=1}^{2}\|f_{k_{m},j_{m}}\|_{L^{2}}
\nonumber\\   &\leq C
\sum\limits_{k_{1},k_{2}>0,\>k}
2^{-k_{2}(2s_{1}+\frac{9}{4}-8\epsilon)
+k(\frac{1}{4}+2\epsilon)}\|f\|_{L^{2}}
\prod_{m=1}^{2}\|f_{m}\|_{L^{2}}
\nonumber\\
&\leq C
\|f\|_{L^{2}}
\prod_{m=1}^{2}\|f_{m}\|_{L^{2}}
.\label{4.038}
\end{align}
When $j=j_{\rm max},$ from (\ref{4.037}), if $-\frac{3}{4}\leq s_{1}<0,$ we have
\begin{align}
{\rm Int_{5}}
&\leq  C\sum\limits_{k_{1},k_{2}>0,\>k}
\sum\limits_{j_{1},j_{2},j\geq0}
2^{2j\epsilon-(j_{1}+j_{2})\epsilon
-\frac{j_{\rm max}}{2}-k_{2}(2s_{1}+\frac{1}{4})
+k\left(\frac{3}{4}+s_{1}\right)}
\|f_{k,j}\|_{L^{2}}
\prod_{m=1}^{2}\|f_{k_{m},j_{m}}\|_{L^{2}}
\nonumber\\
&\leq C
\sum\limits_{k_{1},k_{2}>0}
2^{-k_{2}(s_{1}+\frac{3}{2}-8\epsilon)}
\|f\|_{L^{2}}
\prod_{m=1}^{2}\|f_{m}\|_{L^{2}}
\nonumber\\
&\leq C
\|f\|_{L^{2}}
\prod_{m=1}^{2}\|f_{m}\|_{L^{2}}
.\label{4.039}
\end{align}
When $j_{1}=j_{\rm max},$ from (\ref{4.038}), if $-\frac{9}{8}+16\epsilon\leq s_{1}<-\frac{3}{4},$ we have
\begin{align}
{\rm Int_{5}}
&\leq   C\sum\limits_{k_{1},k_{2}>0,\>k}
\sum\limits_{j_{1},j_{2},j\geq0}
2^{2j\epsilon-(j_{1}+j_{2})\epsilon
-\frac{j_{\rm max}}{2}-k_{2}(2s_{1}
+\frac{1}{4})+k(\frac{3}{4}+s_{1})}
\|f_{k,j}\|_{L^{2}}
\prod_{m=1}^{2}\|f_{k_{m},j_{m}}\|_{L^{2}}
\nonumber\\
&\leq C
\sum\limits_{k_{1},k_{2}>0,\>k}
2^{-k_{2}(2s_{1}+\frac{9}{4}-4\epsilon)}
\|f\|_{L^{2}}
\prod_{m=1}^{2}\|f_{m}\|_{L^{2}}
\nonumber\\
&\leq C
\|f\|_{L^{2}}
\prod_{m=1}^{2}\|f_{m}\|_{L^{2}}
.\label{4.040}
\end{align}
When $j_{2}=j_{\rm max},$ from (\ref{4.038}), if $-\frac{3}{4}\leq s_{1}<0,$ we have
\begin{align}
{\rm Int_{5}}
&\leq  C\sum\limits_{k_{1},k_{2}>0,\>k}
\sum\limits_{j_{1},j_{2},j\geq0}
2^{2j\epsilon-(j_{1}+j_{2})\epsilon
-\frac{j_{\rm max}}{2}-k_{2}(2s_{1}+\frac{1}{4})
+k\left(\frac{3}{4}+s_{1}\right)}
\|f_{k,j}\|_{L^{2}}
\prod_{m=1}^{2}\|f_{k_{m},j_{m}}\|_{L^{2}}
\nonumber\\
&\leq C
\sum\limits_{k_{1},k_{2}>0}
2^{-k_{2}(s_{1}+\frac{3}{2}-4\epsilon)}
\|f\|_{L^{2}}
\prod_{m=1}^{2}\|f_{m}\|_{L^{2}}
\nonumber\\
&\leq C
\|f\|_{L^{2}}
\prod_{m=1}^{2}\|f_{m}\|_{L^{2}}.\label{4.041}
\end{align}
When $j_{2}=j_{\rm max},$ this case can be proved similarly to case $j_{1}=j_{\rm max}.$

When (\ref{4.024}) is valid,  by using (1) of  Lemma 3.1, since $-\frac{9}{8}+16\epsilon\leq s_{1}<0,$ we have
\begin{align}
{\rm Int_{5}}
&\leq  C\sum\limits_{k,k_{1},k_{2}>0}
\sum\limits_{j_{1},j_{2},j\geq0}
2^{-j(\frac{1}{2}-2\epsilon)
-(j_{1}+j_{2})(\frac{1}{2}+\epsilon)
-2k_{2}s_{1}+k(1+s_{1})}
\int_{\SR^{6}}f_{k,j}\prod_{m=1}^{2}f_{k_{m},j_{m}}dV\nonumber\\   &\leq C
\sum\limits_{k,k_{1},k_{2}>0}
\sum\limits_{j_{1},j_{2},j\geq0}
2^{2j\epsilon-(j_{1}+j_{2})\epsilon
-k_{2}(2s_{1}+\frac{7}{2})+k(\frac{3}{2}+s_{1})}
\|f_{k,j}\|_{L^{2}}
\prod_{m=1}^{2}\|f_{k_{m},j_{m}}\|_{L^{2}}
\nonumber\\
&\leq C
\sum\limits_{k,k_{1},k_{2}>0}
2^{-k_{2}(s_{1}+\frac{7}{2}-8\epsilon)
+k(\frac{3}{2}+2\epsilon)}
\|f\|_{L^{2}}
\prod_{m=1}^{2}\|f_{m}\|_{L^{2}}\nonumber\\
&\leq C
\sum\limits_{k,k_{1},k_{2}>0}
2^{-k_{2}(s_{1}+2-10\epsilon)}
\|f\|_{L^{2}}\prod_{m=1}^{2}\|f_{m}\|_{L^{2}}\nonumber\\
&\leq C
\|f\|_{L^{2}}\prod_{m=1}^{2}\|f_{m}\|_{L^{2}}.\label{4.042}
\end{align}

(VI) Region $\Omega_{6}.$
In this region, we consider (\ref{4.015}),  (\ref{4.016}),  respectively.

  When (\ref{4.015}) is valid,  one of (\ref{4.017})-(\ref{4.019}) must occur.

 When (\ref{4.017}) is valid, since  $-\frac{9}{8}+16\epsilon\leq s_{1}<0$, we have
\begin{align}
K_{1}(\xi_{1},\mu_{1},\tau_{1},\xi,\mu,\tau)
&\leq C\frac{|\xi|\langle\xi\rangle^{s_{1}}}
{\langle\sigma\rangle^{\frac{1}{2}-2\epsilon}
\prod\limits_{j=1}^{2}\langle\xi_{j}\rangle^{s_{1}}\langle\sigma_{j}\rangle^{\frac{1}{2}+\epsilon}}
\leq C\frac{|\xi_{2}|^{-\frac{3}{2}-s_{1}+10\epsilon}}
{\prod\limits_{j=1}^{2}\langle\sigma_{j}\rangle^{\frac{1}{2}+\epsilon}}
\nonumber\\
&\leq C\frac{|\xi_{2}|^{-\frac{3}{8}}}
{\prod\limits_{j=1}^{2}\langle\sigma_{j}\rangle^{\frac{1}{2}+\epsilon}}\leq
C\frac{|\xi_{1}|^{-\frac{1}{2}}|\xi_{2}|}
{\prod\limits_{j=1}^{2}\langle\sigma_{j}\rangle^{\frac{1}{2}+\epsilon}}
.\label{4.043}
\end{align}
Thus, combining (\ref{2.017}) with (\ref{4.043}),  we have
\begin{align*}
|{\rm Int_{6}}|\leq C\|F\|_{L_{\tau\xi\mu}^{2}}
\prod\limits_{j=1}^{2}\|F_{j}\|_{L_{\tau\xi\mu}^{2}}.
\end{align*}

When (\ref{4.018}) is valid, since
$-\frac{9}{8}+16\epsilon\leq s_{1}<0$ and
$$
\langle \sigma\rangle^{-\frac{1}{2}+2\epsilon}\langle \sigma_{1}\rangle^{-\frac{1}{2}-\epsilon}
\leq \langle \sigma\rangle^{-\frac{1}{2}-\epsilon}\langle \sigma_{1}\rangle^{-\frac{1}{2}+2\epsilon},
$$
we have
\begin{align}
K_{1}(\xi_{1},\mu_{1},\tau_{1},\xi,\mu,\tau)
&\leq C\frac{|\xi|\langle\xi\rangle^{s_{1}}}
{\langle\sigma\rangle^{\frac{1}{2}-2\epsilon}
\prod\limits_{j=1}^{2}\langle\xi_{j}\rangle^{s_{1}}\langle\sigma_{j}\rangle^{\frac{1}{2}+\epsilon}}
\leq C\frac{|\xi_{2}|^{-\frac{3}{2}-s_{1}+10\epsilon}}
{\langle\sigma\rangle^{\frac{1}{2}+\epsilon}\langle\sigma_{2}\rangle^{\frac{1}{2}+\epsilon}}
\nonumber\\
&\leq
C\frac{|\xi_{2}|^{-\frac{3}{8}}}
{\prod\limits_{j=1}^{2}\langle\sigma_{j}\rangle^{\frac{1}{2}+\epsilon}}
\leq C\frac{|\xi|^{-\frac{1}{2}}|\xi_{2}|}
{\langle\sigma\rangle^{\frac{1}{2}+\epsilon}\langle\sigma_{2}\rangle^{\frac{1}{2}+\epsilon}}
.\label{4.044}
\end{align}
Thus, combining (\ref{2.019}) with (\ref{4.044}),  we have
\begin{align*}
{\rm |Int_{6}|}\leq C\|F\|_{L_{\tau\xi\mu}^{2}}
\prod\limits_{j=1}^{2}\|F_{j}\|_{L_{\tau\xi\mu}^{2}}.
\end{align*}

When (\ref{4.019}) is valid, this case can be proved similarly to (\ref{4.018}) with the aid of (\ref{4.020}).

  When (\ref{4.016}) is valid, this case can be proved similarly to case (\ref{4.016}) of Region $\Omega_{5}.$

  (VII) Region $\Omega_{7}.$ This case can be proved similarly to Region $\Omega_{6}.$

This completes the proof of Lemma 4.1.

\noindent{\bf Remark 3.} In the case (\ref{4.023})  of  Region $\Omega_4$ it leads to the  requirement $-\frac{9}{8} < s_{1}<0.$

 \begin{Lemma}\label{Lemma4.2}
Let $-1+10\epsilon\leq s<0$.
Then, we have
\begin{align}
\|\partial_{x}
  \left[I_{N}(u_{1}u_{2})-I_{N}u_{1}I_{N}u_{2}
   \right]
\|_{X_{-\frac{1}{2}+2\epsilon}^{0,0}}
 \leq CN^{-2+10\epsilon}
  \prod_{j=1}^{2}\|I_{N}u_{j}
   \|_{X_{\frac{1}{2}+\epsilon}^{0,0}} .
  \label{4.045}
\end{align}
\end{Lemma}
\noindent{\bf Proof.}
To prove (\ref{4.045}),  by duality, it suffices to  prove that
\begin{align}
&\left|\int_{\SR^{3}}\bar{h}\partial_{x}
  \left[I_{N}(u_{1}u_{2})-I_{N}u_{1}I_{N}u_{2}\right]dxdydt\right|
   \nonumber\\
\leq&  CN^{-2+10\epsilon}
     \|h\|_{X_{\frac{1}{2}-2\epsilon}^{0,0}}
      \prod_{j=1}^{2}
       \|I_{N}u_{j}\|_{X_{\frac{1}{2}+\epsilon}^{0,0}}.
        \label{4.046}
\end{align}
for $h\in X_{\frac{1}{2}-2\epsilon}^{0,0}.$
Let
\begin{align}
&F(\xi,\mu,\tau)=
\langle \sigma\rangle^{\frac{1}{2}-2\epsilon}M(\xi)\mathscr{F}h(\xi,\mu,\tau),\nonumber\\   &
F_{j}(\xi_{j},\mu_{j},\tau_{j})=M(\xi_{j})
\langle \sigma_{j}\rangle^{\frac{1}{2}+\epsilon}
\mathscr{F}u_{j}(\xi_{j},\mu,\tau_{j})\qquad(j=1,2).
 \label{4.047}
\end{align}
To obtain (\ref{4.046}), from (\ref{4.047}), it suffices to prove
\begin{align}
&\int_{\SR^2}\!\int_{\
 \begin{array}{l}
   \xi=\xi_1+\xi_2,\\
     \mu=\mu_1+\mu_2,\\
   \tau = \tau_1+\tau_2
   \end{array}}
   \frac{|\xi|G(\xi_{1},\xi_{2})F(\xi,\mu,\tau)
    \prod\limits_{j=1}^{2}F_{j}(\xi_{j},\mu_{j},\tau_{j})}
    {\langle\sigma_{j}\rangle^{\frac{1}{2}-2\epsilon}
     \prod\limits_{j=1}^{2}\langle\sigma_{j}\rangle^{\frac{1}{2}+\epsilon}}
      d\xi_{1}d\mu_{1}d\tau_{1}d\xi d\mu d\tau\nonumber\\
\leq & CN^{-1+10\epsilon}
     \|F\|_{L_{\xi\mu\tau}^{2}}
      \prod_{j=1}^{2}\|F_{j}\|_{L_{\xi\mu\tau}^{2}},
       \label{4.048}
\end{align}
where
\begin{align*}
G(\xi_{1},\xi_{2})=\frac{M(\xi_{1})M(\xi_{2})-M(\xi)}{M(\xi_{1})M(\xi_{2})}.
\end{align*}
Without loss of generality,  we assume that
 $F(\xi,\mu,\tau)\geq 0,F_j(\xi_{j},\mu_{j},\tau_{j})\geq 0(j=1,2)$.
By symmetry, we can assume that $|\xi_{1}|\geq |\xi_{2}|.$

 \noindent
We define
\begin{align*}
&\hspace{-0.8cm}A_1=\left\{(\xi_1,\mu_{1},\tau_1,\xi,\mu,\tau)\in D^{*},
  |\xi_{2}|\leq |\xi_{1}|\leq \frac{N}{2}\right\},\\
&\hspace{-0.8cm} A_2=\{ (\xi_1,\mu_{1},\tau_1,\xi,\mu,\tau)\in D^{*},
|\xi_1|>\frac{N}{2},|\xi_{1}|\geq |\xi_{2}|, |\xi_{2}|\leq 2A\},\\
&\hspace{-0.8cm} A_3=\{ (\xi_1,\mu_{1},\tau_1,\xi,\mu,\tau)\in D^{*},
|\xi_1|>\frac{N}{2},|\xi_{1}|\geq |\xi_{2}|, 2A<|\xi_{2}|\leq N\},\\
&\hspace{-0.8cm}A_4=\{(\xi_1,\mu_{1},\tau_1,\xi,\mu,\tau)\in D^{*},
|\xi_{1}|>\frac{N}{2},|\xi_{1}|\geq |\xi_{2}|,|\xi_{2}|>N\}.
\end{align*}
Here $D^{*}$ is defined as in Lemma 3.1.
Obviously, $D^{*}\subset\bigcup\limits_{j=1}^{4}A_{j}.$
We define
\begin{equation}
    K_{2}(\xi_{1},\mu_{1},\tau_{1},\xi,\mu,\tau):=\frac{|\xi|G(\xi_{1},\xi_{2})}
    {\langle\sigma_{j}\rangle^{\frac{1}{2}-2\epsilon}
\prod\limits_{j=1}^{2}\langle\sigma_{j}\rangle^{\frac{1}{2}+\epsilon}}\label{4.049}
\end{equation}
and
\begin{align*}
J_{k}:=\int_{A_{j}} K_{2}(\xi_{1},\mu_{1},\tau_{1},\xi,\mu,\tau)F(\xi,\mu,\tau)
\prod_{j=1}^{2}F_{j}(\xi_{j},\mu_{j},\tau_{j})
d\xi_{1}d\mu_{1}d\tau_{1}d\xi d\mu d\tau,
\end{align*}
$1\leq k\leq 4, k\in\mathbf{ N}.$

 We consider (\ref{4.015}) and (\ref{4.016}), respectively.

 When (\ref{4.015}) is valid,
one of (\ref{4.017})-(\ref{4.019}) must occur,
from \cite[Lemma 1.4]{ILM-CPAA}, we have
\begin{align*}
\sum\limits_{k=1}^{4}{\rm J}_{k}\leq CN^{-(2-10\epsilon)}
\|F\|_{L_{\tau\xi\mu}^{2}}
\prod_{j=1}^{2}\|F_{j}\|_{L_{\tau\xi\mu}^{2}}
.
\end{align*}
Thus, we only consider the case  (\ref{4.016}).

Now we consider the integrals over the above four regions one by one.

 (I) Region $A_{1}.$ In this case, since $M(\xi_{1},\xi_{2})=0$, thus we have  ${\rm J}_{1}=0$.

 (II) Region $A_{2}$. From \cite[Page 902]{ILM-CPAA}, we have
\begin{align}
G(\xi_{1},\xi_{2})\leq C\frac{|\xi_{2}|}{|\xi_{1}|}\label{4.050}.
\end{align}
Inserting (\ref{4.046}) into (\ref{4.047}) yields
\begin{align}
    K_{2}(\xi_{1},\mu_{1},\tau_{1},\xi,\mu,\tau)\leq C\frac{|\xi|G(\xi_{1},\xi_{2})}
    {\langle\sigma\rangle^{\frac{1}{2}-2\epsilon}
\prod\limits_{j=1}^{2}\langle\sigma_{j}\rangle^{\frac{1}{2}+\epsilon}}\leq
\frac{C|\xi_{2}|}{\langle\sigma\rangle^{\frac{1}{2}-2\epsilon}
\prod\limits_{j=1}^{2}\langle\sigma_{j}\rangle^{\frac{1}{2}+\epsilon}}\label{4.051}.
\end{align}
We dyadically decompose the spectra as
\begin{align*}
\langle\sigma\rangle\sim 2^{j},\quad
\langle\sigma_{1}\rangle\sim 2^{j_{1}},\quad
\langle\sigma_{2}\rangle\sim 2^{j_{2}},\quad
|\xi|\sim 2^{k},\quad
|\xi_{1}|\sim 2^{k_{1}},\quad
|\xi_{2}|\sim 2^{k_{2}}.
\end{align*}
We define
\begin{align*}
&f_{k_m, j_{m}}:=\eta_{k_{m}}(\xi_{m})\eta_{j_{m}}(\sigma_{m})F_{j}(\xi_{m},\mu_{m},\tau_{m})(m=1,2),\\
&f_{k, j}:=\eta_{k}(\xi)\eta_{j}(\sigma)\left|F(\xi,\mu,\tau))\right|.
\end{align*}
Thus, by using (\ref{2.04}),  we have
\begin{align}
&\hspace{-2.5cm}{\rm J}_{2}\leq C\sum\limits_{k_{1},k_{2}>0,\>k}
\sum\limits_{j_{1},j_{2},j\geq0}2^{-j(\frac{1}{2}-2\epsilon)-(j_{1}+j_{2})(\frac{1}{2}+\epsilon)+k_{2}}
 \int_{\SR^{6}}
f_{k,j}\prod\limits_{m=1}^{2}f_{k_{m},j_{m}}dV .\label{4.052}
\end{align}
In this case, we consider (\ref{4.023}), (\ref{4.024}), respectively.

\noindent When (\ref{4.023}) is valid,  we consider  $j=j_{\rm max},j_{1}=j_{\rm max},j_{2}=j_{\rm max}$, respectively.

\noindent When $j=j_{\rm max}$ is valid, from (\ref{4.052}), we have
\begin{align}
{\rm J}_{2}
&\leq C\sum\limits_{k_{1},k_{2}>0,\>k}
      \sum\limits_{j_{1},j_{2}}
      2^{- k(\frac{9}{4}-8\epsilon)
      -(j_{1}+j_{2})\epsilon
      +k_{2}(\frac{1}{4}+2\epsilon)}\|f\|_{L^{2}}
\prod\limits_{j=1}^{2}\|f_{j}\|_{L^{2}} \nonumber\\
&\leq CN^{-(\frac{9}{4}-8\epsilon)}\|f\|_{L^{2}}
\prod\limits_{j=1}^{2}\|f_{j}\|_{L^{2}}.\label{4.053}
\end{align}
 When $j_{1}=j_{\rm max}$ is valid, from (\ref{4.052}), we have
\begin{align}
     {\rm J}_{2}
&\leq C\sum\limits_{k_{1},k_{2}>0,\>k}
       \sum\limits_{j_{1},j_{2},j\geq0}
       2^{2j\epsilon-(j_{1}+j_{2})\epsilon-\frac{j_{1}}{2}
       +\frac{3}{4}k_{2}}\|f\|_{L^{2}} \prod\limits_{j=1}^{2}\|f_{j}\|_{L^{2}}\nonumber\\
&\leq C\sum\limits_{k_{1},k_{2}>0,\>k}
       \sum\limits_{j_{1},j_{2}\geq0}
       2^{-j_{2}\epsilon-j_{1}(\frac{1}{2}-\epsilon)
       -\frac{k_{1}}{4}+\frac{3}{4}k_{2}}\|f\|_{L^{2}}
       \prod\limits_{j=1}^{2}\|f_{j}\|_{L^{2}}\nonumber\\
&\leq C\sum\limits_{k_{1},k_{2}>0,\>k}
      2^{-(\frac{9}{4}-4\epsilon)k_{1}
      +k_{2}(\frac{1}{2}+\epsilon)}\|f\|_{L^{2}}
      \prod\limits_{j=1}^{2}\|f_{j}\|_{L^{2}}\nonumber\\
&\leq CN^{-(\frac{9}{4}-4\epsilon)}
      \|f\|_{L^{2}}
      \prod\limits_{j=1}^{2}\|f_{j}\|_{L^{2}}.\label{4.054}
\end{align}
 When $j_{2}=j_{\rm max}$ is valid, this case can be proved similarly to $j_{1}=j_{\rm max}$ of Region $A_{2}$.

 When (\ref{4.024}) is valid, by using (2) of Lemma 3.1, we have
 \begin{align}
{\rm J}_{2}
&\leq C\sum\limits_{k_{1},k_{2}>0,\>k}
\sum\limits_{j_{1},j_{2},j\geq0}2^{2j\epsilon-(j_{1}
  +j_{2})\epsilon-2k_{1}+k_{2}\epsilon}\|f\|_{L^{2}} \prod\limits_{j=1}^{2}\|f_{j}\|_{L^{2}}
  \nonumber\\
&\leq C\sum\limits_{k_{1},k_{2}>0,\>k}
  \sum\limits_{j_{1},j_{2}\geq0}2^{-(j_{1}+j_{2})\epsilon
   -(2-8\epsilon)k_{1}+3k_{2}\epsilon}\|f\|_{L^{2}}
    \prod\limits_{j=1}^{2}\|f_{j}\|_{L^{2}} \nonumber\\
&\leq C\sum\limits_{k_{1},k_{2}>0,\>k}
   2^{-(\frac{7}{2}-8\epsilon)k_{1}
    +k_{2}(\frac{1}{2}+2\epsilon)}\|f\|_{L^{2}}
    \prod\limits_{j=1}^{2}\|f_{j}\|_{L^{2}}\nonumber\\
&\leq CN^{-(\frac{7}{2}-8\epsilon)}\|f\|_{L^{2}}
      \prod\limits_{j=1}^{2}\|f_{j}\|_{L^{2}}.\label{4.055}
\end{align}

(III) Region $A_{3}$.
From of \cite[Page 902]{ILM-CPAA}, we know that  (\ref{4.050}) is valid.
Combining (\ref{4.050}) with (\ref{4.049}), we have
\begin{equation}
    K_{2}(\xi_{1},\mu_{1},\tau_{1},\xi,\mu,\tau)\leq C\frac{{\rm min}
    \left\{|\xi|,|\xi_{1}|,|\xi_{2}|\right\}}{\langle\sigma\rangle^{\frac{1}{2}-2\epsilon}
\prod\limits_{j=1}^{2}\langle\sigma_{j}\rangle^{\frac{1}{2}+\epsilon}}.\label{4.056}
\end{equation}
We dyadically decompose the spectra as
\begin{align*}
\langle\sigma\rangle\sim 2^{j},\quad
\langle\sigma_{1}\rangle\sim 2^{j_{1}},\quad
\langle\sigma_{2}\rangle\sim 2^{j_{2}},\quad
|\xi|\sim 2^{k},|\xi_{1}|\sim 2^{k_{1}},\quad
|\xi_{2}|\sim 2^{k_{2}}.
\end{align*}
We define
\begin{align*}
&f_{k_{m},j_m}:
  =\eta_{k_{m}}(\xi_{m})\eta_{j_{m}}(\sigma_{k})
    F_{m}(\xi_{m},\mu_{m},\tau_{m})\qquad (m=1,2),\\
&f_{k,j}:=\eta_{k}(\xi)\eta_{j}(\sigma)\left|F(\xi,\mu,\tau))\right|.
\end{align*}
Thus, we have
\begin{align}
  {\rm J}_{3}\leq C\sum\limits_{m_{1},m_{2}>0,\>m}
  \sum\limits_{j_{1},j_{2},j\geq0}
& 2^{-j(\frac{1}{2}-2\epsilon)
-(j_{1}+j_{2})(\frac{1}{2}+\epsilon)
+{\rm min}\{k,k_{1},k_{2}+k_{2}\}}\nonumber\\
  & \times  \int_{\SR^{6}}
f_{m,j}\prod\limits_{m=1}^{2}f_{k_{m},j_{m}}dV.\label{4.057}
\end{align}
In this case, we consider (\ref{4.023}), (\ref{4.024}), respectively.

 When (\ref{4.023}) is valid, by using Lemma 3.2,  we have
\begin{align}
    {\rm J_{3}}
\leq  C\sum\limits_{k,k_{1},k_{2}>0}
      \sum\limits_{j_{1},j_{2},j\geq0}
& 2^{-j(\frac{1}{2}-2\epsilon)-(j_{1}+j_{2})\epsilon
    +\frac{j}{2}-\frac{j_{\rm max}}{2} -\frac{k_{1}}{4}
    +\frac{3}{4}{\rm min}\{k,k_{1},k_{2}\}} \nonumber\\
& \times \int_{\SR^{6}}f_{k,j}
    \prod_{m=1}^{2}f_{k_{m},j_{m}}dV.
\label{4.058}
\end{align}
When $j=j_{\rm max},$ from (\ref{4.058}),  we have
\begin{align}
{\rm J_{3}}
&\leq  C\sum\limits_{k_{1},k_{2}>0,\>k}
       \sum\limits_{j_{1},j_{2},j\geq0}
       2^{-j(\frac{1}{2}-2\epsilon)
       -(j_{1}+j_{2})\epsilon
      -\frac{k_{1}}{4}
      +\frac{3}{4}{\rm min}\{k,k_{1},k_{2}\}}
      \|f_{k,j}\|_{L^{2}}
      \prod_{m=1}^{2}\|f_{k_{m},j_{m}}\|_{L^{2}}\nonumber\\
      &\leq C\sum\limits_{k_{1},k_{2}>0,\>k}
       2^{-k_{2}(\frac{9}{4}-8\epsilon)
    +{\rm min}\{k,k_{1},k_{2}\}
     (\frac{1}{4}+2\epsilon)} \|f\|_{L^{2}}
       \prod_{m=1}^{2}\|f_{m}\|_{L^{2}}\nonumber\\
&\leq CN^{-2+10\epsilon}\|f\|_{L^{2}}
      \prod_{m=1}^{2}\|f_{m}\|_{L^{2}}.\label{4.059}
\end{align}
When $j_{1}=j_{\rm max},$ from (\ref{4.059}),  we have
\begin{align}
       {\rm J_{3}}
&\leq  C\sum\limits_{k_{1},k_{2}>0,\>k}
       \sum\limits_{j_{1},j_{2},j\geq0}
       2^{2j\epsilon-(j_{1}+j_{2})\epsilon
        -\frac{j_{1}}{2}+\frac{k_{1}}{4}
        +\frac{3}{4}{\rm min}\{k,k_{1},k_{2}\}}
        \|f_{k,j}\|_{L^{2}}
        \prod_{m=1}^{2}\|f_{k_{m},j_{m}}\|_{L^{2}}\nonumber\\
&\leq C\sum\limits_{k_{1},k_{2}>0,\>k}
            2^{-k_{2}(\frac{9}{4}-4\epsilon)
            +{\rm min}\{k,k_{1},k_{2}\}
            (\frac{1}{4}+\epsilon)}\|f\|_{L^{2}}
             \prod_{m=1}^{2}\|f_{m}\|_{L^{2}} \nonumber\\
&\leq CN^{-2+5\epsilon}\|f\|_{L^{2}}
       \prod_{m=1}^{2}\|f_{m}\|_{L^{2}}.\label{4.060}
\end{align}
When $j_{2}=j_{\rm max},$ this case can be proved similarly to the case $j_{1}=j_{\rm max}.$

 When (\ref{4.024}) is valid,   combining (2) of  Lemma 3.1 with (\ref{4.060}), we have
\begin{align}
      {\rm J_{3}}
&\leq  C\sum\limits_{k,k_{1},k_{2}>0}
      \sum\limits_{j_{1},j_{2},j\geq0}
    2^{-j(\frac{1}{2}-2\epsilon)
   -(j_{1}+j_{2})(\frac{1}{2}+\epsilon)
   +{\rm min}\{k,k_{1},k_{2}\}}
   \int_{\SR^{6}}f_{k,j}
   \prod_{m=1}^{2}f_{k_{m},j_{m}}dV \nonumber\\
 &\leq C\sum\limits_{k,k_{1},k_{2}>0}
        \sum\limits_{j_{1},j_{2},j\geq0}2^{2j\epsilon
        -(j_{1}+j_{2})\epsilon
        -2k_{1}}\|f_{k,j}\|_{L^{2}}
        \prod_{m=1}^{2}\|f_{k_{m},j_{m}}\|_{L^{2}}\nonumber\\
&\leq C\sum\limits_{k,k_{1},k_{2}>0}
    2^{-k_{1}(2-10\epsilon)}\|f\|_{L^{2}}
    \prod_{m=1}^{2}\|f_{m}\|_{L^{2}}\nonumber\\
&\leq CN^{-2+11\epsilon}\|f\|_{L^{2}}
     \prod_{m=1}^{2}\|f_{m}\|_{L^{2}}.\label{4.061}
\end{align}

 (IV) Region $A_{4}$. In this case, we have
\begin{align}
     M(\xi_{1},\xi_{2})
\leq C N^{2s_{1}} \left(|\xi_{1}|\> |\xi_{2}|
       \right)^{-s_{1}}.\label{4.062}
\end{align}
We dyadically decompose the spectra as
\begin{align*}
\langle\sigma\rangle\sim 2^{j},\quad
\langle\sigma_{1}\rangle\sim 2^{j_{1}},\quad
\langle\sigma_{2}\rangle\sim 2^{j_{2}},\quad
|\xi|\sim 2^{k},\quad
|\xi_{1}|\sim 2^{k_{1}},\quad
|\xi_{2}|\sim 2^{k_{2}}.
\end{align*}
 We define
\begin{align*}
&f_{k_{m},j_m}:=\eta_{k_{m}}(\xi_{m})
\eta_{j_{m}}(\sigma_{m})
F_{m}(\xi_{m},\mu_{m},\tau_{m}),\qquad(m=1,2),\\
&f_{k,j}:=\eta_{k}(\xi)\eta_{j}(\sigma)\left|F(\xi,\mu,\tau))\right|.
\end{align*}
Thus, we have
\begin{align}
{\rm J}_{4}\leq CN^{2s}\sum\limits_{m_{1},m_{2}>0,\>m}
\sum\limits_{j_{1},j_{2},j\geq0}
&2^{-j(\frac{1}{2}-2\epsilon)-(j_{1}+j_{2})
(\frac{1}{2}+\epsilon)+(k_{1}+k_{2})|s|+k}
\nonumber\\   & \quad\times
 \int_{\SR^{6}}f_{k,j}\prod_{m=1}^{2}f_{k_{m},j_{m}}dV
 .\label{4.063}
\end{align}
In this case, we consider (\ref{4.023}), (\ref{4.024}), respectively.

When (\ref{4.023}) is valid,  by using Lemma 3.2, from (\ref{4.063}), we have
\begin{align}
{\rm J}_{4}\leq CN^{2s}\sum\limits_{m_{1},m_{2}>0,\>m}
\sum\limits_{j_{1},j_{2},j\geq0}
& 2^{2j\epsilon-(j_{1}+j_{2})\epsilon
    -\frac{j_{\rm max}}{2}+(m_{1}+m_{2})|s|
    -\frac{1}{4}k_{1}
    -\frac{1}{4}{\rm min}\{k,k_{1},k_{2}\}+k}
    \nonumber\\
&\ \times \|f\|_{L^{2}}
    \prod_{m=1}^{2}\|f_{m}\|_{L^{2}}.\label{4.064}
\end{align}
In this case, we consider $k=k_{\rm min}, k_{2}=k_{\rm min},$ respectively.

  When $k=k_{\rm min}$ is valid, we consider $j=j_{\rm max},j_{1}=j_{\rm max},j_{2}=j_{\rm max},$ respectively.

When $j=j_{\rm max},$ from (\ref{4.064}), since $s_{1}\geq-1+6\epsilon,$ we have
\begin{align}
      {\rm J_{4}}
&\leq  CN^{2s_{1}}\sum\limits_{k_{1},k_{2}>0,\>k}
       \sum\limits_{j_{1},j_{2},j\geq0}
       2^{-j(\frac{1}{2}-2\epsilon)
       -(j_{1}+j_{2})\epsilon
       -k_{1}(2s_{1}+\frac{1}{4})
       +\frac{3}{4}k}\|f_{k,j}\|_{L^{2}}
 \prod_{m=1}^{2}\|f_{k_{m},j_{m}}\|_{L^{2}} \nonumber\\
&\leq C
\sum\limits_{k_{1},k_{2}>0,\>k}
2^{-k_{2}(2s_{1}+\frac{9}{4}-8\epsilon)
+k(\frac{1}{4}+2\epsilon)}\|f\|_{L^{2}}
\prod_{m=1}^{2}\|f_{m}\|_{L^{2}} \nonumber\\
&\leq C\sum\limits_{k_{1}>0}
2^{-k_{2}(2s_{1}+2-10\epsilon)}\|f\|_{L^{2}}
 \prod_{m=1}^{2}\|f_{m}\|_{L^{2}} \nonumber\\
&\leq CN^{-2+10\epsilon}
\|f\|_{L^{2}}
\prod_{m=1}^{2}\|f_{m}\|_{L^{2}}.\label{4.065}
\end{align}
 When $j_{1}=j_{\rm max}$ is valid, from (\ref{4.064}),  since $s_{1}\geq-1+6\epsilon,$  we have
\begin{align}
{\rm J_{4}}
&\leq  C\sum\limits_{k_{1},k_{2}>0,\>k}
        \sum\limits_{j_{1},j_{2}\geq0}
    2^{-(\frac{1}{2}-\epsilon)j_{1}-j_{2}\epsilon
     -(k_{1}+k_{2})s_{1}-\frac{k_{1}}{4}+\frac{3}{4}k}\|f_{k,j}\|_{L^{2}}
 \prod_{m=1}^{2}\|f_{k_{m},j_{m}}\|_{L^{2}} \nonumber\\
 &\leq C
\sum\limits_{k_{1},k_{2}>0,\>k}
2^{-k_{2}(2s_{1}+\frac{9}{4}-4\epsilon)
  +k(\frac{1}{4}+\epsilon)}\|f\|_{L^{2}}
   \prod_{m=1}^{2}\|f_{m}\|_{L^{2}} \nonumber\\
&\leq C\sum\limits_{k_{1}>0}
2^{-k_{2}(2s_{1}+2-5\epsilon)}\|f\|_{L^{2}} \prod_{m=1}^{2}\|f_{m}\|_{L^{2}} \nonumber\\
&\leq CN^{-2+11\epsilon}
\|f\|_{L^{2}} \prod_{m=1}^{2}\|f_{m}\|_{L^{2}} .\label{4.066}
\end{align}

When $j_{2}=j_{\rm max}$ is valid, this case can be proved similarly to $j_{1}=j_{\rm max}$.

When $k_{2}=k_{\rm min}$ is valid, we consider $j=j_{\rm max},j_{1}=j_{\rm max},j_{2}=j_{\rm max},$ respectively.

 When $j=j_{\rm max},$ from (\ref{4.064}), since $s_{1}\geq-1+10\epsilon,$ we have
\begin{align}
     {\rm J_{4}}
&\leq  CN^{2s_{1}}\sum\limits_{k_{1},k_{2}>0,\>k}\sum\limits_{j_{1},j_{2},j\geq0}
2^{-j(\frac{1}{2}-2\epsilon)-(j_{1}+j_{2})\epsilon-k_{1}(s_{1}-\frac{3}{4})-k_{2}(s_{1}+\frac{1}{4})}
\nonumber\\
&\qquad\qquad\qquad\qquad\qquad\qquad
\times\|f_{k,j}\|_{L^{2}}
 \prod_{m=1}^{2}\|f_{k_{m},j_{m}}\|_{L^{2}} \nonumber\\   &\leq CN^{2s_{1}}
\sum\limits_{k_{1},k_{2}>0,\>k}
2^{-k_{1}(s_{1}+\frac{5}{4}-8\epsilon)-k_{2}(s_{1}+\frac{3}{4}-2\epsilon)}\|f\|_{L^{2}} \prod_{m=1}^{2}\|f_{m}\|_{L^{2}} \nonumber\\
&\leq CN^{2s_{1}}\sum\limits_{k_{1},k_{2}>0}
2^{-k_{1}(s_{1}+1-8\epsilon)-k_{2}(s_{1}+1-2\epsilon)}\|f\|_{L^{2}} \prod_{m=1}^{2}\|f_{m}\|_{L^{2}} \nonumber\\
&\leq CN^{-2+10\epsilon}
\|f\|_{L^{2}} \prod_{m=1}^{2}\|f_{m}\|_{L^{2}} .\label{4.067}
\end{align}

 When $j_{1}=j_{\rm max}$ is valid, from (\ref{4.064}),  since $s_{1}\geq-1+10\epsilon,$  we have
\begin{align}
{\rm J_{4}}
&\leq  CN^{2s_{1}}\sum\limits_{k_{1},k_{2}>0,\>k}
\sum\limits_{j_{1},j_{2}\geq0}
2^{-(\frac{1}{2}-\epsilon)j_{1}-j_{2}\epsilon-k_{1}(s_{1}-\frac{3}{4})-k_{2}(s_{1}+\frac{1}{4})}\|f_{k,j}\|_{L^{2}}
 \prod_{m=1}^{2}\|f_{k_{m},j_{m}}\|_{L^{2}} \nonumber\\
 &\leq CN^{2s_{1}}
\sum\limits_{k_{1},k_{2}>0,\>k}
2^{-k_{1}(s_{1}+\frac{5}{4}-4\epsilon)-k_{2}(s_{1}+\frac{3}{4}+\epsilon)}\|f\|_{L^{2}} \prod_{m=1}^{2}\|f_{m}\|_{L^{2}} \nonumber\\
&\leq CN^{2s_{1}}\sum\limits_{k_{1},k_{2}>0}
2^{-k_{2}(s_{1}+1-4\epsilon)-k_{2}(s_{1}+1+\epsilon)}\|f\|_{L^{2}} \prod_{m=1}^{2}\|f_{m}\|_{L^{2}} \nonumber\\
&\leq CN^{-2+5\epsilon}
\|f\|_{L^{2}} \prod_{m=1}^{2}\|f_{m}\|_{L^{2}} .\label{4.068}
\end{align}

When $j_{2}=j_{\rm max}$ is valid, this case can be proved similarly to $j_{1}=j_{\rm max}$.

This completes the proof of Lemma 4.2.

\begin{Lemma}\label{Lem4.3}
Let $s\geq-\frac{9}{8}+16\epsilon,s_{2}\geq0$ and $u_{j}\in X_{\frac{1}{2}+\epsilon}^{s_{1},s_{2}}(j=1,2)$.
Then, we have
\begin{align}
&\|\partial_{x}I(u_{1}u_{2})\|_{X_{-\frac{1}{2}+2\epsilon}^{0,0}}\leq C
\prod_{j=1}^{2}\|Iu_{j}\|_{X_{\frac{1}{2}+\epsilon}^{0,0}}.\label{4.069}
\end{align}
\end{Lemma}
\noindent{\bf Proof.} To prove (\ref{4.069}),  by duality, it suffices to  prove that
\begin{align}
&\left|\int_{\SR^{3}}\bar{u}\partial_{x}I(u_{1}u_{2})dxdydt\right|\leq
C\|u\|_{X_{\frac{1}{2}-2\epsilon}^{0,0}}
\prod_{j=1}^{2}
\|Iu_{j}\|_{X_{\frac{1}{2}+\epsilon}^{0,0}}.\label{4.070}
\end{align}
for $u\in X_{\frac{1}{2}-2\epsilon}^{0,0}.$
Let
\begin{align}
&F(\xi,\mu,\tau)=
\langle \sigma\rangle^{\frac{1}{2}-2\epsilon}\mathscr{F}u(\xi,\mu,\tau),\nonumber\\   &
F_{j}(\xi_{j},\mu_{j},\tau_{j})=M(\xi_{j})
\langle \sigma_{j}\rangle^{\frac{1}{2}+\epsilon}
\mathscr{F}u_{j}(\xi_{j},\mu,\tau_{j})(j=1,2),\label{4.071}
\end{align}
and
\begin{align*}
D:=\left\{(\xi_1,\mu_{1},\tau_1,\xi,\mu,\tau)\in {\rm R^6},
\xi=\sum_{j=1}^{2}\xi_j,\mu=\sum_{j=1}^{2}\mu_{j},\tau=\sum_{j=1}^{2}\tau_j\right\}.
\end{align*}
To obtain (\ref{4.070}), from  (\ref{4.071}), it suffices to prove that
\begin{align}
&\int_{D}\frac{|\xi|M(\xi)
F(\xi,\mu,\tau)F_{j}(\xi_{j},\mu_{j},\tau_{j})}{\langle\sigma_{j}\rangle^{\frac{1}{2}-2\epsilon}
\prod\limits_{j=1}^{2}M(\xi_{j})\langle\sigma_{j}\rangle^{\frac{1}{2}+\epsilon}}
d\xi_{1}d\mu_{1}d\tau_{1}d\xi d\mu d\tau\nonumber\\   \leq &C
\|F\|_{L_{\tau\xi\mu}^{2}}
\prod_{j=1}^{2}\|F_{j}\|_{L_{\tau\xi\mu}^{2}}
.\label{4.072}
\end{align}
From (2.4) of \cite{IMEJDE}, we have
\begin{align}
\frac{M(\xi)}{\prod\limits_{j=1}^{2}M(\xi_{j})}\leq C\frac{\langle\xi\rangle^{s}}
{\prod\limits_{j=1}^{2}\langle\xi_{j}\rangle^{s}}\label{4.073}.
\end{align}
Inserting  (\ref{4.073}) into the left hand side of (\ref{4.072}), we have
\begin{align}
&\int_{D}\frac{|\xi|\langle\xi\rangle^{s}
F(\xi,\mu,\tau)F_{j}(\xi_{j},\mu_{j},\tau_{j})}{\langle\sigma_{j}\rangle^{\frac{1}{2}-2\epsilon}
\prod\limits_{j=1}^{2}\langle\xi_{j}\rangle^{s}\langle\sigma_{j}\rangle^{\frac{1}{2}+\epsilon}}
d\xi_{1}d\mu_{1}d\tau_{1}d\xi d\mu d\tau.\label{4.074}
\end{align}
By using (\ref{4.04}),  we have  (\ref{4.074}) can be bounded by
$ C
\|F\|_{L_{\tau\xi\mu}^{2}} \prod\limits_{j=1}^{2}\|F_{j}\|_{L_{\tau\xi\mu}^{2}}.$

This completes the proof of Lemma 4.3.

\bigskip
\bigskip
\noindent {\large\bf 5. Proof of Theorem  1.1}

\setcounter{equation}{0}

 \setcounter{Theorem}{0}

\setcounter{Lemma}{0}

\setcounter{section}{5}

This section is devoted to  proving  Theorem 1.1.

\noindent We define
\begin{align}
&\Phi_{1}(u):=\psi(t)W(t)u_{0}+\frac{1}{2}\psi\left(\frac{t}{\tau}\right)
\int_{0}^{t}W(t-\tau)\partial_{x}(u^{2})d\tau,\label{5.01}\\
&B_{1}(0,2C\|u_{0}\|_{H^{s_{1},s_{2}}}):=\left\{u:\|u\|_{X_{\frac{1}{2}
+\epsilon}^{s_{1},s_{2}}}\leq 2C\|u_{0}\|_{H^{s_{1},s_{2}}}\right\}.\label{5.02}
\end{align}
Combining  Lemmas 2.2, 4.1 with (\ref{5.01})-(\ref{5.02}),  we derive that
\begin{align}
\left\|\Phi_{1}(u)\right\|_{X_{b}^{s_{1},s_{2}}}
\leq & \left\|\eta(t)W(t)u_{0}\right\|_{X_{\frac{1}{2}+\epsilon}^{s_{1},s_{2}}}
+\left\|\frac{1}{2}\eta\left(\frac{t}{\tau}\right)\int_{0}^{t}W(t-\tau)\partial_{x}(u^{2})d\tau\right\|_{X_{b}^{s_{1},s_{2}}}\nonumber\\
\leq& C\|u_{0}\|_{H^{s_{1},s_{2}}}+CT^{\epsilon}\left\|\partial_{x}(u^{2})\right\|_{X_{-\frac{1}{2}+2\epsilon}^{s_{1},s_{2}}}\nonumber\\
\leq & C\|u_{0}\|_{H^{s_{1},s_{2}}}+CT^{\epsilon}\left\|u\right\|_{X_{b}^{s_{1},s_{2}}}^{2}\nonumber\\
\leq & C\|u_{0}\|_{H^{s_{1},s_{2}}}+4C^{3}T^{\epsilon}\left\|u_{0}\right\|_{H^{s_{1},s_{2}}}^{2}.\label{5.03}
\end{align}
We define
\begin{align}
T^{\epsilon}:=\left[16C^{2}(\|u_{0}\|_{H^{s_{1},s_{2}}}+1)\right]^{-1}.\label{5.04}
\end{align}
From (\ref{5.03})-(\ref{5.04}),  we have
\begin{align}
&\left\|\Phi_{1}(u)\right\|_{X_{\frac{1}{2}+\epsilon}^{s_{1},s_{2}}}
\leq C\|u_{0}\|_{H^{s_{1},s_{2}}}+C\left\|u_{0}\right\|_{H^{s_{1},s_{2}}}=2C\left\|u_{0}\right\|_{H^{s_{1},s_{2}}}.\label{5.05}
\end{align}
Thus, $\Phi_{1}$ maps $B_{1}(0,2C\|u_{0}\|_{H^{s_{1},s_{2}}})$ into $B_{1}(0,2C\|u_{0}\|_{H^{s_{1},s_{2}}})$.
Combining  Lemmas 2.2, 4.1 with (\ref{5.04})-(\ref{5.05}), we have
\begin{align}
\left\|\Phi_{1}(u_{1})-\Phi_{1}(u_{2})\right\|_{X_{\frac{1}{2}+\epsilon}^{s_{1},s_{2}}}
\leq & C\left\|\frac{1}{2}\eta\left(\frac{t}{\tau}\right)
\int_{0}^{t}W(t-\tau)\partial_{x}(u_{1}^{2}-u_{2}^{2})d\tau\right\|_{X_{\frac{1}{2}+\epsilon}^{s_{1},s_{2}}}\nonumber\\
&\leq CT^{\epsilon}\left\|u_{1}-u_{2}\right\|_{X_{\frac{1}{2}+\epsilon}^{s_{1},s_{2}}}
\left[\left\|u_{1}\right\|_{X_{\frac{1}{2}+\epsilon}^{s_{1},s_{2}}}+\left\|u_{2}\right\|_{X_{\frac{1}{2}+\epsilon}^{s_{1},s_{2}}}\right]\nonumber\\
&\leq 4C^{2}T^{\epsilon}\left\|u_{0}\right\|_{H^{s_{1},s_{2}}}\left\|u_{1}-u_{2}\right\|_{X_{\frac{1}{2}+\epsilon}^{s_{1},s_{2}}}\nonumber\\
&\leq \frac{1}{2}\left\|u_{1}-u_{2}\right\|_{X_{\frac{1}{2}+\epsilon}^{s_{1},s_{2}}}.\label{5.06}
\end{align}
Thus, $\Phi_{1}$ is a contraction mapping  in the closed ball $B_{1}(0,2C\|u_{0}\|_{H^{s_{1},s_{2}}})$.
 Consequently, $u$ is the fixed point of $\Phi$ in the closed ball
$B_{1}(0,2C\|u_{0}\|_{H^{s_{1},s_{2}}})$. Then $v:=u|_{[0,T]}\in X_{\frac{1}{2}+\epsilon}^{s_{1},s_{2}}([0,T])$
is a solution to the Cauchy problem for (\ref{1.01}) with the initial data $u_{0}$ in the interval $[0,T]$.
 For the facts that  uniqueness of the solution and
the solution to the Cauchy problem for (\ref{1.01}) is continuous
with respect to the initial data, we refer the readers  to Theorems II, III of \cite{IMS}.

This completes the proof of Theorem 1.1.
\bigskip
\bigskip

\noindent {\large\bf 6.Proof of  Theorem  1.2}

\setcounter{equation}{0}

 \setcounter{Theorem}{0}

\setcounter{Lemma}{0}

\setcounter{section}{6}
We firstly  prove  Lemma 6.1  which is a variant of Theorem 1.1, then we apply Lemma 6.1 to prove Theorem 1.2.

\begin{Lemma}\label{Lemma6.1}
Let $s_{1}>-\frac{9}{8}$ and $R:=\frac{1}{8(C+1)^{3}}$, where  $C$  is the largest of
 those constants which appear in (\ref{2.07})-(\ref{2.08}),  (\ref{4.042}), (\ref{4.066}).
  Then, the Cauchy problem for (\ref{1.01}) is  locally well-posed for data satisfying
  \begin{align}
 \left\|I_{N}u_{0}\right\|_{L^{2}}\leq R.\label{6.01}
  \end{align}
Moreover, the solution to the Cauchy problem for (\ref{1.01}) exists on a time interval $[0,1]$.
\end{Lemma}
\noindent{\bf Proof.} We define $v:=I_{N}u$. Let $u$ be the solution to the Cauchy problem for (\ref{1.01}), then
$v$ is the solution to the following equations
\begin{align}
v_{t}+\partial_{x}^{5}v+\partial_{x}^{-1}
\partial_{y}^{2}v+\frac{1}{2}I_{N}\partial_{x}(I_{N}^{-1}v)^{2}=0.\label{6.02}
\end{align}
Thus,  $v$ satisfies the following equations
\begin{align}
v=W(t)v_{0}+\frac{1}{2}\int_{0}^{t}W(t-\tau)I_{N}\partial_{x}(I_{N}^{-1}v)^{2}.\label{6.03}
\end{align}
We define
\begin{align}
\Phi_{2}(v):=\psi(t)W(t)I_{N}u_{0}+\frac{1}{2}\psi(t)
\int_{0}^{t}W(t-\tau)I_{N}\partial_{x}(I_{N}^{-1}v)^{2}.\label{6.04}
\end{align}
Combining Lemma 2.2 with  4.3, we have
\begin{align}
\left\|\Phi_{2}(v)\right\|_{X_{\frac{1}{2}+\epsilon}^{0,0}}
\leq&
 \left\|\psi(t)W(t)I_{N}u_{0}\right\|_{X_{\frac{1}{2}+\epsilon}^{0,0}}
+C\left\|\psi(t)\int_{0}^{t}W(t-\tau)I_{N}\partial_{x}
(I_{N}^{-1}v)^{2}\right\|_{X_{\frac{1}{2}+\epsilon}^{0,0}}\nonumber\\
\leq & C\left\|I_{N}u_{0}\right\|_{L^{2}}+C\left\|I_{N}\partial_{x}
(I_{N}^{-1}v)^{2}\right\|_{X_{-\frac{1}{2}+2\epsilon}^{0,0}}\nonumber\\
\leq & C\left\|I_{N}u_{0}\right\|_{L^{2}}+C\left\|I_{N}\partial_{x}
(I_{N}^{-1}v)^{2}\right\|_{X_{-\frac{1}{2}+2\epsilon}^{0,0}}\nonumber\\
\leq & C\left\|I_{N}u_{0}\right\|_{L^{2}}+C\|v\|_{X_{\frac{1}{2}+\epsilon}^{0,0}}^{2}\nonumber\\
\leq & CR+C\|v\|_{X_{\frac{1}{2}+\epsilon}^{0,0}}^{2}.\label{6.05}
\end{align}
We define
\begin{align}
B_{2}(0,2CR):=\left\{v:\|v\|_{X_{\frac{1}{2}+\epsilon}^{0,0}}\leq
2CR\right\}.\label{6.06}
\end{align}
Combining   (\ref{6.05})-(\ref{6.06}) with the definition of  $R$,  we have
\begin{align}
&\left\|\Phi_{2}(v)\right\|_{X_{\frac{1}{2}+\epsilon}^{0,0}}
\leq CR+4C^{3}R^{2}
=2CR.\label{6.07}
\end{align}
Thus, $\Phi_{2}$ is a map from  $B_{2}(0,2CR)$ to
$B_{2}(0,2CR)$. We define
\begin{align}
v_{j}:=I_{N}u_{j}\>(j=1,2),\quad
w_{1}=I_{N}^{-1}v_{1}-I_{N}^{-1}v_{2},\quad
w_{2}:=I_{N}^{-1}v_{1}+I_{N}^{-1}v_{2}.\label{6.08}
\end{align}
Combining  Lemmas 2.2, 3.1, 3.2,  (\ref{6.05})-(\ref{6.06}) with the definition of $R$, we have
\begin{align}
\left\|\Phi_{2}(v_{1})-\Phi_{2}(v_{2})\right\|_{X_{\frac{1}{2}+\epsilon}^{0,0}}
\leq &
  C\left\|\psi(t)
\int_{0}^{t}W(t-\tau)\partial_{x}I_{N}\left[(I_{N}^{-1}v_{1})^{2}-(I_{N}^{-1}v_{2})^{2}\right]
d\tau\right\|_{X_{\frac{1}{2}+\epsilon}^{0,0}}\nonumber\\
&\leq C\left\|\partial_{x}I_{N}(w_{1}w_{2})\right\|_{X_{-\frac{1}{2}+2\epsilon}^{0,0}}\nonumber\\
&\leq C\|v_{1}-v_{2}\|_{X_{\frac{1}{2}+\epsilon}^{0,0}}
\left[\|v_{1}\|_{X_{\frac{1}{2}+\epsilon}^{0,0}}+\|v_{2}\|_{X_{\frac{1}{2}+\epsilon}^{0,0}}\right]\nonumber\\
&\leq 4C^{2}R
\|v_{1}-v_{2}\|_{X_{\frac{1}{2}+\epsilon}^{0,0}}\leq \frac{1}{2}\|v_{1}-v_{2}\|_{X_{\frac{1}{2}+\epsilon}^{0,0}}
.\label{6.09}
\end{align}
Thus, $\Phi_{2}$ is a contraction mapping from  $B_{2}(0,2CR)$ to $B_{2}(0,2CR)$.
 Consequently, $u$ is the fixed point of $\Phi_{2}$ in the closed ball
$B_{2}(0,2CR)$. Then $v:=Iu|_{[0,1]}\in X_{\frac{1}{2}+\epsilon}^{0,0}([0,1])$
is a solution to the Cauchy problem for (\ref{5.03}) with the initial data $I_{N}u_{0}$ on the interval $[0,1]$.
For the uniqueness of the solution and the fact that
the solution  is continuous
with respect to the initial data, we  refer the readers to Theorem II, III of \cite{IMS}.

This completes the proof of Lemma 6.1.

Now we apply the idea of  \cite{ILM-CPAA} and  Lemmas 2.7, 4.2, 6.1 to prove Theorem 1.2.

For $\lambda>0$, we define
\begin{align}
u_{\lambda}(x,y,t):=\lambda^{\frac{4}{5}}u
\left(\lambda^{\frac{1}{5}}x,\lambda^{\frac{3}{5}}y,\lambda t\right),
 u_{0\lambda}(x,y):=\lambda^{\frac{4}{5}}
u\left(\lambda^{\frac{1}{5}}x,\lambda^{\frac{3}{5}}y\right).\label{6.010}
\end{align}
Thus, $u_{\lambda}(x,y,t)\in X_{\frac{1}{2}+\epsilon}^{s_{1},0}([0,\frac{T}{\lambda}])$ is the solution to
\begin{align}
&\partial_{t}u_{\lambda}+\partial_{x}^{5}u_{\lambda}
+
\partial_{x}^{-1}\partial_{y}^{2}u_{\lambda}+u_{\lambda}\partial_{x}u_{\lambda}=0,\label{6.011}\\
&u_{\lambda}(x,y,0)=u_{0\lambda}(x,y),\label{6.012}
\end{align}
if and only if $u(x,y,t)\in X_{\frac{1}{2}+\epsilon}^{s,0}([0,T])$ is the solution to the
Cauchy problem for (\ref{1.01}) in $[0,T]$ with the initial data $u_{0}.$
By using a direct computation, for $\lambda \in (0,1),$  we have
\begin{align}
\|I_{N}u_{0\lambda}\|_{L^{2}}\leq CN^{-s}\lambda^{\frac{2}{5}+\frac{s}{5}}\|u_{0}\|_{H^{s,0}}.\label{6.013}
\end{align}
For $u_{0}\neq0$ and $u_{0}\in H^{s,0}(\R^{2})$, we choose $\lambda,N$ such that
\begin{align}
\|I_{N}u_{0\lambda}\|_{L^{2}}\leq CN^{-s}\lambda^{\frac{2}{5}+\frac{s}{5}}
\|u_{0}\|_{H^{s,0}}:=\frac{R}{4}.\label{6.014}
\end{align}
Then there exists $w_{3}$ which satisfies that  $\|w_{3}\|_{X_{\frac{1}{2}+\epsilon}^{s,0}}\leq  2CR$ such that $v:=w_{3}\mid_{[0,1]}$
is a solution to the Cauchy problem for (\ref{6.011}) with $u_{0\lambda}$.
Multiplying (\ref{6.011}) by $2I_{N}u_{\lambda}$ and integrating with respect to $x,y$
 yield
\begin{align}
\frac{d}{dt}\int_{\SR^{2}}(I_{N}u)^{2}dxdy+\int_{\SR^{2}}I_{N}u\partial_{x}I_{N}
\left[(u)^{2}\right]dxdy=0.\label{6.015}
\end{align}
Inserting
\begin{align*}
\int_{\SR^{2}}I_{N}u\partial_{x}\left[(I_{N}u)^{2}\right]dxdy=0
\end{align*}
into (\ref{6.015}) yields
\begin{align}
\frac{d}{dt}\int_{\SR^{2}}(I_{N}u)^{2}dxdy=-
\int_{\SR^{2}}I_{N}u\partial_{x}\left[I_{N}\left((u)^{2}\right)-(I_{N}u)^{2}\right]dxdy.\label{6.016}
\end{align}
Combining  (\ref{6.016}) with Lemmas 2.6, 4.2,  we have
\begin{align}
&\quad
\int_{\SR^{2}}(I_{N}u(x,y,1))^{2}dxdy-\int_{\SR^{2}}(I_{N}u_{0\lambda})^{2}dxdy\nonumber\\   &
=-\int_{0}^{1}\int_{\SR^{2}}I_{N}u_{\lambda}\partial_{x}\left[I_{N}\left((u_{\lambda})^{2}
\right)-(I_{N}u_{\lambda})^{2}\right]dxdydt\nonumber\\
&=-\int_{\SR}\int_{\SR^{2}}\left(\chi_{[0,1]}(t)I_{N}u_{\lambda}\right)
\left(\chi_{[0,1]}(t)\partial_{x}\left[I_{N}\left((u_{\lambda})^{2}\right)-(I_{N}u_{\lambda})^{2}\right]\right)dxdydt
\nonumber\\
&\leq C\left\|\chi_{[0,1]}(t)I_{N}u_{\lambda}\right\|_{X_{\frac{1}{2}-\epsilon}^{0,0}}
\left\|\chi_{[0,1]}(t)\partial_{x}\left[I_{N}\left((u_{\lambda})^{2}\right)-(I_{N}u_{\lambda})^{2}\right]
\right\|_{X_{-\frac{1}{2}+\epsilon}^{0,0}}\nonumber\\
&\leq C\left\|I_{N}u_{\lambda}\right\|_{X_{\frac{1-\epsilon}{2}}^{0,0}}
\left\|\partial_{x}\left[I_{N}\left((u_{\lambda})^{2}\right)-(I_{N}u_{\lambda})^{2}\right]\right\|_{X_{-\frac{1}{2}+2\epsilon}^{0,0}}\nonumber\\
&\leq CN^{-2+10\epsilon}\|I_{N}u_{\lambda}\|_{X_{\frac{1}{2}+\epsilon}^{0,0}}^{3}.\label{6.017}
\end{align}
From (\ref{6.014}) and  (\ref{6.015}) and the definition of $R$, we have
\begin{align}
\int_{\SR^{2}}(I_{N}u(x,y,1))^{2}dxdy
&\leq \frac{R^{2}}{16}+ CN^{-2+10\epsilon}
\|I_{N}u_{\lambda}\|_{X_{\frac{1}{2}+\epsilon}^{0,0}}^{3}\nonumber\\
&\leq \frac{R^{2}}{16}+8C^{4}N^{-2+10\epsilon}R^{3}\leq \frac{R^{2}}{16}+CN^{-2+10\epsilon}.\label{6.018}
\end{align}
Let $N$ be sufficiently large such that such that $8C^{4}N^{-2+10\epsilon}R^{3}\leq \frac{3}{4}R^{2},$
then
\begin{align}
\left[\int_{\SR^{2}}(I_{N}u(x,y,1))^{2}dxdy\right]^{\frac12}\leq R.\label{6.019}
\end{align}
We consider $I_{N}u(x,y,1)$ as the initial data and repeat the above argument, from Lemma 6.1, we obtain that (\ref{6.011})-(\ref{6.012}) possess
a solution in $\R^{2}\times [1,2]$. In this way, we can extend the solution to (\ref{6.011})-(\ref{6.012}) to the time interval
$[0,2].$ The above argument can be repeated $L$ steps, where $L$ is the maximal positive integer  such that
\begin{align}
CN^{-2+10\epsilon}L\leq \frac{3}{4}R^{2}.\label{6.020}
\end{align}
More precisely, the solution to (\ref{6.011})-(\ref{6.012}) can be extended to the time interval $[0,L]$. Thus, we can prove that (\ref{6.011})-(\ref{6.012})
are globally well-posed in $[0,\frac{T}{\lambda}]$ if
\begin{align}
L\geq \frac{T}{\lambda}.\label{6.021}
\end{align}
From (\ref{6.020}), we know that
\begin{align}
L\sim N^{2-10\epsilon}.\label{6.022}
\end{align}
We know  that (\ref{6.021}) is valid provided that the following inequality is valid
\begin{align}
CN^{2-10\epsilon}\geq \frac{T}{\lambda}\sim CTN^{\frac{-5s}{2+s}}.\label{6.023}
\end{align}
In fact,  (\ref{6.023}) is valid if
\begin{align}
N^{2}> N^{\frac{-5s}{2+s}}\label{6.024}
\end{align}
which is equivalent to $-\frac{4}{7}<s<0.$

This completes the proof of Theorem 1.2.

\bigskip
\bigskip

\leftline{\large \bf Acknowledgments}

\bigskip

\noindent

 This work is supported by the Natural Science Foundation of China
 under grant numbers 11571118 and 11771127. The first author is also
 supported by
 the Fundamental Research Funds for the
Central Universities of China under the grant
number 2017ZD094, while second author is also
 supported by the Young core Teachers program of Henan Normal University
 under grant number 15A110033.

\end{document}